\documentclass[11pt, reqno]{article}
\usepackage{amsfonts}
\usepackage{amsmath}
\usepackage{amsthm}
\usepackage{amssymb}
\usepackage{fancyhdr}
\usepackage{fancybox}
\usepackage{mathrsfs}
\usepackage{enumerate}
\usepackage{amsfonts}
\usepackage{fancyhdr}
\usepackage{fancybox}
\usepackage{mathrsfs}

\def\phi{\varphi }
\def\epsilon{\varepsilon}
\def\cal{\mathcal }
\swapnumbers

\theoremstyle{plain}
\newtheorem{theorem}{Theorem}[section]
\newtheorem{lemma}[theorem]{Lemma}
\newtheorem{proposition}[theorem]{Proposition}
\newtheorem{corollary}[theorem]{Corollary}
\theoremstyle{definition}
\newtheorem{definition}[theorem]{Definition}
\newtheorem{example}[theorem]{Example}
\newtheorem{examples}[theorem]{Examples}
\newtheorem{remark}[theorem]{Remark}
\newtheorem{remarks}[theorem]{Remarks}
\newtheorem{facts}[theorem]{Facts}
\newtheorem{deffacts}[theorem]{Definitions and facts}
\newtheorem{orbit-schemes}[theorem]{Orbit schemes}
\newtheorem{orbit-schemes-mod}[theorem]{Orbit schemes modified by positive semicharacters}
\newtheorem{directp}[theorem]{Direct products}
\newtheorem{joinconst}[theorem]{Joins of CAS}
\newtheorem{cas-dist-gr}[theorem]{CAS associated with infinite distance-transitive graphs}
\newtheorem{exposmult}[theorem]{Examples of positive multiplicative pairs}
\newtheorem{random-walk-hy}[theorem]{Random walks on hypergoups}
\newtheorem{random-walk-CAS}[theorem]{Random walks on continuous association schemes}
\numberwithin{equation}{section}

\begin{document}

\title{Continuous Association Schemes and  Hypergroups}

\author{Michael Voit\\
Fakult\"at Mathematik, Technische Universit\"at Dortmund\\
          Vogelpothsweg 87,
          D-44221 Dortmund, Germany\\
e-mail:  michael.voit@math.tu-dortmund.de}

\date{\today}

\maketitle



\begin{abstract}
Classical finite association schemes lead to a finite-dimensional algebras which are generated by 
finitely many stochastic matrices. Moreover, there  exist associated finite hypergroups.
The notion of classical discrete  association schemes can be easily extended to the possibly infinite case. 
Moreover, the notion  of association schemes can be relaxed slightly by using suitably deformed families of stochastic
 matrices by skipping the integrality conditions.
 This leads to larger class of examples which are again associated to discrete hypergroups.

In this paper we propose a topological generalization of the notion of association schemes by using a
locally compact basis space $X$ and a family of  Markov-kernels on $X$ indexed
 by a further locally compact space $D$ where the supports of the associated probability measures
 satisfy some partition property.
These objects, called  continuous association schemes, will be 
 related to hypergroup structures on $D$. 
We study some basic results for this new notion and present several classes of examples.
It turns out that for a given commutative hypergroup the existence of an associated
  continuous association scheme implies that the hypergroup has many features of a double coset hypergroup.
We in particular show
that commutative hypergroups, which are associated with commutative 
continuous association schemes, carry dual positive product formulas for the
characters. On the other hand, we prove some rigidity results in particular in the compact case
which say that for given spaces $X,D$ there are only a few
continuous association schemes.
\end{abstract}

\maketitle

\section{Introduction}

In this paper we study a topological generalization of the notion of classical finite
association schemes by using the notion of hypergroups in the sense of Dunkl, Jewett, and Spector.
To explain this, let us start with the  notion of a finite
  association scheme which is common in algebraic combinatorics; see e.g. the monographs  \cite{BI}, \cite{B}.

\begin{definition}
 Let $X, D$ be finite nonempty
sets and $(R_i)_{i\in D}$ a disjoint
  partition of $X\times X$ with $R_i\ne\emptyset$ ($i\in D$) and with:
\begin{enumerate}
\item[\rm{(1)}] There exists $e\in D$ with $R_e=\{(x,x): x\in X\}$.
\item[\rm{(2)}] There exists an involution $i\mapsto \bar i$ on $D$ such that
  for $i\in D$,
$R_{\bar i}=\{(y,x):\> (x,y)\in R_i\}$.
\item[\rm{(3)}] For all $i,j,k\in D$ and $(x,y)\in R_k$, the number
$$p_{i,j}^k:=|\{z\in X:\> (x,z)\in R_i\>\> {\rm and }\>\> (z,y)\in R_j\}|$$
is independent of $(x,y)\in R_k$.
\end{enumerate}
Then $\Lambda:=(X,D, (R_i)_{i\in D})$ is called a finite association scheme with
intersection numbers $(p_{i,j}^k)_{i,j,k\in D}$ and identity $e$. 
\end{definition}

Now let $\Lambda:=(X,D, (R_i)_{i\in D})$ be  a finite association scheme.
Form the 
 adjacency 
matrices $A_i\in \mathbb R^{X\times X}$ ($i\in D$) with
$$(A_i)_{x,y}:=
\left\{ \begin{array}{cc}
\displaystyle 1 &
 {\rm if}\>\> (x,y)\in R_i
\\     0&  {\rm otherwise}\\     \end{array} \right. 
 \quad\quad\quad \quad(i\in D,\>
x,y\in X).$$
Then $A_e$ is the identity matrix $I_X$, and the transposed matrices satisfy
$A_i^T = A_{\bar i}$ for $i\in D$.
Moreover, for $i,j\in D$, we have for the usual matrix product 
 $A_iA_j= \sum_{k\in D} p_{i,j}^k A_k$.

Define the valencies
\begin{equation}
\omega_i:=p_{i,\bar i}^e=|\{z\in X:\> (x,z)\in R_i\}|\in\mathbb N
\end{equation}
 of $R_i$ (or $i\in D$), where these numbers are independent of $x\in X$.
Then the  renormalized  adjacency matrices
  $S_i:= \frac{1}{\omega_i}A_i\in \mathbb R^{X\times X}$ are stochastic, i.e.,
  all  row sums are equal to 1. Moreover, the products
\begin{equation}\label{product-rule-si}
S_iS_j= \sum_{k\in D} \frac{\omega_k}{\omega_i\omega_j} p_{i,j}^k S_k  \quad(i,j\in D)
\end{equation}
are convex combinations of the $S_i$, and the linear span of the $S_i$ is a finite dimensional 
algebra. This algebra is isomorphic with the algebra of measures of some finite hypergroup structure on $D$ in the sense
of Dunkl, Jewett, and Spector, where the $S_i$ are identified with the point measures $\delta_i$ of $i\in D$.
For this we recapitulate the definition of a finite hypergroup;
 see \cite{BH}, \cite{Du}, \cite{J}, and in the finite case, \cite{W1}, \cite{W2}. 
We point out that we  do not use another definition of hypergroups
 where products of sets are considered, and which runs under the subject classification 20N20.

\begin{definition}\label{def-finite-hy}
A finite hypergroup $(D,*)$ is a  finite non-empty set  $D$ with an
 associative, bilinear, probability-preserving multiplication $*$ (called convolution) 
on the vector space $M_b(D)$
  of all complex  measures on $D$ with the following properties:
\begin{enumerate}
\item[\rm{(1)}] There exists a neutral element $e\in D$ with $\delta_x*\delta_e=
\delta_e*\delta_x=\delta_x$ for  $x\in D$. 
\item[\rm{(3)}] There exists an involution $x\mapsto\bar x$ on $D$
  such that for all $x,y\in D$, $e\in supp\> (\delta_x*\delta_y)$ 
  if and only if $y=\bar x$. 
\item[\rm{(4)}] If for $\mu\in M_b(D)$, $\mu^-$ is the
  image of $\mu$ under the involution, then
  $(\delta_x*\delta_y)^-= \delta_{\bar y}*\delta_{\bar x}$ for all $x,y\in D$.
\end{enumerate}
\end{definition}

Now let  $\Lambda:=(X, D, (R_i)_{i\in D})$ be a finite association
  scheme as above. It is easy to see that then the unique bilinear extension
of the convolution
$$\delta_i*\delta_j:=\sum_{k\in D}\frac{\omega_{ k}}{\omega_{ i}\omega_{ j}} 
\cdot p_{i,j}^k \delta_k  \quad\quad(i,j\in D)$$
of point measures leads to a finite hypergroup $(D,*)$, 
the so-called hypergroup associated with $\Lambda$.

Classical examples of finite association schemes and hypergroups appear from groups:

\begin{example}\label{ex-finite-double-coset}
Let $H$ be a subgroup of a finite group $G$ with identity $e$.
Consider the set  $X:=G/H:=\{gH:\> g\in G\}$ of cosets as well as the
set $D:=G//H:=\{HgH:\> g\in G\}$ of double cosets. It can be easily cheked and is well-known that
 the partition
$R_{HgH}:=\{(xH, yH)\in X\times X: \> Hx^{-1}yH=HgH\}$ ($HgH\in D$) of $X\times X$ leads
 to a finite association scheme with identity $HeH$ and involution $HgH\mapsto Hg^{-1}H$.
The associated hypergroup is the so-called double coset hypergroup $(D,*)$ with
the double coset convolution
$$\delta_{HxH}*\delta_{HyH}= \frac{1}{|H|} \sum_{h\in H}\delta_{HxhyH} \quad (x,y\in G).$$
The associated convolution algebra $M_b(D)$ is often also called a Hecke-algebra. 
\end{example}

Typical commutative examples for \ref{ex-finite-double-coset} appear,
 if one considers so-called distance-transitive graphs
$X$ on which the group $G$ of all graph automorphisms acts
 where $H$ is the fix group of some  vertex. We then have $G/H\equiv X$ in a canonical way,
 and $G//H$ can be identified with $\{0,1,\ldots,N\}$ with the diameter $N$ of $X$.
We do not give further details here and refer
 to \cite{BI}. The set of  distance-transitive graphs
is a proper subset of the set of  distance-regular graphs where again  canonical associated 
 commutative association schemes and  commutative hypergroups exist, and where 
the construction of  \ref{ex-finite-double-coset} via groups usually is not longer available.
We here skip details of the theory of  distance-regular graphs and refer
 to the  monographs \cite{BI}, \cite{BCN} and the recent survey \cite{DKT}.
\medskip

Let us now consider generalizations of association schemes.
 We first skip the condition that $X$ and $D$ are finite,
 where we keep a finiteness condition for the partition; see 
Definition \ref{def-discrete-asso} below. It turns out that then most statements for  
 finite association schemes remain valid. In particular, there exist associated discrete 
hypergroups $(D,*)$  as  above. Typical examples appear when we consider 
 totally disconnected, locally compact groups $G$ with a compact, open subgroup $H$.
Then the spaces $X:=G/H$ and $D:=G//H$ as above are discrete with respect to the quotient topology, and
$X,D$ and a partition as in \ref{ex-finite-double-coset} lead to a possibly infinite
association scheme. The associated hypergroup  $(D,*)$ is then just the double coset hypergroup
 studied in hypergroup theory; see \cite{J}. The details are worked out in \cite{Voit16}.
 Typical examples for this are the infinite association schemes associated with homogeneous trees and, slightly more general,  with
infinite distance-transitive graphs; see  \cite{Voit16}.
We also mention that there exists  examples of higher rank with
 $X$ as sets of vertices of affine  buildings; see e.g.~\cite{APM} and references there.

We next consider a further discrete extension from \cite{Voit16}. Fix some (possibly infinite) 
association scheme 
with an associated  partition and the associated stochastic, renormalized adjacency matrices $S_i$; 
  assume that we in addition have a further  algebra of 
matrices generated by ``deformed'' stochastic matrices  $\tilde S_i\in \mathbb R^{X\times X}$ for $i\in D$
where  any entry of any $S_i$ is positive if and only if so is the corresponding entry of $\tilde S_i$.
We add some further technical axioms like $\tilde S_e=S_e$ and that there is a measure on $X$ which
 replaces the counting measure of an association scheme and which satisfies some adjoint relation; 
see Definition  \ref{def-gen-discrete-asso}.
It turns out that these so-called generalized association schemes with the matrices $\tilde S_i$ instead of the $S_i$
 also admit associated hypergroups $(D,*)$ 
as  above.

In this paper we  use this  notion of generalized association schemes from \cite{Voit16} and present
a topological extension  in Definition \ref{central-def} by using families of 
Markov kernels on locally compact spaces $X$ which are indexed by some locally compact space $D$
 instead of stochastic matrices as before. We require that the supports of the measures associated with these kernels
admit partition properties similar to those of  association schemes, and we require that  the kernels generate an 
algebra  such that again the product linearizations of the kernels fit to some 
 hypergroup structure $(D,*)$. In addition, some  natural topological conditions are added.
We point out that we here require from the beginning that there exists an associated  hypergroup structure $(D,*)$
(different from the discrete case). We have done this as we otherwise would run into technical topological problems
  (which we want to avoid in this paper), 
and as for all known examples this hypergroup 
property is available from the beginning. This is in particular the case 
for  standard classes of examples of such continuous association schemes (CAS for short). 
Here is a short incomplete list of examples:
\begin{enumerate}
\item[\rm{(1)}] If $H$ is a compact subgroup of a locally compact group $G$,  then $X:=G/H$ and $D:=G//H$ lead to canonical
CAS associated with groups analog to the finite case or the case where $H\subset G$ is compact and open.
\item[\rm{(2)}] All (unimodular) association schemes and all generalized association schemes as above are CAS.
\item[\rm{(3)}] If a non-compact commutative  CAS is given, then it often can be deformed via so-called pairs $(\alpha,\phi)$
of positive multiplicative functions on $D$ and $X$; see Section 8 for details. This construction often  leads
 to plenty of interesting families of deformed CAS
with deformed Markov kernels, where the spaces $X,D$ remain unchanged.
On the other hand, in the compact and in particular finite case, the situation is much more
 rigid. It turns out that 
for given compact spaces $X,D$, there is at most one associated CAS structure; see Corollary \ref{equal-compact}.
Moreover, each finite CAS is  automatically an association scheme, i.e., there is in fact no freedom  in the choice
of the stochastic matrices  $\tilde S_i$ of a generalized association scheme  in the finite case.
This difference between the compact and non-compact case is  remarkable.
\item[\rm{(4)}] Besides the examples indicated above we point out that there are several further standard constructions 
to get new CAS from given ones; see Section 11.
\end{enumerate}
%\medskip

This paper is organized as follows. 
In Section 2 we recapitulate some facts about  hypergroups in the sense of Dunkl, Jewett and Spector
with a focus on the commutative case;
the main references are the monograph of Bloom and Heyer \cite{BH} and  Jewett \cite{J}.
Some technical details of Section 2 may be skipped at a first reading.
In Section 3 we recapitulate some notations and facts on possibly infinite classical 
association schemes and their discrete generalizations 
 mentioned above.
This  discrete generalization  motivates the definition of
continuous association schemes (CAS for short) 
on the basis of Markov kernels and associated transition operators in Section 4.
The central Section 4 contains the discussion of basic properties and some  natural classes of examples. 
In Section 5 we add some further axioms to the basic definition,
 called translation properties (T1) and (T2), which are needed to  get stronger 
interrelation between the analysis on $X$ and $D$. It turns out that all
 compact CAS as well as all CAS associated with 
groups and all classical discrete  association schemes have these properties.
As a byproduct we obtain
some rigidity result, e.g., that all 
finite CAS are in fact association schemes.

Section 6 is then devoted to positive definite functions on $D$ and $X$ for commutative CAS.
We in particular obtain that each commutative hypergroup $(D,*)$
 which is associated with some CAS with property (T2) admits a dual positive 
convolution on the support of the Plancherel measure of $(D,*)$; see Theorem 
\ref{main-cont}.  This central result will be improved in Section 7 where we consider
 two possibly different commutative  CAS structures with the same basic spaces $X,D$
 where we assume that one of them has property (T2) and where the schemes are related in some way.
The central positive definiteness result in Theorem \ref{cont-main-g} 
will also lead to further rigidity results.

The sections 8-11 are  mainly 
 devoted to examples of CAS and  construction principles of  examples
 beyond the group cases and discrete association schemes.
We start in Section 8 with
 nontrivial functions $\phi$ on $X$ which are eigenfunctions under
 all transition operators of the given commutative CAS. It turns out that these $\phi$ are always
 related with multiplicative functions $\alpha$ of the hypergroup $(D,*)$. The interrelations between  $\phi$ and
$\alpha$ will
lead to further results regarding the properties (T1) and (T2) in the commutative case.
Moreover, if $\phi$ and
$\alpha$ are in addition positive, we shall  construct
 a deformed CAS  with the same spaces $X,D$ but deformed Markov kernels. 
On the level of hypergroups this deformation is just the known deformation of a hypergroup by
 a positive semicharachter in \cite{V1}, \cite{BH}.
In the case of commutative CAS associated with
 non-compact symmetric spaces $X$, the eigenfunctions $\phi$ 
are closely related with the joint eigenfunctions of the invariant differential operators on $X$ which 
are completely classified; see \cite{Hel}, \cite{Kash}. In Section 9 we shall mainly study the deformation
of commutative CAS which appear via orbits when some compact group acts continuously on 
some locally compact abelian group.

Section 10 is  devoted to the deformation of a concrete class of examples,
 namely of the infinite association schemes associated with infinite distance-transitive graphs.
For this recapitulate that the set of these graphs extend the class of all homogeneous trees
 only slightly and is parametrized by two parameters. We show how
 boundary points of these graphs lead to  deformations.
In Section 11 we present several further standard  constructions which lead from given CAS to
 new ones. Typical examples are direct products and  joins, 
which are well-known in the theory of hypergroups by \cite{J}, \cite{BH}.

Section 12 contains an introduction into random walks on $X$ for a CAS $(X,D,K)$; we in particular show that the canonical
projections of these random walks to $D$ are random walks on the hypergroup $(D,*)$. This observation may be used 
to transfer limit theorems for 
 random walks  on $(D,*)$ like (strong) laws of large numbers and central limit theorems
 (see Ch.~7 of  \cite{BH} and references there) to random walks on $X$ in future.
This seems to be interesting 
 in particular for examples which appear as deformations of 
group CAS, as here random walks  on $X$ may be seen as ``radial random walks with aditional drift'' on the
homogeneous space $X$. 
Finally, Section 13 contais a short list of central open problems for CAS.

\section{Hypergroups}

In this section we recapitulate some facts on  hypergroups with an focus on the commutative case 
mainly from
\cite{Du}, \cite{J},  \cite{BH}. Only some results in the end of this section are new.

Hypergroups  form an extension  of locally compact   groups. For
 this,  remember that the group multiplication
 on a locally compact   group $(G, \cdot)$ leads to the convolution
$\delta_x*\delta_y=\delta_{xy}$ ($x,y\in G$) of point measures. 
Bilinear, weakly continuous extension  of this convolution together with the canonical involution
with $\delta_x\mapsto \delta_{x^{-1}}$  then lead to a  
 Banach-$*$-algebra structure on the Banach space $M_b(G)$ of all
signed bounded regular Borel measures with the total variation norm
$\|.\|_{TV}$ as norm.

In the case of hypergroups we usually do not have an algebraic operation on the basis space,
and we only require 
 a convolution  $*$ for bounded complex measures which admits 
most properties of a group
 convolution:

\begin{definition}\label{def-hy-gen}
A hypergroup $(D,*)$ is a locally compact Hausdorff space $D$ with a
  weakly continuous, associative, bilinear convolution $*$ on the Banach space $M_b(D)$
  of all bounded, complex regular Borel measures with the following properties:
\begin{enumerate}
\item[\rm{(1)}] For all $x,y\in D$,  $\delta_x*\delta_y$  is a compactly
  supported probability measure on $D$ such that the support 
  $supp\>(\delta_x*\delta_y)$ depends continuously on $x,y$ w.r.t.~the
  so-called Michael topology on the space of all compacta in $X$
 (see \cite{J} for details).  
\item[\rm{(2)}] There exists a neutral element $e\in D$ with $\delta_x*\delta_e=
\delta_e*\delta_x=\delta_x$ for  $x\in D$. 
\item[\rm{(3)}] There exists a continuous involution $x\mapsto\bar x$ on $D$
  such that for all $x,y\in D$, $e\in supp\> (\delta_x*\delta_y)$ holds
  if and only if $y=\bar x$. 
\item[\rm{(4)}] If for $\mu\in M_b(D)$, $\mu^-$ denotes the
  image of $\mu$ under the involution, then
  $(\delta_x*\delta_y)^-= \delta_{\bar y}*\delta_{\bar x}$ for all $x,y\in D$.
\end{enumerate}
A hypergroup is called commutative if the convolution $*$ is commutative. 
It is called symmetric if the involution is the identity. 
\end{definition}

If $D$ is finite, then Definition \ref{def-hy-gen} agrees with that of the introduction.

\begin{remarks}
\begin{enumerate}
\item[\rm{(1)}]  The identity $e$ and the involution $.^-$ above
are  unique.
\item[\rm{(2)}] Each symmetric hypergroup is  commutative.
\item[\rm{(3)}] For each hypergroup $(D,*)$, $(M_b(D),*)$ is a
  Banach-$*$-algebra with the involution $\mu\mapsto\mu^*$ with
$\mu^*(A):=\overline{\mu(A^-)}$   for Borel sets $A\subset D$.
\item[\rm{(4)}] For a second countable locally compact space $D$, the Michael
  topology agrees with the well-known Hausdorff topology; see \cite{KS}.
\end{enumerate}
\end{remarks}

The most prominent examples of hypergroups are  
double coset hypergroups $G//H:=\{HgH:\> g\in G\}$
for compact subgroups $H$  of 
locally compact groups $G$. This extends the discussion in the introduction:

\begin{example}\label{Doublecoset} Let $H$ be  a compact
 subgroup of a locally compact group
  $G$ with identity $e$ and 
 unique normalized Haar measure $\omega_H\in M^1(H)\subset M^1(G)$.
Then the space
$$M_b(G||H):=\{\mu\in M_b(G):\> \mu=\omega_H*\mu*\omega_H\}$$
of all $H$-biinvariant measures in $M_b(G)$ is a Banach-$*$-subalgebra of
$M_b(G)$. With the quotient topology, 
 $G//H$  is  locally compact, and the
canonical
projection $p_{G//H}:G\to G//H$ is continuous, proper and open. 
Now consider the push forward 
  $\tilde p_{G//H}:M_b(G)\to M_b(G//H)$  with 
$\tilde p_{G//H}(\mu)(A)=\mu(p_{G//H}^{-1}(A))$ for $\mu\in M_b(G)$
and Borel sets $A\subset G//H$. Then
 $\tilde p_{G//H}$ is  an isometric isomorphism between the Banach
 spaces  $M_b(G||H)$ and $M_b(G//H)$ w.r.t.~the total variation
 norms, and the  transfer of the  convolution on
 $M_b(G||H)$  to   $M_b(G//H)$
 leads to a hypergroup $(G//H, *)$ with  identity $HeH$ 
 and  involution $HgH\mapsto Hg^{-1}H$, cf.~\cite{J}.

The pair $(G,H)$ is called a Gelfand pair if the double coset hypergroup is commutative. 
For the theory of Gelfand pairs we refer to \cite{Di} and \cite{F}.
\end{example}

The notion of  Haar
measures on hypergroups generalizes  that on  groups:

\begin{definition} Let $(D,*)$ be a hypergroup,
 $x,y\in D$, and $f\in C_c(D)$ a continuous function with compact support.
We write $_xf(y):=f(x*y):=\int_K f\>
d(\delta_x*\delta_y)$ and $f_x(y):=f(y*x)$  where, 
by the hypergroup axioms, $f_x,\>_xf  \in C_c(D)$ holds by \cite{J}.

A non-trivial positive  Radon measure $\omega\in M^+(D)$ 
is  a left or right Haar measure if 
$${\int_D} {_xf}\> d\omega=
\int_D f\> d\omega     \quad\text{or}\quad
{\int_D} {f_x}\> d\omega=
\int_D f\> d\omega \quad\quad (f\in C_c(D), \> x\in D)$$ 
  respectively. $\omega$ is 
 called a Haar measure if
it is a left and  right Haar measure.
If $(D,*)$  admits a  Haar measure, then it is called unimodular. 
\end{definition} 

The uniqueness of left and right  Haar measures
 and their existence for
particular classes  are known for a long time by 
Dunkl, Jewett, and Spector;
see \cite{BH} for details. The general existence was settled only
 recently by Chapovsky \cite{Ch}:

\begin{theorem} Each hypergroup admits a left and a right Haar measure.
Both  are unique up to normalization.
\end{theorem}

\begin{examples}\label{example-haar}
\begin{enumerate} 
\item[\rm{(1)}] Let  $(D,*)$ be a discrete hypergroup. Then, by \cite{J},
 left and right Haar measures
 are given by 
$$\omega_l(\{x\})= \frac{1}{(\delta_{\bar x}*\delta_x)(\{e\})}, \quad
\omega_r(\{x\})= \frac{1}{(\delta_{ x}*\delta_{\bar x})(\{e\})} \quad
\quad(x\in D).$$
Notice that discrete hypergroups are not necessarily unimodular; see e.g. \cite{KW}
 for examples of double coset hypergroups.
\item[\rm{(2)}] If $(G//H,*)$ is a double coset hypergroup and  $\omega_G$
 a left Haar measure of $G$, then its canonical projection to
 $G//H$ is a left Haar measure of $(G//H,*)$. 
 \end{enumerate}
\end{examples}

We next recapitulate some facts on Fourier analysis on commutative hypergroups
from \cite{BH}, \cite{J}.
For the rest of Section 2 let $(D,*)$ be a commutative hypergroup with Haar
measure $\omega$. For $p\ge1$ consider the $L^p$-spaces
$L^p(D):=L^p(D,\omega)$. Moreover $C_b(D)$ and $C_0(D)$ are the Banach spaces
of all bounded continuous functions on $D$ and 
those which vanish at infinity respectively. For a function $f:D\to\mathbb C$ and $x\in X$ we put $f^-(x):=f(\bar x)$ and $f^*(x):=\overline{f(\bar x)}$.

\begin{deffacts}\label{deffacts}
\begin{enumerate} 
\item[\rm{(1)}] The spaces of all (bounded) non-trivial multiplicative continuous functions on $(D,*)$ are
$$\chi(D,*):=\{\alpha\in C(D):\>\> \alpha\not\equiv 0,\>\>
\alpha(x* y)=\alpha(x)\cdot \alpha(y) \>\>\text{for all}\>\>
x,y\in D\}$$
and $\chi_b(D,*):=\chi(D,*)\cap C_b(D)$.
$$\hat D:=(D,*)^\wedge:= \{\alpha\in \chi_b(D,*):\>\> \alpha(\bar x)=\overline{\alpha(x)} 
 \>\>\text{for all}\>\>
x\in D\}$$
is the dual space of $(D,*)$. Its elements are called characters.

All spaces will be equipped with the topology of compact-uniform
convergence. $\chi_b(D,*)$ and $\hat D$ are then locally compact.
 
If $D$ is discrete, then $\hat D$ is compact, and if  $D$ is compact, then $\hat D$ is discrete.

 All
characters $\alpha\in\hat D$ satisfy $\|\alpha\|_\infty=1$ and $\alpha(e)=1$.
\item[\rm{(2)}] For $f\in L^1(D)$ and $\mu\in M_b(D)$, their Fourier(-Stieltjes) transforms
  are defined by
$$\hat f(\alpha):=\int_D f(x)\overline{\alpha(x)}\> d\omega(x),\quad
\hat \mu(\alpha):=\int_D \overline{\alpha(x)}\> d\mu(x) 
\quad (\alpha\in\hat D).$$
We have $\hat f\in C_0(\hat D)$, $\hat\mu\in C_b(\hat D)$ and 
$\|\hat f\|_\infty\le\|f\|_1$, $\|\hat \mu\|_\infty\le\|\mu\|_{TV}$.
\item[\rm{(3)}] There exists a unique positive measure $\pi\in M^+(\hat D)$
  such that the  Fourier transform
$.^\wedge: L^1(D)\cap L^2(D) \to C_0(\hat D)\cap L^2(\hat D,\pi)$
is an isometry. $\pi$ is called the Plancherel measure on $\hat D$.

Notice that, different from locally compact abelian groups, the support
$S:=supp\> \pi$ may be a proper closed subset of $\hat D$.
Quite often, we even have ${\bf 1}\not\in S$.
\item[\rm{(4)}] For $f\in L^1(\hat D,\pi)$, $\mu\in M_b(\hat D)$, their
  inverse
Fourier transforms are given by
$$\check f(x):=\int_S f(\alpha) {\alpha(x)}\> d\pi(\alpha),\quad
\check \mu(x):=\int_{\hat D} {\alpha(x)}\> d\mu(\alpha) 
\quad (x\in D)$$
with $\check f\in C_0(D)$, $\check\mu\in C_b(D)$
and 
$\|\check f\|_\infty\le\|f\|_1$, $\|\check \mu\|_\infty\le\|\mu\|_{TV}$.
\item[\rm{(5)}] $f\in C_b(D)$ is called positive definite 
on the hypergroup $D$ if for all $n\in\mathbb N$, $x_1,\ldots, x_n\in D$ and
  $c_1,\ldots,c_n\in\mathbb C$,
$\sum_{k,l=1}^n c_k \bar c_l \cdot f(x_k*\bar x_l)\ge0.$
Obviously, all characters $\alpha\in\hat D$ are positive definite.
 \end{enumerate}\end{deffacts}

We collect further results:

\begin{facts}\label{facts-hy}
\begin{enumerate} 
\item[\rm{(1)}] (Theorem of Bochner, \cite{J}) A function $f\in C_b(D)$ is positive
  definite if and only if $f=\check \mu$ for some $\mu\in M_b^+(\hat D)$.
In this case, $\mu$ is a probability measure if and only if
$\check\mu(e)=1$.
\item[\rm{(2)}] For $f,g\in L^2(D)$, the convolution product 
$f*g(x):=\int f(x*\bar  y)g(y)\> d\omega(y)$ ($x\in D$)
 satisfies $f*g\in C_0(D)$. Moreover,
for $f\in L^2(D)$,  $f^*(x)=\overline{f(\bar x)}$ satisfies
  $f^*\in L^2(D)$, and 
 $f*f^*\in C_0(D)$ is positive definite; see \cite{J}, \cite{BH}. 
\item[\rm{(3)}]  (Refining of the Theorem of Bochner, \cite{V2})
For a positive definite function  $f\in C_b(D)$ with $f=\check \mu$ for some $\mu\in M_b^+(\hat D)$,
 the following statements are equivalent
\begin{enumerate} 
\item[\rm{(i)}] $supp\>\mu\subset S$;
\item[\rm{(ii)}]  $f$ is the compact-uniform limit of positive definite functions in $C_c(D)$;
\item[\rm{(iii)}]  $f$ is the compact-uniform limit of  functions of
  the form $h*h^*$ with $h\in C_c(D)$.
\end{enumerate}
\item[\rm{(4)}] There exists precisely one positive character $\alpha_0\in S$
by \cite{V1}, \cite{BH}.
\item[\rm{(5)}] If $\mu\in M^1(\hat D)$ satisfies $\check\mu\ge0$ on $D$, then 
its support $supp\>\mu$ contains at least one positive character; see \cite{V3}. 
\end{enumerate}
\end{facts}

In contrast to l.c.a.~groups, pointwise products of positive definite
functions on $D$ are not necessarily positive definite; see e.g.~Section
9.1C of \cite{J} for an example with $|D|=3$. However, in some cases
 positive definiteness is preserved under pointwise products. 

If for all $\alpha,\beta\in \hat D$ (or a subset of  $\hat D$ like
$S$) the products $\alpha\beta$ are positive definite, 
then by Bochner's theorem \ref{facts-hy}(1), there are probability measures 
$\delta_\alpha\hat*\delta_\beta\in M^1(\hat D)$ with 
$(\delta_\alpha\hat*\delta_\beta)^\vee=\alpha\beta$, i.e., we obtain dual 
positive product formulas as claimed in Section 1. Under additional
conditions, $(\hat D,\hat*)$ then  carries a dual hypergroup structure
with ${\bf 1}$ as identity and complex conjugation as involution. This  for
instance holds for all compact commutative double coset hypergroups $G//H$
by  \cite{Du}. For non-compact Gelfand pairs $(G,H)$ it is known that
 there are  dual positive  
convolutions on $S$; see \cite{J}, \cite{V3} for details. These convolutions usually do not
generate a dual hypergroup structure. Moreover, it is possible here that $\alpha\beta$
is not positive definite on $D$ for some  $\alpha,\beta\in\hat D$;
 see \cite{V4} for discrete examples.

The following result extends Theorem 2.1(4) of \cite{V3} and will be needed below:

\begin{proposition}\label{support-condition-pl}
Let $(D,*)$ be a commutative hypergroup. Let $\alpha\in C_b(D)$ be a function on $D$ such that  $\alpha\cdot\beta$
 is positive definite for each character $\beta\in S$ in the support of the Plancherel measure. 
Then for each  $\beta\in S$ there is a unique measure $\mu\in M_b^+(S)$ with  $\alpha\cdot\beta=\check \mu$.
\end{proposition}

\begin{proof}
Fix  $\beta\in S$. By \ref{facts-hy}(3) there exists a sequence of positive definite functions $f_n$ in $C_c(D)$
which tend locally uniformly to $\beta$. Moreover, again by  \ref{facts-hy}(3), each $f_n$ has the form
 $f_n=\check \mu_n$ for some  $\mu_n\in M_b^+(S)$. We conclude easily 
from the assumption of the proposition that the functions
$\alpha\cdot f_n= \alpha\cdot\check \mu_n=\int_S \alpha\gamma\> d\mu(\gamma)\in C_c(D) $
 are positive definite for all $n$. As these functions
 tend locally uniformly to $\alpha\cdot\beta$, we obtain from \ref{facts-hy}(3) that  $\alpha\cdot\beta=\check \mu$
for some $\mu\in M_b^+(S)$ as claimed.
\end{proof}

\section{Discrete association schemes}

In this section we briefly recapitulate two discrete generalizations of classic
finite association schemes
from \cite{Voit16} as announced in the introduction.  For classic
finite association schemes we refer to the monographs \cite{BI}, \cite{B}.
This section is useful to understand the non-discrete generalization
in the next section which is technically more involved.

The first extension from the finite to the possibly infinite case is  canonical:

\begin{definition}\label{def-discrete-asso}
 Let $X, D$ be nonempty, at most countable 
sets and $(R_i)_{i\in D}$ a disjoint
  partition of $X\times X$ with $R_i\ne\emptyset$ for $i\in D$ and the
  following properties:
\begin{enumerate}
\item[\rm{(1)}] There exists $e\in D$ with $R_e=\{(x,x): x\in X\}$.
\item[\rm{(2)}] There exists an involution $i\mapsto \bar i$ on $D$ such that
  for $i\in D$,
$R_{\bar i}=\{(y,x):\> (x,y)\in R_i\}$.
\item[\rm{(3)}] For all $i,j,k\in D$ and $(x,y)\in R_k$, the number
$$p_{i,j}^k:=|\{z\in X:\> (x,z)\in R_i\>\> {\rm and }\>\> (z,y)\in R_j\}|$$
is finite and independent of $(x,y)\in R_k$.
\end{enumerate}
Then $\Lambda:=(X,D, (R_i)_{i\in D})$ is called an association scheme with
intersection numbers $(p_{i,j}^k)_{i,j,k\in D}$ and identity $e$. 

An association scheme is called commutative if $p_{i,j}^k=p_{j,i}^k$
 for all $i,j,k\in D$. It is called symmetric (or hermitian) if
the involution on $D$ is the identity.
Moreover, it is called finite, if so are $X$ and $D$.
\end{definition}

\begin{facts}
Now let $\Lambda:=(X,D, (R_i)_{i\in D})$ be  an association scheme according to Definition 
\ref{def-discrete-asso}.
Following \cite{BI}, \cite{Voit16}, we form the 
 adjacency 
matrices $A_i\in \mathbb R^{X\times X}$ ($i\in D$) with
$$(A_i)_{x,y}:=
\left\{ \begin{array}{cc}
\displaystyle 1 &
 {\rm if}\>\> (x,y)\in R_i
\\     0&  {\rm otherwise}\\     \end{array} \right. 
 \quad\quad\quad \quad(i\in D,\>
x,y\in X).$$
The adjacency matrices  $A_i$  have the following obvious properties: 
\begin{enumerate}
\item[\rm{(1)}] $A_e$ is the identity matrix $I_X$.
\item[\rm{(2)}] $\sum_{i\in D} A_i$ is the matrix $J_X$ whose entries are all
  equal to 1.
\item[\rm{(3)}] $A_i^T = A_{\bar i}$ for $i\in D$.
\item[\rm{(4)}] For all $i\in D$ and all rows and columns of $A_i$, all
  entries are equal to
  zero except for finitely many cases.
\item[\rm{(5)}] For $i,j\in D$, the usual matrix product $A_iA_j$ exists, and
 $A_iA_j= \sum_{k\in D} p_{i,j}^k A_k$.
\item[\rm{(6)}]  $\Lambda$ is commutative if and only
  if $A_iA_j=A_jA_i$ for all  $i,j\in D$.
\item[\rm{(7)}]  $\Lambda$ is symmetric if and
  only if  all $A_i$ are symmetric. 
\item[\rm{(8)}] For $i,j,k\in D$, $p_{i,j}^{k}=p_{\bar j, \bar i}^{\bar k}$.
\end{enumerate}
\end{facts}

The valency of $R_i$ or $i\in D$ is
defined as 
\begin{equation}\label{deutung-omegai}
\omega_i:=p_{i,\bar i}^e=|\{z\in X:\> (x,z)\in R_i\}|\in\mathbb N
\end{equation}
where $\omega_i$ is independent of $x\in X$.
Therefore, the  renormalized matrices
  $S_i:= \frac{1}{\omega_i}A_i\in \mathbb R^{X\times X}$ are stochastic, i.e.,
  all  row sums are equal to 1.
 The stochastic matrices $S_i$ satisfy
\begin{equation}\label{rel-stoch-matrices}
\textstyle S_iS_j= \sum_{k\in D} \frac{ \omega_k}{\omega_i\omega_j}
p_{i,j}^k S_k  \quad\quad\text{ for}  \quad\quad i,j\in D.
\end{equation}
We next discuss a property of association schemes which is always valid in the
finite case,  but not necessarily in infinite cases.

\begin{definition} An association scheme  with
 valencies $\omega_i$ is called unimodular if $\omega_i=\omega_{\bar i}$
 for all $i\in D$. 
\end{definition}

We collect some facts from Section 3  of \cite{Voit16}:

\begin{facts}\label{unimod}
\begin{enumerate}
\item[\rm{(1)}] If an association scheme is commutative or finite,
 then it is unimodular.
\item[\rm{(2)}] An association scheme  is unimodular if
  and only if
 the associated discrete hypergroup  is unimodular.
\item[\rm{(3)}] If  $(X,D,(R_i)_{i\in D})$ is unimodular, then
$S_i^T=S_{\bar i}$ for all $i\in D$.
\item[\rm{(4)}] There exist non-unimodular association schemes.
\end{enumerate}
\end{facts}

The observations above and in particular Eq.~(\ref{rel-stoch-matrices}) 
were used in Section 5 of \cite{Voit16} for the following
generalization of Definition \ref{def-discrete-asso} of association schemes:

\begin{definition}\label{def-gen-discrete-asso}
 Let $X, D$ be nonempty, at most countable 
sets and $(R_i)_{i\in D}$ a disjoint
  partition of $X\times X$ with $R_i\ne\emptyset$ for $i\in D$. 
Let $\tilde S_i\in \mathbb R^{X\times X}$ for $i\in D$ be stochastic
matrices. 
Assume that:
\begin{enumerate}
\item[\rm{(1)}] For all $i,j,k\in D$ and $(x,y)\in R_k$, the number
$$p_{i,j}^k:=|\{z\in X:\> (x,z)\in R_i\>\> {\rm and }\>\> (z,y)\in R_j\}|$$
is finite and independent of $(x,y)\in R_k$.
\item[\rm{(2)}] For all $i\in D$ and $x,y\in X$, $\tilde S_i(x,y)>0$ 
if and only if $(x,y)\in R_i$.
\item[\rm{(3)}] For all $i,j,k\in D$ there exist (necessarily nonnegative)  numbers
$\tilde p_{i,j}^k$ with $\tilde S_i\tilde S_j=\sum_{k\in D} \tilde p_{i,j}^k\tilde S_k$.
\item[\rm{(4)}] There exists an identity $e\in D$ with $\tilde S_e=I_X$ as identity matrix.
\item[\rm{(5)}] There exists a positive measure $\omega_X\in M^+(X)$ with
  $supp\>\omega_X =X$ and an involution $i\mapsto \bar i$ on $D$ 
such that
for all $i\in D$, $x,y\in X$,
  $$\omega_X(\{y\})\tilde  S_{\bar i}(y,x)=\omega_X(\{x\})\tilde   S_{ i}(x,y).$$
\end{enumerate}
Then $\Lambda:=(X,D, (R_i)_{i\in D},(\tilde  S_i)_{i\in D})$ 
is called a generalized association scheme. 

 $\Lambda$ is called commutative if $\tilde S_i\tilde S_j=\tilde S_j\tilde S_i$
 for all $i,j\in D$. It is called symmetric  if
the involution  is the identity. $\Lambda$ is called finite, if so are $X$ and  $D$.
\end{definition}

\begin{remarks}\label{ex-ass-gen-ass}
\begin{enumerate}
\item[\rm{(1)}] If  $\Lambda=(X,D, (R_i)_{i\in D})$ 
is an unimodular association scheme with the associated
stochastic matrices $ (S_i)_{i\in D}$ as above,
 then 
$(X,D, (R_i)_{i\in D}, (S_i)_{i\in D})$  is a generalized association
scheme. In fact, axioms (1)-(4) are clear, and for axiom (5) we take the
involution of $\Lambda$ and $\omega_X$ as the counting measure on $X$. 
Fact \ref{unimod}(3)
and  unimodularity then imply axiom (5). Clearly, notions like commutativity and symmetry are preserved.
For details we refer to \cite{Voit16}.
\item[\rm{(2)}] If $(X,D, (R_i)_{i\in D}, (\tilde S_i)_{i\in D})$ is a  generalized association
scheme, then  $(X,D, (R_i)_{i\in D})$ is an association
scheme. If this scheme is unimodular, then we may form the two generalized association
schemes  $(X,D, (R_i)_{i\in D}, (\tilde S_i)_{i\in D})$ and  $(X,D, (R_i)_{i\in D}, (\tilde S_i)_{i\in D})$
 on the same spaces $X,D$
where the second one is formed according to (1).

\end{enumerate}
\end{remarks}

For examples of infinite commutative association schemes and of generalized association schemes,
 which are no association schemes, we refer to \cite{Voit16}, \cite{Voit17} and to Section 10 below. 

Generalized association schemes always lead to discrete  hypergroups; see Prp.~5.4 of \cite{Voit16}:

 \begin{proposition}\label{prp-gen-asso-to-hypergroup}
 Let $\Lambda:=(X, D, (R_i)_{i\in D}, (\tilde S_i)_{i\in D})$ be a generalized association
  scheme with deformed intersection numbers $\tilde p_{i,j}^k$.
Then the product $\tilde*$ with
$$\textstyle\delta_i\tilde*\delta_j:=\sum_{k\in D}
 \tilde p_{i,j}^k \delta_k$$
can be extended uniquely to an associative, bilinear, 
$\|\>.\>\|_{TV}$-continuous mapping on $M_b(D)$.
 $(D,\tilde*)$  is a discrete hypergroup with identity $e$ and the involution on $D$ 
from Definition \ref{def-gen-discrete-asso}(5).  $(D,\tilde*)$ has
 the left and right Haar measure
$$\textstyle\Omega_l:=\sum_{i\in D} \omega_i \delta_i
\quad\text{and}\quad\Omega_r:=
\sum_{i\in D} \omega_{\bar i} \delta_i \quad\text{with}\quad
\omega_i:=\frac{1}{\tilde p_{i,\bar i}^e}>0 \quad(i\in D) $$
 respectively. $(D,\tilde*)$   is commutative or  
 symmetric  if and only if so is $\Lambda$.
\end{proposition}

For association schemes there is a corresponding result; see Proposition 3.8 of \cite{Voit16}.
In fact, the associated hypergroup convolution algebras 
 are just the Bose-Mesner
algebras for finite association schemes in \cite{BI}.

\begin{proposition}\label{prp-asso-to-hypergroup}
 Let $\Lambda:=(X, D, (R_i)_{i\in D})$ be an association
  scheme with intersection numbers $p_{i,j}^k$ and valencies $\omega_i$.
Then the product $*$ with
$$\delta_i*\delta_j:=\sum_{k\in D}\frac{\omega_{ k}}{\omega_{ i}\omega_{ j}} 
\cdot p_{i,j}^k \delta_k$$
can be extended uniquely to an associative, bilinear, 
$\|\>.\>\|_{TV}$-continuous mapping on $M_b(D)$.
 $(D,*)$  is a discrete hypergroup with the left and right Haar measure
\begin{equation}\label{const-haar}
\Omega_l:=\sum_{i\in D} \omega_i \delta_i
\quad\quad\text{and}\quad\quad \Omega_r:=\Omega_l^*:=
\sum_{i\in D} \omega_{\bar i} \delta_i   
\end{equation}
 respectively.  $(D,*)$ is commutative, 
 symmetric, or unimodular  if and only if so is $\Lambda$.
\end{proposition}

A generalized association scheme is called unimodular if so is the associated hypergroup. Clearly,  unimodular
association schemes always lead
to unimodular generalized association schemes in Remark \ref{ex-ass-gen-ass}(1).

\section{Continuous association schemes}

In this section we  propose and discuss a system of axioms which extends the
notion of generalized association schemes above to a
continuous setting where we replace the stochastic matrices $(S_i)_{i\in D}$
 by Markov kernels.  We briefly recapitulate some well-known notations on Markov kernels:

\begin{definition}\label{def-makov-kernel}
Let $X,Y$ be locally compact spaces equipped with the
 associated Borel $\sigma$-algebras ${\mathcal B}(X),{\mathcal B}(Y)$.
 A Markov kernel $K$
 from $X$ to $Y$ is  as
 a mapping $K:X\times {\mathcal B}(Y)\to [0,1]$ such that
\begin{enumerate}
\item[\rm{(1)}] for all $x\in X$, the mapping $K(x,.): {\mathcal B}(Y)\to [0,1]$,
$A\mapsto K(x,A)$ is a probability measure on $(Y,{\mathcal B}(Y))$, and 
\item[\rm{(2)}]  for  $A\in {\mathcal B}(Y)$,  the mapping   $K(.,A): X\to [0,1]$,
$x\mapsto K(x,A)$ is measurable.
\end{enumerate}
Consider the Banach spaces $F_b(X),F_b(Y)$ of all $\mathbb C$-valued
bounded measurable functions on $X,Y$ with the supremum norm. Then for any 
 Markov kernel $K$ from $X$ to $Y$ we define the associated transition operator
$$T_K:F_b(Y)\to F_b(X) \quad\quad\text{with} \quad\quad 
T_Kf(x):=\int_Y f(y)\> K(x,dy) \quad\quad(x\in X).$$
Clearly, $K$ is determined uniquely by the operatator $T_K$.

We say that a Markov kernel $K$ is continuous if $T_K(C_b(Y))\subset C_b(X)$. 
If $ T_K(C_0(Y))\subset C_0(X)$, then  $K$ is called a Feller kernel.

Let us also recapitulate the  composition
$$K_1\circ K_{2}(x,A):=\int_{X} K_{2}(y,A)\> K_1(x,dy) \quad\quad (x\in X,\> A\in {\mathcal B}(X))$$
of Markov-kernels $K_1,K_2$ on $X$, i.e., from $X$ to $X$. We then have 
 $T_{K_1\circ K_{2}}=T_{K_1}\circ T_{K_{2}}$, and  the composition of
  kernels and transition operators are associative.

A (nontrivial) positive Radon measure $\omega\in M^+(X)$ is called
$K$-invariant w.r.t. a Markov kernel $K$ on $X$, if
$\int_B K(x,A)\> d\omega(x)=\omega(A)\in [0,\infty]$ for $A\in {\mathcal B}(X)$.
\end{definition}

We now turn to continuous association schemes. We here consider two second countable,
 locally compact spaces $X,D$ together with 
 a continuous 
Markov kernel $K$ from $X\times D$ to $X$ with  transition
operator $T_K:C_b(X)\to C_b(X\times D)$. For each $h\in D$, we then define
the Markov kernels
$$K_h(x,A):= K(x,h;A):= K((x,h), A) \quad\quad (x\in X, A\in {\mathcal B}(X))$$
 on $X$ with transition operators
$$T_h:C_b(X)\to C_b(X),\quad
T_hf(x):=T_{K_h}f(x)=\int_X f(y)\> K(x,h;dy)=T_Kf(x,h).$$

With these notions we now define continuous association schemes.
Unfortunately, the definition is more involved than in the discrete case due to the additional
continuity assumptions and some other restrictions which are satisfied in the discrete case automatically.
%For simplicity we also incorporate some hidden unimodularity condition; cf. Remark \ref{ex-ass-gen-ass}(1)
% for a discussion in the discrete case.

\begin{definition}\label{central-def}
Let $X$ and $D$ be second countable, 
locally compact spaces, and  $K$  a continuous 
Markov kernel from
$X\times D$ to $X$. Then  $(X,D,K)$ is called a continuous 
 association scheme (or, for short, {\bf CAS}), if the following holds:
\begin{enumerate}
\item[\rm{(1)}] {\bf(Compact support:)} For $x\in X$, $h\in D$, the support  $supp \>
  K(x,h;.)$ is compact, and the mapping $(x,h)\mapsto supp \> K(x,h;.)$
from $X\times D$ into the space $\mathcal C(X)$  of all compact subsets of $X$ is
continuous w.r.t.~the Hausdorff topology on $\mathcal C(X)$.
\item[\rm{(2)}]  {\bf(Partition property:)} For each $x\in X$,
 the compact sets
$supp \> K(x,h;.)$ ($h\in D$) form a partition of $X$, and the associated unique map $\pi:
X\times X\to D$ with  $y\in supp \> K(x,\pi(x,y);.)$ for $x,y\in X$
is continuous.
\item[\rm{(3)}] {\bf(Hypergroup property:)}
$D$ carries a hypergroup structure $(D,*)$ such that for all $h_1,h_2\in D$,
$x\in X$,  and $A\in{\mathcal B}(X)$
\begin{equation}\label{ess-hg-eq}
K_{h_1}\circ K_{h_2}(x,A)= \int_D K_h(x,A)\> d(\delta_{h_1}*\delta_{h_2})(h).
\end{equation}
Moreover, the identity $e\in D$ satisfies $K(x,e;.)=\delta_x$ for $x\in X$.
\item[\rm{(4)}]  {\bf(Invariant measure on $X$:)}
There exists a positive Radon measure $\omega_X\in M^+(X)$ with $supp\>\omega_X=X$,
such that for the continuous hypergroup involution $\bar.:D\to D$
and all
$h\in D$, $f,g\in C_c(X)$,
\begin{equation}\label{adjoint-operator}
\int_X T_{\bar h}f \cdot g \> d\omega_X = \int_Xf \cdot T_{h}g\> d\omega_X.
\end{equation}
\end{enumerate}
 A CAS $(X,D,K)$ is called commutative, symmetric, or unimodular if so is the hypergroup $(D,*)$.
It is called compact or discrete if so are $D$ and $X$.
\end{definition}

Clearly, a CAS $(X,D,K)$ is commutative if and only  the Markov kernels $K_h$ ($h\in D$) commute.

We  present some standard examples of CAS. 
Further examples are given later on.

\begin{proposition}\label{discrete-gen-ass-to-cont-asso-scheme}
Let $(X,D,(R_h)_{h\in D}, (\tilde S_h)_{h\in D})$ be a generalized  association sche\-me 
as in Definition \ref{def-gen-discrete-asso}.
Then
$$K((x,h),A):=\sum_{y\in A} \tilde S_h(x,y) \quad\quad(x\in X,\> h\in D,\> A\subset X)$$
defines a Markov kernel from $X\times D$ to $X$, and $(X,D,K)$ is a CAS.
\end{proposition}

\begin{proof}
$K$ is a obviously a Markov kernel from $X\times D$ to $X$. Moreover, as
$X$ and $D$ are discrete, all topological axioms are trivial.
Furthermore, fact (4) after Definition \ref{def-discrete-asso} in combination with Definition 
\ref{def-gen-discrete-asso}(2) show that $supp\> K((x,h),.)$ is finite for all $x\in X$ and $h\in D$.
This shows \ref{central-def}(1). Moreover, \ref{central-def}(2) is obvious, and \ref{central-def}(3)
follows from Proposition \ref{prp-gen-asso-to-hypergroup}. Finally, \ref{central-def}(3) follows from 
\ref{def-gen-discrete-asso}(5).
\end{proof}

Remark \ref{ex-ass-gen-ass}(1) and  Proposition \ref{discrete-gen-ass-to-cont-asso-scheme} imply:

\begin{corollary}\label{discrete-ass-to-cont-asso-scheme}
Let $(X,D,(R_h)_{h\in D})$ be an unimodular  association scheme 
with the stochastic matrices $(S_h)_{h\in D}$  as defined after Definition \ref{def-discrete-asso}.
Then
$$K((x,h),A):=\sum_{y\in A}  S_h(x,y) \quad\quad(x\in X,\> h\in D,\> A\subset X)$$
defines a Markov kernel from $X\times D$ to $X$, and $(X,D,K)$ is an unimodular CAS.
\end{corollary}

\begin{remark} Proposition \ref{discrete-gen-ass-to-cont-asso-scheme} admits the following partial converse statement:

Let $(X,D,K)$ be a discrete CAS. Define the stochastic matrices
$(\tilde S_h)_{x,y}:=K_h(x,\{y\})$ for $h\in D$, $x,y\in X$ as well as the sets 
$R_h:=\{(x,y)\in X\times X:\> (\tilde S_h)_{x,y}>0\}$. Then $(X,D, (R_h)_{h\in D}, (\tilde S_h)_{h\in D})$ 
satisfies almost all axioms
of a generalized association scheme in \ref{def-gen-discrete-asso}. In fact  $(R_h)_{h\in D}$ forms a partition of 
$ X\times X$, and the axioms (2)-(5) of \ref{def-gen-discrete-asso} hold. We do not know at the moment 
 whether also (1) in
 \ref{def-gen-discrete-asso} holds automatically. We come back to this problem later on in the finite case.
\end{remark}

Here is a further standard class of examples of CAS.

\begin{proposition}\label{group-to-cont-asso-scheme}
Let $H$ be a compact subgroup of a locally compact unimodular 
group $G$ with 
normalized Haar measure $\omega_H\in M^1(H)\subset M^1(G)$. Then the quotient
$X:=G/H$ and the double coset space $D:=G//H$ are locally compact 
w.r.t.~the quotient topology, and the canonical projections
$$p_G:G\to G/H, \quad p_G(x):=xH, 
\quad
 p_{G/H}:G/H\to G//H, \quad p_{G/H}(xH):=HxH$$
are continuous, open, and closed.
Moreover,
\begin{equation}\label{def-K-allg-group-case}
K((xH, HhH), A):= p_G(\delta_x*\omega_H*\delta_h *\omega_H)(A)
\quad\quad(x,h\in G,\> A\in{\mathcal B}(X))
\end{equation}
establishes a well-defined Markov kernel from $X\times D$ to $X$, and 
$(X,D,K)$ is an unimodular CAS.
\end{proposition}

\begin{proof} The topological statements about $p_G, p_{G/H}$ are well-known from the
 theory of locally compact groups.
We  next check that $K$ is a well-defined continuous Markov kernel,
and that $(X,D,K)$ is a CAS.
Clearly, by its construction, $K$ is a probability measure w.r.t. the variable $A$, and the definition
 of $K((xH, HhH), A)$ is
independent of the representatives $x,h$ of $xH$ and $HhH$ respectively. 
Before checking that the maps
$$D_A:X\times D\to[0,1], \quad (x,h)\mapsto K((x,h), A)$$ 
are measurable for all Borel set $A$, we investigate the associated transition operator
$T_K$. For  $f\in C_b(X)$ we  have that
\begin{align}\label{comp-TK}
T_Kf(xH,HhH)&= \int_G f(yH)\> d(\delta_x*\omega_H*\delta_h *\omega_H)(y)
\\
&=\int_G \int_Gf(xz_1hz_2H)\>d\omega_H (z_1)\>d\omega_H (z_2)
\notag\end{align}
is continuous in $x,h\in G$. As the projections $p_G,p_{G/H}$ are open,
 it follows  that the map $(xH, HhH)\mapsto T_Kf(xH,HhH)$ is continuous.
If we have proved that the maps $D_A$ are measurable for all Borel set $A$, we  conclude that 
 $K$ is  a continuous Markov kernel as claimed. To prove measurabilty of the maps $D_A$, we first 
choose a compact set $A$. As the characteristic function $1_A$ of $A$ is a monotone limit of 
functions $f_n\in C_c(X)$, we obtain from the theorem of monotone convergence 
that the continuous functions  $(xH, HhH)\mapsto T_Kf_n(xH,HhH)$ tend to 
 $(xH, HhH)\mapsto T_K1_A(xH,HhH)= K((xH,HhH), A)=D_A(xH,HhH)$. This proves that $D_A$ is
 measurable for $A$ compact. For the general case we use Dynkin systems. In fact, 
$$\cal D:=\{ A\in \cal B(X): \> D_A\quad\text{is  measurable}\}$$
is a Dynkin system which contains the set $\cal K$ of all compacta in $X$,
 where  $\cal K$ is closed under intersections.
Therefore, the $\sigma$-algebra $\sigma(\cal K)$ and the Dynkin system $\cal D(\cal K)$ generated by  $\cal K$ 
satisfy $\cal B(X)=\sigma(\cal K)=\cal D(\cal K)\subset\cal D$. Hence, $D_A$ is measurable for all
 Borel set $A$ as claimed.

It is now standard to check the axioms (1)--(4) of Definition \ref{central-def}. In fact, 
 the compact support in (1) is clear, and the continuity w.r.t.~the Hausdorff topology 
follows in the same way as in \cite{J} for double coset hypergroups.
 Moreover, the projection $\pi$ in (2) is  given by
$\pi(xH,yH):=Hx^{-1}yH$ and thus continuous, while the partition property in (2) is clear.
 Furthermore, it is well-known $D=G//H$ is a double coset hypergroup  with the convolution
$$\delta_{Hh_1H}*\delta_{Hh_2H}=\int_H \delta_{Hh_1hh_2H}\>d\omega_H(h);$$
see \ref{Doublecoset} and \cite{J} for details. Moreover, for $x,h_1,h_2\in G$ and Borel sets $A\subset X$,
\begin{align}
K_{Hh_1H}&\circ K_{Hh_2H}(xH, A)= \int_X  K_{Hh_2H}(w, A)\> K_{Hh_1H}(xH, dw)\notag\\
&=
\int_G (\delta_y*\omega_H*\delta_{h_2}*\omega_H)(p_G^{-1}(A)) \> d(\delta_x*\omega_H*\delta_{h_1}*\omega_H)(y)
\notag\\
&=
(\delta_x*\omega_H*\delta_{h_1}*\omega_H*\delta_{h_2}*\omega_H)(p_G^{-1}(A))
\notag\end{align}
and 
\begin{align}
\int_D K_w(xH,A)\> &d(\delta_{Hh_1H}*\delta_{Hh_2H})(w)= \int_H K_{Hh_1hh_2H}(xH,A)\> d\omega_H(h)
\notag\\
&= \int_H (p_G(\delta_x*\omega_H*\delta_{h_1hh_2}*\omega_H))(A)\> d\omega_H(h)\notag\\
&=
(\delta_x*\omega_H*\delta_{h_1}*\omega_H*\delta_{h_2}*\omega_H)(p_G^{-1}(A)),
\notag\end{align}
which proves  Eq.~(\ref{ess-hg-eq}).

We next turn to \ref{central-def}(4). Let  $\omega_G$ be some Haar measure $\omega_G$ of $G$ and take
 its projection $\omega_{X}:=p_G(\omega_G)$ as invariant measure on $X$.
 To check (\ref{adjoint-operator}), we take $f,g\in C_c(X)$, $h\in G$ 
and observe 
from (\ref{comp-TK}) and unimodularity that
\begin{align}
\int_X f\cdot T_hg\> d\omega_X&=\int_G\int_G  f(xH) g(xyH) \>
 d(\omega_H*\delta_h*\omega_H)(y) \> d\omega_G(x) \notag\\
&=\int_G\int_G f(xyH) g(xH)\>
 d(\omega_H*\delta_{h^{-1}}*\omega_H)(y) \> d\omega_G(x)\notag\\
&=\int_X  T_{\bar h}f\cdot g\> d\omega_X.
\notag\end{align}
This completes the proof of \ref{central-def}(4). Finally, as the canonical projection of the Haar measure 
$\omega_G$ to $G//H=D$ is a Haar measure of the double coset hypergroup $(D;*)$ (see \cite{J}), $(X,D,K)$ is unimodular.
\end{proof}

We next proceed with the theory of CAS.
We first collect some obvious consequences from Definition \ref{central-def} 
for a CAS $(X,D,K)$.

\begin{facts}\label{facts-basic-gen-ass}
\begin{enumerate}
\item[\rm{(1)}] Property \ref{central-def}(2) ensures that for
  $h_1,h_2\in D$ with $K_{h_1}=K_{h_2}$ we have $h_1=h_2$. 
Therefore, the convolution $*$ of $(D,*)$ is determined uniquely by the kernels $K_h$, $h\in D$.
\item[\rm{(2)}]  $T_e$ is the identity
  operator. 
\item[\rm{(3)}] It can be easily shown that the adjoint relation (\ref{adjoint-operator}) 
holds for further classes of functions. In particular, it can be easily seen that
(\ref{adjoint-operator}) holds for all $f\in C_b(X)$ and 
 $g\in  L^1(X,\omega_X)$. Taking in particular $f={\bf 1}$, we obtain 
\begin{equation}\label{invariant-omega}
\int_X g \> d\omega_X=\int_X T_hg \> d\omega_X
\end{equation}
for all $h\in D$ and $g\in  L^1(X,\omega_X)$. Taking
 $g={\bf  1}_A$ for measurable sets $A\subset X$, we conclude that $\omega_X$ is
$K_h$-invariant for all  $h\in D$.
\item[\rm{(4)}] For all $h\in D$ and $p\in[1,\infty[$, the operator $T_h$ associated with the kernel $K_h$ on $X$ 
is
 a continuous linear operator on $L^p(X,\omega_X)$ with 
$\|T_h\|\le1$. In fact, for $f\in L^p(X,\omega_X)$, the H\"older inequality
    and the invariance of $\omega_X$ imply
\begin{align}\|T_hf\|_p^p&=\int_X \Bigl|\int_X f(y)\> K_h(x,dy)\Bigr|^p\>
  d\omega_X(x) \notag\\
&\le
\int_X\Bigl(\int_X |f(y)|^p\> K_h(x,dy)\Bigr)
 \Bigl(\int_X 1\> K_h(x,dy)\Bigr)^{p/q}d\omega_X(x) =\|f\|_p^p.
\notag\end{align}
\item[\rm{(5)}] For all $h\in D$, $T_h$ clearly is a continuous operator on $(C_b(X),\|.\|_\infty) $.
\item[\rm{(6)}] As $C_c(X)$ 
is $\|.\|_2$-dense in $L^2(X,\omega_X)$, and as $T_h$ is $\|.\|_2$-continuous, the adjoint
relation (\ref{adjoint-operator}) implies that $T_{\bar h}$
 is the  adjoint operator $T_h^*$ on $L^2(X,\omega_X)$.
\item[\rm{(7)}] For all $x,y\in X$, $\pi(y,x)=\overline{\pi(x,y)}$. 

In fact, for each $h\in D$, $\pi(x,y)\ne h$ means that there is a
 neighborhood $U_y$ of $y$ in $x$ such that for all $g\in C_c(X)$ with $supp\>g\subset U_y$ we have
$T_hg(x)=0$. As $T_hg$ is continuous by our assumptions, we conclude that $\pi(x,y)\ne h$ is equivalent to the fact that
for all $\epsilon>0$ there are neighborhoods $U_x,U_y$ of $x,y$ respectively such that for all
$f,g\in C_c(X)$ with $supp\>f\subset U_x$,  $supp\>g\subset U_y$, and with $\|f\|_{1,\omega_X}=\|g\|_{1,\omega_X}=1$
we have
$|\int_X f(w)\> T_hg(w)\> d\omega_X(w)|\le\epsilon$. This  and  the adjoint relation \ref{central-def}(4) lead 
to the claim.
\end{enumerate}
\end{facts}

We next study some topological properties of $\pi$. For this recapitulate that in our setting,
 the  Hausdorff topology
on $\mathcal C(X)$ agrees with the so-called Michael topology in \cite{J} or Section 1.1 of \cite{BH};
 see  e.g.~\cite{KS}. In particular, by these references (see in particular (2.5F) of  \cite{J}),
  $\mathcal C(X)$ is locally compact, and for each compactum $\Omega\subset \mathcal C(X)$, the set
$\bigcup_{A\subset \Omega}A\subset X$ is compact. With these preparation we obtain:

\begin{lemma}\label{compact-union} 
\begin{enumerate}
\item[\rm{(1)}] For compact sets $K\subset X$ and $L\subset D$, the set 
$\bigcup_{h\in L, x\in K}supp\> K_h(x,.)\subset X$
is compact.
\item[\rm{(2)}] For each $x\in X$, the projection $\pi_x:X\to D$,  $\pi_x(y):=\pi(x,y)$ is open,
 closed and proper, i.e., $\pi_x^{-1}(A)\subset X$ is compact for each compactum $A\subset D$.
\end{enumerate}
\end{lemma}

\begin{proof} Part (1) follows from \ref{central-def}(1) and the remark about  $\mathcal C(X)$ above.

These facts also imply that  $\pi_x$ is proper, as  for each compactum $A\subset D$, the set
$$\pi_x^{-1}(A)= \bigcup_{y\in A} \pi_x^{-1}(y) =\bigcup_{y\in A} supp\> K_h(x,.)$$
is compact.  Problem 5 of Section XI.6 of \cite{Dug} now implies that  $\pi_x$ is also closed.

We finally show that $\pi_x$ is open. For this assume that there is some neighborhood
 $U\subset X$ of some $y\in X$ such that $\pi_x(U)\subset D$ is no  neighborhood of $\pi_x(y)$.
 This means that there is a sequence $(h_n)_n\subset D\setminus \pi_x(U)$ with $h_n\to \pi(x,y)$.
 Hence, by \ref{central-def}(1), the compacta $\pi_x^{-1}(h_n)= supp\> K_{h_n}(x,.)$ tend to 
 $\pi_x^{-1}(\pi(x,y))$. On the other hand, $\pi_x^{-1}(h_n)\cap U=\emptyset$, and  $\pi_x^{-1}(\pi(x,y))$ contains
$x\in U$. This leads to a contradiction. Hence,  $\pi_x$ is open. 
\end{proof}

\begin{lemma}\label{uniform-cont}
Euch $g\in C_0(X)$ is uniformly continuous in the  sense that for each $\epsilon>0$ there exists a
neighborhood $U\subset D$ of the identity $e$ such that for all
$x,y\in X$ with $\pi(x,y)\in U$, $|g(x)-g(y)|\le\epsilon$.
\end{lemma}

\begin{proof} The proof is similar to  a corresponding result for hypergroups; see e.g.~1.2.28 of \cite{BH}.

Fix some $\epsilon>0$.   Choose some compactum $G\subset X$ such that $|g(x)|\le \epsilon/2$ for
 $x\in X\setminus G$. For each $x\in X$ we take some open neighborhood $W_x\subset X$ with 
$|g(y)-g(x)|\le\epsilon/2$ for $y\in W_x$. Now choose open  neighborhoods $\tilde U_x\subset D$ of $e$
with $\{y\in X: \> \pi(x,y)\in \tilde U_x\}\subset W_x$.
 By a basic result on hypergroups we find open symmetric  neighborhoods $U_x\subset D$
of $e$ with $U_x*U_x\subset \tilde U_x $. We now consider the open set $V_x:=\{y\in X: \> \pi(x,y)\in U_x\}$
which cover the compactum $G$. Choose $n\in \mathbb N$ and $x_1,\ldots , x_n\in X$ with 
$G\subset \bigcup_{l=1,\ldots,n} V_{x_l}$, and define the open  neighborhood $U:=\bigcap_{l=1,\ldots,n}  U_{x_l}\subset D$.

Now consider $x\in G$ and $y\in X$ with $\pi(x,y)\in U$. We find $l$ with $x\in V_{x_l}\subset W_{x_l} $. As
then $\pi(x_l,x)\in  U_{x_l}$ and  $\pi(x,y)\in  U_{x_l}$, we obtain
  $\pi(x_l,y)\in  U_{x_l}* U_{x_l}\subset \tilde U_{x_l}$, which implies  $y\in  W_{x_l} $. Hence, by the definition of 
$ W_{x_l} $, $|g(y)-g(x)|\le\epsilon$. As this also holds
 for all $x,y\in X\setminus G$, the proof is complete.
\end{proof}

\begin{lemma}\label{compact-preserved} Let $h\in D$. Then:
\begin{enumerate}
\item[\rm{(1)}] If $f\in C_c(X)$, then $T_hf\in C_c(X)$.
\item[\rm{(2)}] If $f\in C_0(X)$, then $T_hf\in C_0(X)$; in other words, the kernels $K_h$ on $X$ are Feller kernels.
\end{enumerate}
\end{lemma}

\begin{proof} The continuity of  $T_hf$ is clear in both cases.

Now let $f\in C_c(X)$. Let $x\in X$ with $T_hf(x)\ne0$. Then, by the definition of the projection $\pi$,
$supp\> f \cap \{y\in X:\> \pi(x,y)=h\}\ne\emptyset$, and thus, by
\ref{facts-basic-gen-ass}(7), $supp\> f \cap \{y\in X:\> \pi(y,x)=\bar h\}\ne\emptyset$. This yields
$x\in\bigcup_{y\in supp \> f}supp \> K_{\bar h}(y,.)$
As this set is compact by Lemma \ref{compact-union}, we obtain $T_hf\in C_c(X)$.

Part (2) follows from part (1) and the continuity of $T_h$ w.r.t.~$\|.\|_\infty$.
\end{proof}

We next study  integrals over the operators $T_h$.

\begin{lemma}\label{op-t-mu} Let $\mu\in M_b(D)$. 
\begin{enumerate}
\item[\rm{(1)}] For each $f\in C_b(X)$,
$$T_\mu f(x):=\int_D T_h  f(x)\> d\mu(h) \quad\quad(x\in X)$$
defines a function $T_\mu f\in C_b(X)$. The operator $T_\mu$ is a continuous linear operator on $C_b(X)$
with $\|T_\mu\|\le\|\mu\|_{TV} $.  
\item[\rm{(2)}] If $supp\>\mu$ is compact, then   $T_\mu$  maps $C_c(X)$ into $C_c(X)$.
\item[\rm{(3)}] The operator   $T_\mu$  maps $C_0(X)$ into $C_0(X)$.
\item[\rm{(4)}] For each $p\in[1,\infty[$, the  operator   $T_\mu$ from (1) may be also regarded as 
 a continuous linear operator on $ L^p(X,\omega_X)$ with  $\|T_\mu\|\le\|\mu\|_{TV} $.  
 \end{enumerate}
\end{lemma}

\begin{proof}
\begin{enumerate}
\item[\rm{(1)}] $T_\mu f$ is continuous by Definition \ref{central-def}. The further statements are clear.
\item[\rm{(2)}] follows from Lemma \ref{compact-union} in the same way as Lemma \ref{compact-preserved}(1). 
\item[\rm{(3)}] If $supp\>\mu$ is compact, then (3) follows from (2) and the continuity of $T_\mu$.
On the other hand, for each $\epsilon>0$ and $\mu\in M_b(D)$ there exists a measure $\mu_\epsilon\in M_b(D)$
 with compact support and $\|\mu-\mu_\epsilon\|_{TV}\le\epsilon$. Hence $\|T_\mu- T_{\mu_\epsilon}\|\le\epsilon $.
Thus,  $T_\mu f$ is a uniform limit of functions in $C_0(X)$ which yields the claim. 
\item[\rm{(4)}] This follows from \ref{facts-basic-gen-ass}(4) and standard facts on operator-valued integrals.
\end{enumerate}
\end{proof}

We now consider the $C^*$-algebra $\mathcal B(L^2(X,\omega_X))$ of all bounded linear operators on $L^2(X,\omega_X)$
as well as the closed subspace
$A(X):=\overline{span\{T_h:h\in D\}}$. The space
$A(X)$ is closed under the involution $.^*$ on $\mathcal B(L^2(X,\omega_X))$ 
by  \ref{facts-basic-gen-ass}(5). Moreover, 
by Lemma \ref{op-t-mu}(4), we have $T_\mu\in A(X)$ for all $\mu\in M_b(D)$. In summary:

 \begin{proposition}\label{alg-hom}
\begin{enumerate}
\item[\rm{(1)}] $A(X)$ is a  $C^*$-subalgebra of  $\mathcal B(L^2(X,\omega_X))$.
\item[\rm{(2)}] The map $T:(M_b(D),*,.^*,\|. \|_{tv})\to A(X)\subset \mathcal B(L^2(X,\omega_X))$, $\mu\mapsto T_\mu$,
 is a norm-decreasing Banach-$*$-algebra homomorphism, i.e.,
 $T$ is a $*$-representation of the hypergroup $(D,*)$ on the Hilbert space $\mathcal B(L^2(X,\omega_X))$.
 \end{enumerate}
\end{proposition}

\begin{proof}
\begin{enumerate}
\item[\rm{(1)}] Let $h_1,h_2\in D$. Then, by \ref{central-def}(3),  $T_{h_1}T_{h_2}=T_{\delta_{h_1}*\delta_{h_1}}\in A(X)$.
This yields that  $A(X)$ is closed under multiplication. All further facts  are clear.
\item[\rm{(2)}] is also clear by the same arguments and Lemma \ref{op-t-mu}.
 \end{enumerate}
\end{proof}

We next study a couple of linear operators $A:C_c(X)\to C(X)$. For each such $A$  we form the scalar products
\begin{equation}\label{scalar-pr-x}
\langle Ag_1, g_2\rangle_{X}:=\int_X Ag_1(x)\cdot \overline{g_2(x)}\> d\omega_X(x) 
\quad\quad\text{for}\quad g_1,g_2\in C_c(X).
\end{equation}

Here are some examples:
Let $F\in C(X\times X)$ and form $T^F:C_c(X)\to C(X)$ by
\begin{equation}\label{TF-def}
T^Fg(x):= \int_X F(x,y)\> g(y)\> d\omega_X(y)\quad\quad(x\in X, \> g\in C_c(X)).
\end{equation}

Moreover, for all $\mu\in M_b(D)$, the operators $A:=T_\mu$ are  operators as in (\ref{scalar-pr-x}).

We now fix some left Haar measure $\omega_D$ of the hypergroup $(D,*)$. Then $L^1(D,\omega_D)$ is a
 Banach-$*$-algebra with the convolution and involution
$$f*g(x):=\int_D f(x* y)g(\bar y)\> d\omega_D(y), \quad f^*(x)=\overline{f(\bar x)} \quad\quad (x\in D).$$
 Moreover, the map $L^1(D,\omega_D)\to M_b(D)$, $ f\mapsto f\> \omega_D$ 
is an embedding of the Banach-$*$-algebra $L^1(D,\omega_D)$ into the 
Banach-$*$-algebra $M_b(D)$; see \cite{J}. For each $f\in L^1(D,\omega_D)$ 
we thus may define the linear operators $T_f:=T_{f\>\omega_D}$,
for which the results of Lemma   \ref{op-t-mu} and Proposition \ref{alg-hom} hold.

Moreover, even for $f\in C(D)$, the linear operators $T_f:C_c(X)\to C(X)$ with
\begin{equation}\label{Tuf-def}
 T_fg(x):= \int_D\int_X g(y) \> K_h(x,dy)\> f(h) \> d\omega_D(h)\quad\quad (x\in X, \> g\in C_c(X))
\end{equation}
are well-defined, as by \ref{central-def}(2) for $g\in C_c(X)$, the set 
$\pi(x,supp\> g)\subset D$ is compact. It is also clear that for 
 $f\in C(D)$ and  $f_1,f_2\in C_c(D)$, we have $f_1*f*f_2\in C(D)$ with $T_{f_1*f*f_2}=T_{f_1}T_{f}T_{f_2}$.
We shall use these facts for relations between positive definite functions
 on $(D,*)$ and   $X$ in  Section \ref{section-pos-def}.
For this we need additional properties for  CAS, which we discuss in the next section.

Before doing this, we study the linear operators $T_\alpha$ for multiplicative functions $\alpha\in C(D)$, i.e.,
 $\alpha(h*l)=\alpha(h)\cdot\alpha(l)$ for all $h,l\in D$

\begin{lemma}\label{Tmult}
Let $(X,D,K)$ be a CAS, $g\in C_c(X)$, and let  $\alpha\in C(D)$ be multiplicative. 
Then the function $\phi:=T_\alpha g\in C(X)$ satisfies $T_l\phi(x)=\alpha(\bar l)\cdot\phi(x)$ 
for all $l\in D$ and $x\in X$.
\end{lemma}

\begin{proof}  Let $\omega_D$ be a left Haar measure of $(D,*)$ as above. 
Fix $l\in D$, $x\in X$ and consider the function $g_x(w):=\int_X g(y)\> K_w(x,dy)$ ($w\in D$).
Then $g_x\in C_c(D)$ by the considerations above, and, by
 Lemma 5.5G of \cite{J},
\begin{align}
\int_D\int_D g_x(h)\> d(\delta_l*\delta_w)(h) \> \alpha(w) \> d\omega_D(w)&=
\int_D g_x(l*w) \alpha(w) \> d\omega_D(w)
\notag\\
&
=\int_D g_x(w)\alpha(\bar l*w)\> d\omega_D(w)
\notag\\
&=\alpha(\bar l)\int_D g_x(w)\alpha(w) \> d\omega_D(w) \quad =\quad \alpha(\bar l)T_\alpha g(x).
\notag\end{align}
Hence,
 \begin{align}
T_l(T_{\alpha}g)(x)&= \int_D\int_X\int_X g(y)\> K_h(z,dy)\> K_l(x,dz)\> \alpha(h) \> d\omega_D(h)\notag\\
&= \int_D\int_D\int_X g(y)\> K_h(x,dy)  \> d(\delta_l*\delta_w)(h)\> \alpha(w)\> d\omega_D(w) \notag\\
&=\int_D\int_D g_x(h)\> d(\delta_l*\delta_w)(h) \> \alpha(w) \> d\omega_D(w)
\quad =\quad \alpha(\bar l)T_\alpha g(x)
\notag
\end{align}
as claimed.
\end{proof}

Lemma \ref{Tmult} can be applied to the uniqueness of the adjoint  measure $\omega_X\in M^+(X)$
at least in the compact case:

\begin{lemma}\label{unique-invar-measure} Let $(X,D,K)$ be a compact CAS. Then:
\begin{enumerate}
\item[\rm{(1)}] If $g\in C(X)$ satisfies $T_hg=g$ for all $h\in D$, then $g$ is constant.
\item[\rm{(2)}] For all  $g\in C(X)$, $T_{\bf 1}g$ is constant, and there is a unique measure
$\omega\in M_b^+(X)$ with $\omega(A)=\int_D K_h(x,A)\> d\omega_D(h)$ for all Borel sets $A\subset X$, 
where the right hand side is independent of $x\in X$.
\item[\rm{(3)}] The measure $\omega$ from (2) 
is equal to the adjoint measure $\omega_X\in M^+(X)$ from \ref{central-def}(4)
 up to a positive multiplicative constant. In particular $\omega_X$ is unique up to a
 positive multiplicative constant.
\end{enumerate}
\end{lemma}

\begin{proof} For (1) assume without loss of generality that $g$ is $\mathbb R$-valued with 
 $T_hg=g$ for all $h\in D$. Take $x_0\in X$ such that $g(x_0)$ is maximal. 
As $g(x_0)=\int_X g(y)\> K_h(x_0,dy)$ we see that $g(y)=g(x_0)$ for all $y\in supp\> K_h(x_0,.)$. 
As this holds for all $h\in D$, it follows that $g$ is constant.

For (2) we conclude from  Lemma \ref{Tmult} for $\alpha={\bf 1}$ that for $l\in D$, $T_l(T_{\bf 1}g)=T_{\bf 1}g(x)$.
Hence, by part (1),  $T_{\bf 1}g$ is constant. As
$T_{\bf 1}g(x)=\int_D\int_X g(y)\> K_w(x,dy) \> d\omega_D(w)$, we see from the representation theorem of Riesz that 
there exists  a unique measure
$\omega\in M_b^+(X)$ with $T_{\bf 1}g(x)=\int_X g\> d\omega$ for all $g\in C(X)$ and $x\in X$.

For (3) we assume without loss of generality that  $\omega_X\in M^1(X)$ and $\omega_D\in M^1(D)$. We 
 use the invariance condition (\ref{invariant-omega}) for $\omega_X$ and obtain for $g\in C(X)$ and each $x\in X$ that
$$ \int_X g\> d\omega= T_{\bf 1}g(x)= \int_X T_{\bf 1}g\> d\omega_X=
\int_D\int_X T_hg(x) \> d\omega_X(x) \> d\omega_D(h)= \int_X g d\omega_X.$$
This proves $\omega_X=\omega$ and the claim.
\end{proof}

\section{Strong continuous association schemes}

Let  $(X,D,K)$ be a CAS with associated hypergroup $(D,*)$.
We  first recapitulate
the translations $f_h(l):=f(l*h):=\int_Df \> d(\delta_l*\delta_h)$ of functions  $f\in C(D)$ for
 $h,l\in D$  as 
well as the projection maps $\pi_x:X\to D$ for $x\in X$ from the preceding section.

\begin{definition}\label{strong-central-def}
Let  $(X,D,K)$ be a CAS.
\begin{enumerate}
\item[\rm{(1)}] We say that  $(X,D,K)$ has the translation property (T1) 
if for all $h\in D$, $x\in X$, and $f\in C_c(D)$,
$$f_h\circ\pi_x=T_h(f\circ\pi_x).$$
\item[\rm{(2)}] We say that  $(X,D,K)$ has the translation property (T2), 
if for all $f\in C_c(D)$, $T_f=T^{f\circ\pi}$,
 where we assume that the invariant measure $\omega_X$ and the left Haar measure $\omega_D$ of $(D,*)$
are chosen with suitable normalizations.
\item[\rm{(3)}] We say that  $(X,D,K)$ is strong, if (T1) and (T2) hold.
 \end{enumerate}
\end{definition}

We shall prove below that in the discrete case and in the commutative case, property (T2) implies (T1);
see Theorems \ref{discont-asso-to-asso} and \ref{commT2T1}. This indicates that generally, (T2)
 seems to be the stronger condition. Otherwise we do not know further relations between these conditions.
In several sections below we  present
examples where (T1) and (T2) do not hold. On the other hand, there are several
standard classes of  strong  CAS; here is the first one:

\begin{proposition}\label{group-to-strong-cont-asso-scheme}
Let $H$ be a compact subgroup of a locally compact unimodular 
group $G$. Then the 
associated unimodular  CAS $(X:=G/H,D:=G//H,K)$ as in Proposition 
\ref{group-to-cont-asso-scheme} is strong.
\end{proposition}

\begin{proof}
Let $x,y,h\in G$, and $f\in C_c(G//H)$. The proof of Proposition 
\ref{group-to-cont-asso-scheme} yields 
$$(f_{HhH}\circ\pi_{xH})(yH)=\int_{G//H} f \> d(\delta_{Hx^{-1}yH}*\delta_{HhH})=
\int_H f(Hx^{-1}ywhH)\> d\omega_H(w)$$
and 
\begin{align}
T_{HhH}(f\circ\pi_{xH})(yH)&=\int_X f(Hx^{-1}zH)\> K_{HhH}(y, d(zH))) \notag\\
&=
\int_H f(Hx^{-1}ywhH)\> d\omega_H(w).
 \notag \end{align}
This proves (T1). Moreover, for  $f\in C_c(G//H)$, $g\in C_c(G/H)$, and $x\in G$ 
we have with the notations of 
 Proposition \ref{group-to-cont-asso-scheme} that
\begin{align}
T_f g(xH)&=\int_D\int_X g(yH) K_{HhH}(xH, yH) \> f(HhH)\> d\omega_{G//H}(HhH) \notag\\
&=\int_G\int_X g(yH) \> d(p_G(\delta_x*\omega_H*\delta_h*\omega_H))(yH) \> f(HhH)\>  d\omega_{G}(h) \notag\\
&=\int_H\int_G g(xrhH) \> f(HhH)\>  d\omega_{G}(h)\> d\omega_H(r) \notag\\
&=\int_H\int_G g(xhH)\> f(Hr^{-1}hH)\>  d\omega_{G}(h)\> d\omega_H(r)\notag\\
&=\int_G g(xhH)\> f(HhH)\>  d\omega_{G}(h)\notag\\
&=\int_G f(Hx^{-1}yH) \> g(yH) \> d\omega_{G}(y)\notag\\
&=\int_D f(\pi(xH,yH))\> g(yH)\> d\omega_{G/H}(yH) \quad =\quad 
T^{f\circ\pi}g(xH)
 \notag
\end{align}
which proves (T2).
\end{proof}

Here is a second standard class of strong CAS.

\begin{proposition}\label{asso-to-strong-cont-asso-scheme}
Let $(X,D, (R_i)_{i\in D})$ be an unimodular association scheme. Then the associated unimodular discrete
CAS
$(X,D,K)$  of Corollary \ref{discrete-ass-to-cont-asso-scheme}
is strong.
\end{proposition}

\begin{proof}
By linearity, it suffices to check (T1) for characteristic functions $f=\bf 1_{\{r\}}$ with $r\in D$.
For $h\in D$ and $x,y\in X$ we obtain from the axioms and basic properties of an  association scheme and 
the definition of the kernels $K_h$ that
\begin{align}
T_h({\bf 1}_{\{r\}}\circ\pi_x)(y)&= \int_X {\bf 1}_{\{r\}}(\pi(x,z))\> K_h(y,dz)= K_h(y, \{z\in X:\> \pi(x,z)=r\})
\notag\\
&=\frac{1}{\omega_h}|\{z\in X:\> \pi(x,z)=r, \> \pi(y,z)=h\}|
\notag\\
&=\frac{1}{\omega_h}|\{z\in X:\> \pi(z,x)=\bar r, \> \pi(y,z)=h\}|
\notag\\
&= \frac{1}{\omega_h}p_{h,\bar r}^{\pi(y,x)}=\frac{1}{\omega_h} p_{r,\bar h}^{\pi(x,y)}.
 \notag
\end{align}
On the other hand, we see from Proposition 3.8 and Lemma 3.5(4) of \cite{Voit16} that
\begin{align}
(f_h\circ \pi_x)(y)&=\int_D {\bf 1}_{\{r\}} \> d(\delta_{\pi(x,y)}*\delta_h)= (\delta_{\pi(x,y)}*\delta_h)(\{r\})
\notag\\
&= \frac{\omega_r}{\omega_h \omega_{\pi(x,y)}}  p_{\pi(x,y),h}^r
=\frac{1}{\omega_h} p_{r,\bar h}^{\pi(x,y)},
 \notag
\end{align}
which completes the proof of (T1). For (T2) we again use linearity and check (T2) for  $f=\bf 1_{\{r\}}$ and 
$g=\bf 1_{\{z\}}$
 with $r\in D$ and $z\in X$. Let $x\in X$. With the Kronecker-$\delta$ we  obtain
$$T_fg(x)= \omega_r\cdot\int_X g(y)\> K_r(x,dy)= \delta_{r,\pi(x,z)} =T^{f\circ \pi}g(x)$$
which yields the claim (T2).
\end{proof}

 Proposition \ref{asso-to-strong-cont-asso-scheme} has the following converse statement:

\begin{theorem}\label{discont-asso-to-asso}
Let $(X,D,K)$ be a discrete unimodular  CAS with property (T2). 
Then there is an unimodular association scheme $(X,D, (R_r)_{r\in D})$  such that  $(X,D,K)$ is the associated 
CAS according to Corollary \ref{discrete-ass-to-cont-asso-scheme}. In particular,
 for discrete unimodular  CAS,  (T2) implies (T1).
\end{theorem}

\begin{proof} Assume that the measure $\omega_X\in M^+(X)$ and 
the Haar measure $\omega_D$ are normalized such that (T2) holds. Let $r\in D$ and $x,z\in X$ and put
$f=\bf 1_{\{r\}}$ and $g=\bf 1_{\{z\}}$. Then, as in the proof of the preceding result, 
$T^{f\circ \pi}g(x) = \omega_X(\{z\})\>\delta_{r,\pi(x,z)}$ and 
\begin{equation}
T_fg(x)=  K_r(x,\{z\})\> \omega_D(\{r\})=  K_r(x,\{z\})\> \omega_D(\{r\})\>\delta_{r,\pi(x,z)}.
\end{equation}
Hence, by (T2),
\begin{equation}\label{special-form-kernel}
K_r(x,\{z\}) = \frac{\omega_X(\{z\})}{ \omega_D(\{r\})}\>\delta_{r,\pi(x,z)} \quad\quad(x,z\in X, \> r\in D).
\end{equation}
This in particular shows that for $r\in D$ and $x\in X$,
\begin{equation}\label{special-image-m}
 \omega_D(\{r\})=\sum_{z\in X:\> \pi(x,z)=r} \omega_X(\{z\})
\end{equation}
 and, as $K_r(x,.)$ is a probability measure,
 \begin{equation}\label{cond-critequality}
(K_r\circ K_{\bar r})(x,\{x\})= \sum_{z\in X:\> \pi(x,z)=r}\frac{\omega_X(\{z\})}{ \omega_D(\{r\})}
 \cdot \frac{\omega_X(\{x\})}{ \omega_D(\{\bar r\})}\> \delta_{\bar r,\pi(z,x)}
= \frac{\omega_X(\{x\})}{ \omega_D(\{\bar r\})}.
\end{equation}
As by the definition of a discrete CAS $K_r\circ K_{\bar r}$ is a finite convex combination of the $K_s$ ($s\in D$)
where the identity kernel $K_e$ appears with a positive coefficient, we conclude from (\ref{cond-critequality})
that $\omega_X(\{x\})>0$ is independent of $x\in X$. Therefore, after renormalization of  $\omega_X$ and
 $\omega_D$, we may assume that 
 $\omega_X$ is the counting measure.  We then see from (\ref{special-image-m}) that
$ \omega_D(\{r\})=|supp\> K_r(x,.)|$ for $r\in D$ and all $x\in X$. Moreover, by (\ref{special-form-kernel}), 
 \begin{equation}\label{special-form-kernel2}
K_r(x,\{z\})=\frac{\delta_{r,\pi(x,z)} }{\omega_D(\{r\})} \quad\quad(r\in D, \> x,z\in X).
\end{equation}
We now define the partition $(R_r)_{r\in D}$  of $X\times X$ via $R_r:=\{(x,y):\> y\in supp\> K_r(x,.)\}$.
This is a partition by the partition property of a CAS for which clearly property (1) of \ref{def-discrete-asso}
holds. Moreover,  \ref{def-discrete-asso}(2) follows from (\ref{special-form-kernel2}) and 
$\omega(\{\bar r\})=\omega(\{ r\})$. We finally check \ref{def-discrete-asso}(3).
 For this let $i,j,k\in D$ and $x,y\in X$ with $\pi(x,y)=k$. Then
$$K_i\circ K_j(x,\{y\})= \sum_{z\in X: \pi(x,z)=i, \> \pi(z,y)=j} \frac{1}{\omega_D(\{i\})\omega_D(\{j\})}=
\frac{|\{z\in X: \pi(x,z)=i, \> \pi(z,y)=j\}|}{\omega_D(\{i\})\omega_D(\{j\})}$$
and
$$K_i\circ K_j(x,\{y\})= (\delta_i*\delta_j)(\{k\})\cdot
 K_k(x,\{y\})=\frac{(\delta_i*\delta_j)(\{k\})}{\omega_D(\{k\}}.$$
A comparison of both formulas shows that $|\{z\in X: \pi(x,z)=i, \> \pi(z,y)=j\}|$ depends only on 
$\pi(x,y)=k$ and not on the choice of $x,y$ as claimed.

In summary, we see that $(X,D, (R_r)_{r\in D})$ is an association scheme with $(X,D,K)$ as
 associated CAS  by  (\ref{special-form-kernel2}). 
The last statement of the theorem follows from Proposition \ref{asso-to-strong-cont-asso-scheme}.
\end{proof}

In summary, in the unimodular case,
 strong discrete CAS are precisely classical  association schemes.  Moreover, 
 discrete commutative CAS may be seen as generalizations of generalized  association schemes.

Theorem \ref{discont-asso-to-asso} suggests that  in general (T2) implies (T1).
Unfortunately, we do not see any approach for the proof of this conjecture in the non-discrete case.

We next rewrite (T2) as follows:

\begin{lemma}\label{t2-reformulation}
Let $(X,D,K)$ be an unimodular  CAS with the Haar measure $\omega_D$ of $(D,*)$ and the adjoint measure 
$\omega_X\in M^+(X)$. Then (T2) holds if and only if 
 $\omega_X(A)=\int_D K_h(x,A)\> d\omega_D(h)$ for all Borel sets $A\subset X$ and  $x\in X$.

In particular, for each unimodular CAS with (T2), $\omega_X$ is unique up to a positive  constant.
 \end{lemma}

\begin{proof} Assume first that (T2) holds. It can be easily seen (see Lemma \ref{proj-haar-m-T2} below)
that then for all $g\in C_c(X)$ and $x\in X$, $T^{{\bf 1}\circ\pi}g = T_{\bf 1}g$ and thus
$$\int_D\int_X  g(y)  \> K_h(x,dy)\>d\omega_D(h)= T_{\bf 1}g(x)=T^{{\bf 1}\circ\pi}g(x)
=\int_X  g(y) \> d\omega_X(y).$$ 
This shows that  $\omega_X=\int_D K_h(x,.)\> d\omega_D(h)$ for  $x\in X$.

Conversely, this representation of $\omega_X$  shows for $f\in C_c(D)$, $g\in C_c(X)$, and $x\in X$ that
\begin{align}
T^{f\circ\pi}g(x)&= \int_X g(y) \cdot f(\pi(x,y))\> d\omega_X(y) 
=\int_D\int_X g(y) \cdot f(\pi(x,y))\> K_h(x,dy)\>  d\omega_D(h) \notag\\
&= \int_D \int_X g(y)\cdot f(h) \> K_h(x,dy)\> d\omega_D(h)
=T_fg(x).
\notag\end{align}
 Hence, (T2) holds. The uniqueness assertion is  clear.
\end{proof}

Lemma \ref{t2-reformulation} and Lemma \ref{unique-invar-measure} now lead to:

\begin{proposition}\label{compact-t2} 
Each compact CAS has property (T2).
\end{proposition}

If we combine Proposition  \ref{compact-t2}  with Theorem \ref{discont-asso-to-asso} and
 Proposition \ref{discrete-gen-ass-to-cont-asso-scheme}, we obtain the following classification of finite CAS:

\begin{theorem}\label{classification-finite} 
Each finite  CAS $(X,D,K)$ is  associated with an
 association scheme  according to Corollary \ref{discrete-ass-to-cont-asso-scheme}. 
In particular, each finite generalized association scheme is in fact an  association scheme.
\end{theorem}

We notice that this classification does not hold in the infinite case. 
Examples are given in \cite{Voit17} and below, e.g.~in Section 10.

We now return to  (T1) and (T2) and 
 study CAS with these properties:

\begin{lemma}\label{proj-haar-m-T1}
 Let $(X,D,K)$ be a  CAS with  (T1). Then:
\begin{enumerate}
\item[\rm{(1)}] For all $x\in X$, the push forward 
$\pi_x(\omega_X)\in M^+(D)$  is a right Haar measure of $(D,*)$. 
\item[\rm{(2)}]  For all $\mu\in M_b(D)$, $f\in C_c(D)$ and $x\in X$,
$T_\mu(f\circ \pi_x)=(f*\mu^-)\circ \pi_x.$
\item[\rm{(3)}]   For all $\phi\in C(D)$, $f\in C_c(D)$ and $x\in X$,
$T_\phi(f\circ \pi_x)=(f*\phi^-)\circ \pi_x.$
\end{enumerate}
\end{lemma}

\begin{proof}
\begin{enumerate}
\item[\rm{(1)}] For all $h\in D$ and $f\in C_c(D)$ we obtain from \ref{facts-basic-gen-ass}(3) that, as claimed,
$$(\pi_x(\omega_X))(f_h)= \omega_X(f_h\circ\pi_x)=
\omega_X(T_h(f\circ\pi_x))=\omega_X(f\circ\pi_x)=(\pi_x(\omega_X))(f).$$
\item[\rm{(2)}] and (3) follow simply by integration of the equation in Definition \ref{strong-central-def}.
\end{enumerate}
\end{proof}

\begin{lemma}\label{proj-haar-m-T2}
 Let $(X,D,K)$ be a  CAS  with (T2). Then:
\begin{enumerate}
\item[\rm{(1)}] For all $x\in X$, the push forward 
$\pi_x(\omega_X)\in M^+(D)$ is a right Haar measure of $(D,*)$. 
\item[\rm{(2)}] For all $f\in C(D)$ and  $g\in C_c(X)$, $T_fg=T^{f\circ\pi}g\in C(X)$.
\item[\rm{(3)}] For all $f\in C_c(D)$ and  $g\in C(X)$, $T_fg=T^{f\circ\pi}g\in C(X)$.
\end{enumerate}
\end{lemma}

\begin{proof} (2) is clear, and (3) follows from Lemma \ref{compact-union} similar to the proof of Lemma 
\ref{compact-preserved}. 
For the proof of (1) we use (3) with $g\equiv 1$ and $f\in C_c(D)$. Hence, for $x\in X$,
\begin{align}
\int_X f(\pi(x,y))\> d\omega_X(y)&= T^{f\circ\pi}g(x)=T_fg(x)\notag\\
&=\int_D\int_X 1\> K_h(x,dy)\> f(h)\> d\omega_D(h)
= \int_D f\> d\omega_D(h)
\notag\end{align}
which proves the claim.
\end{proof}

\begin{remark}
There exist commutative CAS without  (T1) and (T2). For this consider
 the discrete  generalized association schemes associated with homogeneous trees of
 \cite{Voit17} with the parameter $c\ne 1$ there. As shown in  Remark 2.4 of \cite{Voit17}, 
there the measure $\omega_X$ is unique up to a positive multiplicative constant for which the
push forward statements of Lemmas 
\ref{proj-haar-m-T1}(1) and \ref{proj-haar-m-T2}(1) are not
 correct, i.e., (T1) and (T2) do not hold there; see also Section 10.
\end{remark}

\begin{lemma}\label{proj-haar-m-T12}
 Let $(X,D,K)$ be a  strong unimodular CAS. 
Let  $f\in C_c(D)$ and $f\in C(D)$, or $f\in C(D)$ and $g\in C_c(D)$, or $f,g\in L^2(D,\omega_D)$.
Then, for all $x,z\in X$,
$$\int_X f(\pi(x,y))\> \overline{g(\pi(z,y))}\> d\omega_X(y)
= \int_D f(h)\> \overline{g(\pi(z,x)*h)} \>d\omega_D(h).$$
\end{lemma}

\begin{proof}
Let  $f\in C_c(D)$ and $f\in C(D)$. Then, by (T2) and Lemma \ref{proj-haar-m-T1}(3),
\begin{align}
\int_X& f(\pi(x,y))\> \overline{g(\pi(z,y))}\> d\omega_X(y) = T^{f\circ\pi}(\overline{g\circ \pi_z})(x)
= T_f(\overline{g\circ \pi_z})(x)
\notag\\
&= ((f^-*\bar g)\circ  \pi_z)(x)= (f^-*\bar g)(\pi(z,x))= \int_D f(h) \> \overline{g(\pi(z,x)*h)} \>d\omega_D(h)
\notag\end{align}
as claimed. The same computation works for  $f\in C(D)$ and $g\in C_c(D)$ as well as for
 $f,g\in L^2(D,\omega_D)$ by density. Notice here that due to  Lemma \ref{proj-haar-m-T1}(1), for all $z\in X$
the map $f\mapsto f\circ \pi_z$ is an $ L^2$-isometry from  $L^2(D,\omega_D)$ into  $L^2(X,\omega_X)$.
\end{proof}

%For simplicity we restrict our attention to the commutative case from now on.
We next present some orthogonality result which is well-known
 in the group case.

\begin{corollary}
 Let $(X,D,K)$ be a  compact, commutative strong CAS.
 Then for $\alpha,\beta\in (D,*)^\wedge$ and $x,z\in X$,
$$\int_X \alpha(\pi(x,y))\> \overline{\beta(\pi(z,y))}\> d\omega_X(y)= 
\delta_{\alpha,\beta}\cdot  \overline{\alpha(\pi(z,x))}\cdot \|\alpha\|_{2,\omega_D}^2.$$
\end{corollary}

\begin{proof} Lemma \ref{proj-haar-m-T12}   yields
\begin{align}\int_X& \alpha(\pi(x,y))\> \overline{\beta(\pi(z,y))}\> d\omega_X(y)=
\int_D \alpha(h)\> \overline{\beta(\pi(z,x)*h)} \>d\omega_D(h)
\notag\\
&=\int_D \alpha(h)\>\overline{\beta(h)} \>d\omega_D(h)\cdot \overline{\beta(\pi(z,x))}.
\notag\end{align}
 As the characters of the compact commutative hypergroup $(D,*)$ form an orthogonal basis
of $L^2(D,\omega_D)$, the proof is complete.
\end{proof}

We finally remark that for given spaces $X,D$ and given projection $\pi:X\times X\to D$, there is at most one
CAS with property (T2), i.e., (T2) is a quite strong condition:

\begin{proposition}\label{uniqueness-complete-t2}
Let  $(X,D, K)$ and $(X,D,\tilde K)$ be  CAS with  property (T2)
with the same $X,D$, $\pi$. Then  $(X,D, K)=(X,D,\tilde K)$  and  $(D,*)=(D,\tilde *)$.
\end{proposition} 

\begin{proof} (T2) implies that for all $f\in C_c(D)$ and $g\in C_c(X)$
\begin{equation}\label{eq-tilde-tf}
T_fg= T^{f\circ\pi}g=\tilde T_fg
\end{equation}
with $\tilde T_fg$ as operator associated with the kernels $\tilde K_h$.
 As for $x\in X$ and $g\in C_c(X)$ the map 
$D\to \mathbb C, \> h\mapsto \int_X g(y)\> K_h(x,dy)$ is continuous, (\ref{eq-tilde-tf})
 implies by a  limit that $T_hg=\tilde T_hg$ for all $g\in C_c(D)$ and $h\in D$.
 This also readily shows that 
$K_h=\tilde K_h$ for all $h$. Fact \ref{facts-basic-gen-ass}(1) finally proves $(D,*)=(D,\tilde *)$.
\end{proof}

Propositions \ref{uniqueness-complete-t2} and \ref{compact-t2} yield:

\begin{corollary}\label{equal-compact}
Let  $(X,D, K)$ and $(X,D,\tilde K)$ be  compact CAS 
with the same $X,D$, $\pi$. Then  $(X,D, K)=(X,D,\tilde K)$  and  $(D,*)=(D,\tilde *)$.
\end{corollary}

Variants of  \ref{uniqueness-complete-t2} and \ref{equal-compact} will be given in Section 7.

\section{Positive definite functions}\label{section-pos-def}

In this section we study several concepts of positive definiteness on CAS.
 We restrict our attention to the commutative case
for simplicity even if some results remain valid in a slightly more general setting.
Therefore, $(X,D,K)$ will always  be a  commutative CAS.

\begin{definition}
\begin{enumerate}
\item[\rm{(1)}] Let $A:C_c(X)\to C(X)$ be a linear operator. $A$ is called positive definite if 
$\langle Ag, g\rangle_{X}\in [0,\infty[$ for all $g\in C_c(X)$.
\item[\rm{(2)}] A continuous function $F:X\times X\to\mathbb C$ is called positive definite, if for all 
 $n\in\mathbb N$, $x_1,\ldots, x_n\in X$ and
  $c_1,\ldots,c_n\in\mathbb C$, $\sum_{k,l=1}^n c_k\bar c_l \> F(x_k,x_l)\ge0$.
 \end{enumerate}
\end{definition}

Both concepts are closely related:

\begin{lemma}\label{crit-pd}
 A continuous function $F:X\times X\to\mathbb C$ is positive definite if and only
 if the linear operator $T^F:C_c(X)\to C(X)$ is  positive definite.
\end{lemma}

\begin{proof} This follows from standard density arguments similar to corresponding results
 for hypergroups; see e.g.~Lemma 4.1.4 of \cite{BH}. \end{proof}

As the pointwise products of positive semidefinite matrices are again positive 
 semidefinite (see e.g.~Lemma 3.2 of \cite{BF}), we have the following well known result:

\begin{lemma}\label{prod-pos}
If  $F,G:X\times X\to\mathbb C$ are positive definite, then the pointwise product
$F\cdot G:X\times X\to\mathbb C$ is also positive definite.
\end{lemma}

We now study for which $f\in C_b(D)$ the operators $T_f$ are positive definite. 
The following more or less obvious result will be needed later on.

\begin{lemma}\label{approx-pd}
For a function $f\in C_b(D)$, the operator $T_f$
 is positive definite if and only if for each step function $g=\sum_{i=1}^n c_i{\bf 1}_{A_i}:X\to\mathbb C$
with $n\in\mathbb N$, $c_1,\ldots, c_n\in\mathbb C$, and disjoint, relatively compact Borel sets 
$A_1,\ldots,A_n\subset X$, the inequality $\langle T_fg, g\rangle_{X}\in [0,\infty[$ holds.
\end{lemma}

\begin{proof} For the only-if-part, we notice that each step function $g$ as required in the lemma is the
pointwise limit of functions in $C_c(X)$ whose supports are contained in some fixed compactum in $X$.
 The result then follows from dominated convergence. The if-part follows by the same arguments.
 \end{proof}

We now collect some relations between positive definite functions on $D$ and positive definiteness on $X$.

\begin{lemma}\label{pd-2-cont}
Let  $f\in C_c(D)$. Then,  
$f*f^*$ is positive definite on $D$, and
 $T_{f*f^*}$ is positive definite.
\end{lemma}

\begin{proof} The first  statement is well-known; see \ref{facts-hy}(2). The second one 
 is clear as $T_{f*f^*}= T_{f}T_{f}^*$ by \ref{alg-hom}.
\end{proof}

\begin{corollary}\label{pd-2a-cont}
For each character $\alpha\in(D,*)^\wedge$ in the
 support $S$ of the Plancherel measure of $(D,*)$, the operator $T_\alpha$ is  positive definite.
Moreover, if $f\in C_b(D)$ is positive definite on $(D,*)$ such that $f$ 
has the form $f=\check\mu$ for some $\mu\in M_b^+(S)$, then $T_f$ is  positive definite.
\end{corollary}

\begin{proof}
By  Fact \ref{facts-hy}(3), each $\alpha\in S$ is a locally uniform
 limit of functions of the form $f*f^*$ with  $f\in C_c(D)$. 
It follows from the axioms of a continuous association scheme and the
 definition of $T_{\alpha}$ that for each $g\in C_c(X)$,  $T_\alpha g$ is a locally uniform limit of 
 $T_{f*f^*} g$. Lemma \ref{pd-2-cont} thus implies that  $T_\alpha$ is  positive definite. 
The second statement follows in the same way.
\end{proof}

If $(X,D,K)$ has property (T2), then the preceding results can be rewritten:

\begin{corollary}\label{pd-2-cont-special}
 Let $(X,D,K)$ be a commutative CAS with (T2). Then:
\begin{enumerate}
\item[\rm{(1)}] For each  $f\in C_c(D)$, $(f*f^*) \circ\pi:X\times X\to\mathbb C$ is  positive definite.
\item[\rm{(2)}] For each character $\alpha\in S\subset \hat D$ in the
 support of the Plancherel measure of $(D,*)$, $\alpha\circ\pi:X\times X\to\mathbb C$ is positive definite.
\end{enumerate}
\end{corollary}

\begin{proof} Part (1) follows from Lemma \ref{pd-2-cont}, property (T2), and Lemma \ref{crit-pd}.

For part (2) we again use that  $\alpha\in S$ is a locally uniform
 limit of functions of the form $f*f^*$ with  $f\in C_c(D)$. Hence, by part (1),  $\alpha\circ\pi$
 is a locally uniform limit of positive definite functions
 on $X\times X$ and thus also positive definite.
\end{proof}

We now turn to a kind of converse statement of Lemma \ref{pd-2-cont}  and Corollary \ref{pd-2a-cont}:

\begin{lemma}\label{pd-1-cont} Let $f\in C(D)$ such that $T_f$ is positive definite. Then: 
\begin{enumerate}
\item[\rm{(1)}] $f(e)\in  [0,\infty[$; 
\item[\rm{(2)}] $f$ is positive definite on $(D,*)$.
\end{enumerate}
\end{lemma}

\begin{proof} For part (1) assume that $f(e)\in\mathbb C\setminus [0,\infty[$ holds, i.e.,
$\arg f(e)\in\mathbb R\setminus 2\pi\mathbb Z$ with a branch of the arg-function on
 $ \mathbb C\setminus \{0\}$ which is continuous in $f(e)$.
Now choose $\epsilon>0$ such that for all $z\in\mathbb C$ with $|z-f(e)|<\epsilon$ we have $z\ne0$, 
$\arg z\not\in  2\pi\mathbb Z$, and
$|\arg z-\arg f(e)|<1/2$ (or another  small positive constant). As $f$ is continuous, 
we find a neighborhood $W_e\subset D$ of $e$ with $|f(x)-f(e)|<\epsilon$ for $x\in W_e$.
We thus obtain that for all $\phi\in C(D)$ with values in $[0,\infty[$ and with $\phi(e)>0$,
\begin{equation}\label{bdgg-not-pos}
\int_{W_e} \phi(x)\> f(x)\> d\omega_D(x)\in\mathbb C\setminus [0,\infty[.
\end{equation}
On the other hand we now fix some $z\in X$. As $\pi:X\times X\to D$ is continuous with $\pi(z,z)=e$, 
we find a neighborhood $U_z\subset X$ of $z$ with $\pi(U_z,U_z)\subset W_e$.
Choose some $g\in C_c(X)$ with values in $[0,\infty[$ and with $g(z)>0$ and $supp\> g \subset U_z$.
Then $K_h(x,U_z)=0$ for all $x\in U_z$ and $h\in D\setminus W_e$. As $T_f$ is positive definite, we obtain
\begin{align}
0\le \langle Ag, g\rangle_{X} &= \int_D \int_X\int_X g(x)g(y) \> K_h(x,dy) \> d\omega_X(x) \> \cdot\>  f(h)\> d\omega_D(h)
\notag\\
 &= \int_{W_e} \int_{U_z}\int_{U_z} g(x)g(y) \> K_h(x,dy) \> d\omega_X(x) \> \cdot\>  f(h)\> d\omega_D(h)
\notag\\
 &=: \int_{W_e} \phi(x)\> f(x)\> d\omega_D(x)
\notag\end{align}
where $\phi$ is continuous with with values in $[0,\infty[$ and  with $\phi(e)>0$. This contradicts (\ref{bdgg-not-pos})
and completes the proof of part (1).

For (2) consider any $\phi\in C_c(D)$ and $g\in C_c(X)$. Then, by \ref{op-t-mu}(2), $T_\phi g\in   C_c(X)$,
and by our preceding considerations,
$$\langle T_{\phi^**f*\phi}    g, g\rangle_{X}= \langle T_{\phi}^*T_fT_\phi g,g\rangle_{X}=
\langle T_fT_\phi g,T_\phi g\rangle_{X}\ge0.$$
This shows that  $T_{\phi^**f*\phi}$ is positive definite, 
and we obtain from  a standard computation for hypergroups and 
part (1)  that
 $$\int_D\int_D f(h_1*\bar h_2) \cdot 
\phi(h_1)\cdot \overline{\phi(h_2)}\> d\omega_D(h_1) \> d\omega_D(h_2) =  \phi^**f*\phi(e)\in  [0,\infty[.$$
As this holds for all $\phi\in C_c(D)$, it follows from 
standard arguments for hypergroups (see Lemma 4.1.4 of \cite{BH})
 that  $f$ is positive definite on $(D,*)$.
\end{proof} 

Corollary \ref{pd-2-cont-special} and Lemmas \ref{prod-pos} and \ref{pd-1-cont} now lead to the
following  result. Is was given for  association schemes in  Theorem 4.6 of \cite{Voit16}.

\begin{theorem}\label{main-cont}
 Let $(D,*)$ be a commutative  hypergroup which is associated with
 some CAS $(X,D,K)$ with  (T2).
Then, for all   $\alpha,\beta\in S\subset\hat D$  in the support
 of the Plancherel measure,  $\alpha\cdot\beta$ 
is  positive definite  on $D$, and 
there is a unique probability 
measure $\delta_\alpha\hat*\delta_\beta\in M^1(\hat D)$ with 
$(\delta_\alpha\hat*\delta_\beta)^\vee=\alpha\cdot\beta$.
The support of
this measure is contained in $S$.

Furthermore, for all $\alpha\in S$, 
the unique positive character $\alpha_0$ in $S$ according to 
 \ref{facts-hy}(4)
is contained in the support of
$\delta_\alpha\hat*\delta_{\bar\alpha}$.
\end{theorem}

\begin{proof}
Corollary \ref{pd-2-cont-special}, property (T2), and Lemmas \ref{prod-pos} and \ref{pd-1-cont}
 show
that  $\alpha\cdot\beta$ 
is  positive definite on $(D,*)$. Bochner's theorem  \ref{facts-hy}(1) now leads
to the probability measure $\delta_\alpha\hat*\delta_\beta\in M^1(\hat D)$.
Furthermore, Proposition \ref{support-condition-pl} ensures that the support of
this measure is contained in $S$.
The assertion about the support of
$\delta_\alpha\hat*\delta_{\bar\alpha}$ follows from Theorem 2.1 of \cite{V3}.
\end{proof}

The methods of the proof of Theorem \ref{main-cont}
 can be used to prove the following equivalence of different concepts of  positive definiteness.

\begin{proposition} Let  $(X,D,K)$ be a  commutative CAS with property (T2)
 such that ${\bf 1}$ is contained in the support
$S$ of the Plancherel measure of the associated hypergroup $(D,*)$. Then for $\alpha\in (D,*)^\wedge$ the following
facts are equivalent:
\begin{enumerate}
\item[\rm{(1)}] $\alpha\in S$; 
\item[\rm{(2)}] The operator $T_\alpha$ is  positive definite;
\item[\rm{(3)}] $\alpha\circ\pi\in C_b(X\times X)$ is positive definite;
\item[\rm{(4)}] For each $\beta\in S$, the product $\alpha\cdot\beta$ is positive definite on $(D,*)$.
\end{enumerate}
\end{proposition}

\begin{proof} $(1)\Longrightarrow(2)$ follows from Corollary \ref{pd-2a-cont}, and
  $(2)\Longrightarrow(3)$ is a consequence of (T2) and Lemma \ref{crit-pd}.  $(3)\Longrightarrow(4)$ follows from 
 Lemma \ref{prod-pos} with the methods of the proof of Theorem \ref{main-cont}. Finally,  $(4)\Longrightarrow(1)$
is a consequence of  ${\bf 1}\in S$ and Corollary 7 of \cite{V2}.
\end{proof}

For compact CAS, Theorem \ref{main-cont} can be improved:

\begin{theorem}\label{main-cont-compact} 
Let $(D,*)$ be a compact commutative  hypergroup which is associated with
 some compact commutative CAS $(X,D,K)$. Then   (T2) holds by Proposition \ref{compact-t2}, and,
 with the dual convolution $\hat *$ of Theorem \ref{main-cont}, $(\hat D, \hat *)$
satisfies all hypergroup axioms possibly 
except for the condition that $supp(\delta_\alpha\hat *\delta_\beta)$ is compact
(i.e.~finite) for all $\alpha,\beta\in\hat D$.
\end{theorem}

\begin{proof} For  compact commutative hypergroups we have $S=\hat D$, $\hat D$ is discrete,
 and the unique
positive character in $S$ is the identity ${\bf 1}$; see e.g. \cite{BH}. Therefore, 
if we take ${\bf 1}$ as identity and 
complex conjugation as involution, almost all hypergroup
axioms  of $(\hat D,\hat *)$ follow from Theorem \ref{main-cont}. In fact, as $\hat D$ is discrete,
 the topological axioms hold automatically. Moreover, 
the bilinear, weakly continuous  extension of the dual convolution $\hat *$ from the set of point measures to
 $M_b(\hat D)$ is associative as the inverse Fourier transform is injective; see \cite{BH}.

We thus only have to check that for $\alpha\ne\beta\in \hat D$, the character
 ${\bf 1}$ is not contained in the support of 
$\delta_\alpha\hat *\delta_{\bar\beta}$. For this we recapitulate that
 for all $\gamma,\rho \in\hat D$, 
$\hat\gamma(\rho)=\int_D\gamma\bar\rho\>d\omega_D=\|\gamma\|_2^2 \delta_{\gamma,\rho}$ with the Kronecker-$\delta$.
Therefore, with 12.16 of \cite{J},
\begin{align} (\delta_\alpha\hat *\delta_{\bar\beta})(\{{\bf 1}\})&=
\int_{\hat D}  {\bf 1}_{\{{\bf 1}\}}\> d(\delta_\alpha\hat *\delta_{\bar\beta})= 
\int_{\hat D} \hat{\bf 1}\> d(\delta_\alpha\hat *\delta_{\bar\beta}) =\notag\\
&=\int_{ D} {\bf 1}\> (\delta_\alpha\hat *\delta_{\bar\beta})^\vee\> d\omega_D = 
\int_{ D} \alpha\bar\beta \> d\omega_D=
\|\gamma\|_2^2 \delta_{\alpha,\beta}=0.
\notag\end{align}
This completes the proof.
\end{proof}

In the finite case, Theorem \ref{main-cont-compact} is a follows; see also Theorem 4.7 of \cite{Voit16}:

\begin{corollary}\label{main-cont-finite} 
Let $(D,*)$ be a finite commutative  hypergroup which is associated with
 some finite commutative CAS $(X,D,K)$. Then  $\hat D$ carries a dual hypergroup structure.
\end{corollary}

We next turn to the problem whether
  the conclusions of  \ref{main-cont}, \ref{main-cont-compact}, and \ref{main-cont-finite} 
on positive dual convolutions also hold for commutative CAS without  (T2).
We here follow  Section 5 of \cite{Voit16} and assume that we have two commutative 
 CAS $(X,D,K)$ and  $(X,D,\tilde K)$ 
with the same spaces $X,D$ and the same projection $\pi:X\times X\to D$.
We denote the associated commutative hypergroups by  $(D,*)$ and $(D,\tilde*)$ and the supports
of the associated Plancherel measures by $S$ and $\tilde S$.
 Assume that $(X,D,K)$
 has property (T2), and that all characters $\tilde\alpha\in\tilde S$ of $(D,\tilde*)$ 
in the support of the Plancherel measure have the form
$$\tilde\alpha(h)=\int_S \alpha(h)\> d\mu(\alpha) \quad\quad\text{for all }\quad h\in D$$
for some  $\mu\in M^1(S)$.
 It was proved in Theorem 5.10 of \cite{Voit16} for commutative generalized association schemes 
 that then  $(D,\tilde *)$ also admits a positive dual convolution on 
$\tilde S$. We shall extend  this result to CAS
in Theorem \ref{cont-main-g} which has some unexpected consequences: It turns out that under
 some additional conditions like property (T2) for {\bf one} of the CAS and a support condition,  
the hypergroup structures  $(D,*)$ and $(D,\tilde*)$ are equal; 
compare this  with the assertions of  \ref{classification-finite}, 
\ref{uniqueness-complete-t2}, and \ref{equal-compact}.

\section{A comparison of different CAS on the same spaces}

As before, let  $(X,D, K)$ and $(X,D,\tilde K)$ be commutative CAS 
with the same spaces $X,D$ and the same projection $\pi:X\times X\to D$.
 Let again   $(D,*)$ and $(D,\tilde*)$
be the associated commutative hypergroups and  $S$ and $\tilde S$ the associated Plancherel measures
 respectively.
The following extension of Theorem  \ref{main-cont}
  is the main result of this section.

\begin{theorem}\label{cont-main-g} 
Assume that $(X,D, K)$ has property (T2) in the setting above. Then for all characters
$\tilde\alpha\in \tilde S$ and $\beta\in S$, the product $\tilde\alpha\cdot\beta$ is positive definite on 
$(D,\tilde *)$, and there is
a unique $\mu\in M^1(\tilde S)$ with 
$$\tilde\alpha(h)\beta(h)=\int_{\tilde S} \alpha(h)\> d\mu(\alpha) \quad\quad\text{for all }\quad h\in D.$$
\end{theorem}

In the discrete case, the proof is quite simple and similar to that of Theorem 5.10 of \cite{Voit16},
 while it will
be more involved  in the continuous case due to some approximation procedure.
To highlight the idea of the proof, we first give the proof in the 
discrete case:

 \begin{proof}[Proof of Theorem \ref{cont-main-g} in the discrete case]
Let $\tilde\alpha\in \tilde S$ and $\beta\in S$. Let $\tilde T_{\tilde\alpha}$ be the linear operator associated with
 $\tilde\alpha$ and the CAS $(X,D,\tilde K)$. 
 Then, by Corollary \ref{pd-2a-cont},
$\tilde T_{\tilde\alpha}$ is positive definite. Now let  $g\in C_c(X)$.
 Choose $x_1,\ldots, x_n\in X$ different 
with $supp\> g=\{x_1,\ldots, x_n\}$, and let $\tilde\omega_D$ be a Haar measure of
$(D,\tilde *)$,  $\tilde\omega_X\in M^+(X)$ the associated measure, and 
$\langle  .,.\rangle_{\tilde X}$ the associated scalar product on $L^2(X,\tilde\omega_X)$. 
The  positive definiteness of $\tilde T_{\tilde\alpha}$ and the 
properties of $supp\> K_h(x,.)$ for $h\in D$ and $x\in X$ imply that
\begin{align}\label{pd-test-al}
0 &\le \langle \tilde T_{\tilde\alpha} g,g\rangle_{\tilde X} \> =\>
\sum_{h\in D}\sum_{k,l=1}^n \overline{g(x_k)} \> g(x_l) \> \tilde K_h(x_k,\{x_l\})\> \tilde \omega_X(\{x_k\}) \cdot
\tilde\alpha(h)
\tilde\omega_D(\{h\}) \\
&=\sum_{k,l=1}^n \overline{g(x_k)} \> g(x_l) \> \tilde K_{\pi(x_k,x_l)}(x_k,\{x_l\})\> \tilde \omega_X(\{x_k\}) 
\tilde\alpha(\pi(x_k,x_l))\tilde\omega_D(\{\pi(x_k,x_l)\}).
\notag
\end{align}
On the other hand, as $(X,D,K)$ has property (T2), 
 $\beta\circ\pi:X\times X\to\mathbb C$ is positive definite, i.e., the matrix
 $(\beta(\pi(x_k,x_l)))_{k,l}$ is positive semidefinite. Eq.~(\ref{pd-test-al})
 and the fact that pointwise products of
 positive semidefinite matrices are positive semidefinite yield that
\begin{align}
\sum_{k,l=1}^n \overline{g(x_k)} \> g(x_l) &\> \tilde K_{\pi(x_k,x_l)}(x_k,\{x_l\})\> \tilde\omega_X(\{x_k\})\cdot \notag\\
&\cdot\tilde\alpha(\pi(x_k,x_l))\beta(\pi(x_k,x_l))\tilde\omega_D(\{\pi(x_k,x_l)\})\ge0.
\notag
\end{align}
As in the computation of  Eq.~(\ref{pd-test-al}), we obtain that
 $\langle \tilde T_{\tilde\alpha\cdot\beta} g,g\rangle_{\tilde X} \ge 0$,
 i.e., $\tilde T_{\tilde\alpha\cdot\beta}$ is positive definite.
 Hence, by Lemma \ref{pd-1-cont}, $\tilde\alpha\beta$ is positive definite
on $(D,\tilde *)$. 
Finally, the support condition follows from Proposition \ref{support-condition-pl}.
\end{proof}

We now turn to the  the general case by using Lemma \ref{approx-pd}.

\begin{proof}[Proof of Theorem \ref{cont-main-g} in the general case]
We keep the notations of the  discrete case. Let $\tilde\alpha\in \tilde S$ and $\beta\in S$. 
Consider some   step function $g=\sum_{i=1}^n c_i{\bf 1}_{A_i}:X\to\mathbb C$
with $n\in\mathbb N$, $c_1,\ldots, c_n\in\mathbb C$, and disjoint, relatively compact, non-empty Borel sets 
$A_1,\ldots,A_n\subset X$ as in Lemma \ref{approx-pd}. Let $\epsilon>0$ be a small constant.

As $\beta\circ\pi:X\times X\to\mathbb C$ is continuous and thus uniformly continuous on the compactum 
$supp \> g\times supp\> g\subset X\times X$, we may decompose the sets  $A_1,\ldots,A_n$ into finitely many,
 disjoint,  non-empty  Borel sets such that, after denoting these finitely many sets again by   $A_1,\ldots,A_n$ 
 we have the following additional property:
\begin{equation}\label{distanz-ab}
\text{For all}\quad i,j=1,\ldots,n, \> u,v\in A_i,\> x,y\in A_j: \quad |\beta(\pi(u,x))-\beta(\pi(v,y))|
\le\epsilon.
\end{equation}
We thus now
 assume without loss of generality that the  step function  $g$ has a representation where this property holds,
and where the sets $A_i$ are fixed.
For the functions $f=\tilde\alpha,\tilde\alpha\cdot\beta\in C_b(D)$ we put 
$$\Phi_f(i,j):= \int_D\int_{A_i} \tilde K_h(x,A_j)\> d\tilde \omega_X(x)\cdot f(h)\> d\tilde\omega_D(h).$$
A short computation and the definition of $g$ then yield that
\begin{equation}\label{scalar-as-sum} 
\langle \tilde T_{f} g,g\rangle_{\tilde X} \> =\>\sum_{i,j=1}^n \overline{c_i} \> c_j \> \Phi_f(i,j).
\end{equation}
As $\tilde T_{\tilde\alpha}$ is positive definite, we obtain that  $\langle \tilde T_{\tilde\alpha} g,g\rangle_{\tilde X}\ge0$
for all choices of  $c_1,\ldots, c_n\in\mathbb C$. This means that the matrix $(\Phi_{\tilde\alpha}(i,j))_{i,j=1,\ldots, n}$
is positive semidefinite. 

On the other hand, as $(X,D,K)$ has property (T2), we know that
 $\beta\circ\pi:X\times X\to\mathbb C$ is positive definite, i.e., the matrix
 $(\beta(\pi(x_i,x_j)))_{i,j}$ is positive semidefinite for all choices of points $x_i\in A_i$, $i=1,\ldots,n$. 
As pointwise products of
 positive semidefinite matrices are positive semidefinite, we obtain that the matrix 
 $(\Phi_f(i,j)\cdot\beta(\pi(x_i,x_j)) )_{i,j=1,\ldots, n}$ is positive semidefinite. This means that for all
 choices of  $c_1,\ldots, c_n\in\mathbb C$,
\begin{align}\label{pd-teil}
0&\le \sum_{i,j=1}^n \overline{c_i} \> c_j \> \Phi_{\tilde\alpha}(i,j)\beta(\pi(x_i,x_j)) \notag\\
&=  \sum_{i,j=1}^n \overline{c_i} \> c_j \>  \int_D\int_{A_i} \tilde K_h(x,A_j)\> d\tilde \omega_X(x)
\cdot \tilde\alpha(h)\>\beta(\pi(x_i,x_j))\> d\tilde\omega_D(h).
 \end{align}
Furthermore, for $i,j=1,\ldots,n$ and $h\in D$, the relation $\int_{A_i} \tilde K_h(x,A_j)\> d\tilde \omega_X(x)>0$ 
implies $\sup_{x\in A_i}\tilde K_h(x,A_j)>0$, and this implies that $h\in \pi(A_i,A_j)$. 
This, the estimate (\ref{distanz-ab}), and $\|\tilde\alpha\|_\infty=1$ imply that
\begin{align}
\Bigl| \sum_{i,j=1}^n &\overline{c_i} \> c_j \>  \int_D\int_{A_i} \tilde K_h(x,A_j)\> d\tilde \omega_X(x)
\cdot \tilde\alpha(h)\>\beta(\pi(x_i,x_j))\> d\tilde\omega_D(h)\quad -\notag\\
&  -\> \sum_{i,j=1}^n\overline{c_i} \> c_j \>  \int_D\int_{A_i} \tilde K_h(x,A_j)\> d\tilde \omega_X(x)
\cdot \tilde\alpha(h)\>\beta(h)\> d\tilde\omega_D(h)\Bigr|\le \epsilon \|g\|_\infty^2 \cdot \tilde\omega_X(supp\> g)^2.
\notag\end{align}
We conclude from (\ref{pd-teil}) and (\ref{scalar-as-sum}) that
$$\langle \tilde T_{\tilde\alpha\>\beta} g,g\rangle_{\tilde X} \> =\>
\sum_{i,j=1}^n \overline{c_i} \> c_j \> \Phi_{\tilde\alpha\>\beta} (i,j)\in\mathbb C$$
has a distance from $[0,\infty[\subset\mathbb C$ which is at most
 $C\epsilon$ for some constant $C\ge0$ depending on $g$ only.
As in our approximation $\epsilon>0$ may be arbitrarily small, we obtain 
$\langle \tilde T_{\tilde\alpha\>\beta} g,g\rangle_{\tilde X} \in [0,\infty[$. As this holds  for all step functions $g$, 
 Lemma \ref{approx-pd} shows that $ \tilde T_{\tilde\alpha\>\beta}$ is positive definite as claimed.
Again, the support condition follows from Proposition \ref{support-condition-pl}.
\end{proof}

We now present some applications of Theorem \ref{cont-main-g}.

\begin{corollary}\label{COR-MAIN-1}
Let  $(X,D, K)$ and $(X,D,\tilde K)$ be commutative CAS with the same  $X,D$ and $\pi$.
Assume that  $(X,D, K)$ has property (T2), and that the identity ${\bf 1}$ is contained
 in  $\tilde S$. Then each character $\beta\in S$ 
 is positive definite on  $(D,\tilde *)$, and there is
a unique  $\mu\in M^1(\tilde S)$ with 
$$\beta(h)=\int_{\tilde S} \tilde\beta(h)\> d\mu(\tilde\beta) \quad\quad\text{for all }\quad h\in D.$$
\end{corollary}

\begin{proof} Use Theorem \ref{cont-main-g} with $\tilde\alpha={\bf 1}$.
 \end{proof}

In Remark \ref{counterex-72} we present an example which shows that the technical 
 condition  ${\bf 1}\in\tilde S $ in Corollary \ref{COR-MAIN-1} is necessary.

Here is a further consequence of Theorem \ref{cont-main-g} which generalizes Theorem 5.10 of \cite{Voit16}.

\begin{corollary}\label{COR-MAIN-5}
Let  $(X,D, K)$ and $(X,D,\tilde K)$ be  commutative CAS with the same 
 $X,D$ and $\pi$. Assume that  $(X,D, K)$ has property (T2), 
 and  that each character $\tilde\beta\in\tilde S$ of $(D,\tilde *)$ has the form
\begin{equation}\label{int-rep-annahme}
\tilde\beta(h)=\int_S \beta(h)\> d\mu(\beta) \quad\quad(h\in D)
\end{equation}
for some  $\mu\in M^1(S)$. 

 Then, for all   $\tilde\alpha,\tilde\beta\in\tilde S$, the product  $\tilde\alpha\cdot\tilde\beta$ 
is  positive definite  on $(D,\tilde *)$, and 
there is a unique probability 
measure $\delta_{\tilde\alpha}\hat*\delta_{\tilde\beta}\in M^1(\tilde S)$ with 
$(\delta_{\tilde\alpha}\hat*\delta_{\tilde\beta})^\vee=\tilde\alpha\cdot\tilde\beta$.
\end{corollary} 

\begin{proof} 
Theorem \ref{cont-main-g} shows that for  $\tilde\alpha\in \tilde S$ and $\beta\in S$ 
the product $\tilde\alpha\beta$ is positive definite on  $(D,\tilde *)$  with
$\tilde\alpha\beta=\int_{\tilde S} \gamma \> d\mu_{\tilde\alpha,\beta}(\gamma)$. Eq.~(\ref{int-rep-annahme})
 now implies that  for all   $\tilde\alpha,\tilde\beta\in\tilde S$, the
  product $\tilde\alpha\tilde\beta$ is positive definite on  $(D,\tilde *)$,
 and that the claimed integral representation holds. 
 \end{proof}

Clearly, Corollary \ref{cont-main-g}  is also a generalization of
 Theorem \ref{main-cont}. However, in practice it does not go far beyond   \ref{main-cont} 
by the following theorem which is closely related with \ref{uniqueness-complete-t2}:

\begin{theorem}\label{COR-MAIN-6}
Let  $(X,D, K)$ and $(X,D,\tilde K)$ be commutative CAS with the same $X,D$ and  $\pi$.
Assume that  $(X,D, K)$ has property (T2) and that  ${\bf 1}\in supp\> \tilde S$ holds.
Assume in addition that each character $\tilde\beta\in\tilde S$ 
 has the form
\begin{equation}\label{int-rep-annahme2}
\tilde\beta(h)=\int_S \beta(h)\> d\mu(\beta) \quad\quad(h\in D)
\end{equation}
for some  $\mu\in M^1(S)$. 
Then, $(D,*)=(D,\tilde *)$.
\end{theorem}

\begin{proof} 
We first recapitulate some facts on commutative hypergroups
 from \cite{J} which are well-known for lca groups and 
Gelfand pairs. For this let $(D,*)$ be any commutative hypergroup.
 Then the dual $\hat D$ can be identified with the symmetric spectrum 
\begin{equation}\label{spectrum}
\Delta_s(L^1(D,\omega_D)):=\{\phi\in L^1(D,\omega_D)^*:\> \phi\>\>\text{ multiplicative,}\>\> \>
\phi(f^*)=\overline{\phi(f)}\>\>\text{ for}\>\>\> f\in L^1(D)\}
\end{equation}
of the commutative Banach-$*$-algebra $(L^1(D,\omega_D), *,.^*)$ via
$$\alpha\mapsto\phi_\alpha  \quad\text{ with}\quad \phi_\alpha (f):=
\int_D \phi(x)\alpha(x)\> d\omega_D(x).$$
In particular, if  $\Delta_s(L^1(D,\omega_D))$ carries the Gelfand topology and  $\hat D$ the topology of 
compact-uniform convergence, then this mapping is a homeomorphism.
 We also recapitulate the well-known fact that
$\Delta_s(L^1(D,\omega_D))$ is the set of all extremal points in the set $P(D,*)$ of all positive linear
functionals on $L^1(D,\omega_D))$ with dual norm equal to 1; see e.g. Rudin \cite{Ru}.

We now apply these facts to our theorem. We first conclude from (\ref{int-rep-annahme2}) 
that each $\tilde\beta\in\tilde S\subset C_b(D)$ 
leads to a positive linear functional $\phi_{\tilde\beta} \in P(D,*)$.
Moreover, by Corollary \ref{COR-MAIN-1}, each $\beta\in S \subset C_b(D)$  has the form
\begin{equation}
\beta(h)=\int_{\tilde S} \tilde\beta(h)\> d\mu(\tilde\beta) \quad\quad(h\in D)
\end{equation}
for some  $\mu_\beta\in M^1(S)$. As $\phi_\beta$ is an extremal point, we conclude that  
 $\mu_\beta$  is a point measure.
In fact, if $\mu_\beta$ fails to be a point measure, then we may write  $\mu_\beta$ as 
 $\mu_\beta=\lambda\mu_1+(1-\lambda)\mu_2$ with
 different measures $\mu_1,\mu_2\in M^1(S)$ and $\lambda\in]0,1[$, i.e.,
$\phi_\beta$ would be a nontrivial convex combination of elements of $P(D,*)$ contradicting extremality.

In summary  $\mu_\beta$  is a point measure for each  $\beta\in S $ which proves $S\subset\tilde S$.
The same arguments also yield  $\tilde S\subset S$, i.e., we have $S=\tilde S$.
 Therefore, for all $x,y\in D$ and $\alpha\in S=\tilde S$,
$$(\delta_x*\delta_y)^\wedge(\alpha)=\int_D \bar\alpha\> d(\delta_x*\delta_y)=
\overline{\alpha(x)}\overline{\alpha(y)}=\ldots=(\delta_x\tilde*\delta_y)^\wedge(\alpha).$$
As the restricted Fourier-Stieltjes transform $ M^1(D)\to C_b(S)= C_b(\tilde S)$, $\mu\mapsto\hat \mu$ is
also injective (see \cite{J}), we obtain that
$\delta_x*\delta_y=\delta_x\tilde*\delta_y$ for all $x,y\in D$ as claimed.
 \end{proof}

\begin{remark}
We briefly discuss some implication of the preceding results.
For this we  recapitulate the original motivation in the introduction of 
  \cite{Voit16} for the study of generalizations of classical commutative association schemes.
Consider a  sequence $(G_n,H_n)_{n\in\mathbb N}$
 of
Gelfand pairs such that the 
double coset spaces  $G_n//H_n$  are homeomorphic with
some fixed locally compact space $D$. Modulo these  homeomorphims, we obtain 
 associated  double coset hypergroups $(D,*_n)$.
The spherical functions of $(G_n,H_n)$ then may be regarded as nontrivial continuous
multiplicative functions on  $(D,*_n)$.
For
many  examples of  series $(G_n,H_n)_n$, these functions  are
parameterized by some  spectral parameter set $\chi(D)$ independent on $n$, and
 the associated functions $\phi_n:\chi(D)\times D\to\mathbb C$
 can be embedded into a 
family of special functions which depend
analytically on  $n$  
in some parameter domain $A\subset\mathbb C$.
 In many cases,  these
 special functions are well known, and the product formulas for spherical functions can
 be written down  explicitly on $D$ with $n\in\mathbb N$  as  parameter.
 Based on
Carleson's theorem, a principle of analytic continuation 
(see e.g.~\cite{Ti}, p.186), one can often  easily
extend these positive product formulas 
 to a continuous range 
of parameters, say $n\in [1,\infty[$ such that for all these $n$ associated
commutative hypergroup structures  $(D,*_n)$ exist.

Besides  positive product formulas for  $\phi_n(\lambda,.)$  on $D$, there also exist dual product formulas
for the functions  $\phi_n(.,x)$ ($x\in D$) 
 on  suitable subsets of $\chi(D)$ for the group cases, i.e., for $n\in \mathbb N$.
In particular, 
 positive dual convolutions on the supports $S_n\subset \chi(D)$
of the Plancherel measures of the double coset hypergroups $G_n//H_n$ exist; see e.g.~ Theorem \ref{main-cont}.
 For many  examples, these dual convolutions
are  known and can be extended again by Carleson's theorem to
 positive dual convolutions for all $n\in [1,\infty[$.
However, for symmetric spaces of rank $\ge2$, this dual convolution  is usually a difficult business,
 and not very much is known in this respect.
When writing (the introduction of) \cite{Voit16}, the author hoped that a theory of continuous association schemes 
might lead to examples of commutative CAS associated with $(D,*_n)$ for all $n\in [1,\infty[$ such that 
 Theorem \ref{main-cont} or Corollary \ref{COR-MAIN-5}
leads at least to the existence of dual positive convolutions on $S_n$ in these cases.

This idea was motivated by 
 natural families $(K_h)_{h\in D}$ of Markov kernels on concrete spaces $X$
 which are are associated with commutative hypergroup structures on $D$ in \cite{Ki} and \cite{Bin}. 
In fact,  Kingman \cite{Ki} studies the Euclidean case $X=\mathbb R^2$ 
with $D=[0,\infty[$ where the associated hypergroups are the 
Bessel-Kingman hypergroups indexed by a continuous parameter.
Moreover, Bingham \cite{Bin} studies the spherical 
case $X=S^2:=\{x\in \mathbb R^3: \|x\|_2=2\}$ with
 $D=[-1,1]$ where the associated hypergroups on $[-1,1]$ are  related with ultraspherical polynomials.
Unfortunately, the  kernels $(K_h)_{h\in D}$  in \cite{Ki} and \cite{Bin} do not lead to
commutative CAS such that the theory of our paper cannot be applied there. This becomes clear from Theorem 
\ref{COR-MAIN-6} without 
discussing any details of these kernels from Theorem 
\ref{COR-MAIN-6}, as almost all conditions of  Theorem \ref{COR-MAIN-6}
 are satisfied for these examples.
 In fact, if the structures in  \cite{Ki} and \cite{Bin} would lead to  commutative CAS, then Theorem 
 \ref{COR-MAIN-6} could be applied to this structure as the CAS $(X,D,\tilde K)$ where we would have to take
 the CAS $(X,D, K)$ as a group case with a suitabele smaller group parameter than for
 $(X,D,\tilde K)$. We here notice that then
in particular (\ref{int-rep-annahme2}) holds by well-known explicit
 positive integral representations of the associated 
Bessel functions and ultraspherical polynomials; see e.g.~the survey of Askey \cite{A}.
\end{remark}

\section{Multiplicative functions and deformations}\label{ex-orbit-ass}

Following the well known notion of multiplicative functions, semicharacters, and
 characters  on commutative hypergroups
 (see \cite{BH}, \cite{J}, and Section 2 above), we here introduce a corresponding concept for 
commutative continuous association schemes. We in particular use it to construct deformed
 continuous association schemes $(X,D,\tilde K)$ from a given scheme $(X,D, K)$ with the
 same spaces $X,D$ and modified kernels $K$. This construction leads to examples of CAS which are 
beyond double coset examples and classical discrete association
 schemes where usually (T1) and (T2) do not hold.

\begin{definition}
A pair $(\alpha,\phi)\in C(D)\times C(X) $ of continuous functions 
 is called multiplicative on a commutative continuous association scheme $(X,D,\tilde K)$ if
 $\alpha\not\equiv 0$ and 
\begin{equation}\label{def-asso-semi}
T_h\phi(x)=\int_X \phi(y) K_h(x,dy)=\phi(x) \cdot \alpha(h) \quad\quad\text{for all}\quad h\in D, \> x\in X.
\end{equation}
A multiplicative pair $(\alpha,\phi)$ is called a semicharacter of  $(X,D,\tilde K)$, if in addition
\begin{equation}\label{def-semi-invo}
\alpha(\bar h)=\overline{\alpha(h)} \quad\quad\text{for all}\quad h\in D.
\end{equation}
A  semicharacter $(\alpha,\phi)$ is called a character, if $\alpha$ and  $\phi$ are bounded, and
 positive, if  $\alpha$ and  $\phi$ are $]0,\infty[$-valued.
\end{definition}

\begin{remarks}
\begin{enumerate}
\item[\rm{(1)}] Eq.~(\ref{def-asso-semi}) means that $\phi$ is a joint eigenfunction of all mean value operators
$T_h$, $h\in D$.
\item[\rm{(2)}] If $\alpha\equiv 1$, then (\ref{def-asso-semi}) is a mean 
value condition, i.e., $\phi$ may be seen as a harmonic function. Notice that for a compact CAS, all
harmonic functions are constant by Lemma \ref{unique-invar-measure}.
 We shall see in Remark \ref{example-unbounded-harmonic} 
that usually there might exist unbounded positive harmonic functions.
 It might be interesting to explore under which conditions on a CAS, all bounded harmonic functions are constant.
\item[\rm{(3)}]  For a multiplicative pair $(\alpha,\phi)$,  $\alpha$ is determined uniquely by $\phi$. 
The converse statement  is not correct as for any $\alpha\in C(D)$, the joint eigenspace
\begin{equation}\label{eigenspace}
E_\alpha:=\{\phi\in C(X):\> (\alpha,\phi) \quad\text{ multiplicative}\}
\end{equation}
is a vector space.
\item[\rm{(4)}]  Let $(X=G/H, D=G//H, K)$ be a commutative continuous association scheme which
 comes from some Gelfand pair $(G,H)$ with a connected Lie group $G$. Then,  any $\phi\in C(X)$ 
is contained in some joint eigenspace $E_\alpha$ for   $\alpha\in C(D)$ if and only if 
$\phi$ is a joint eigenfunction of all $G$-invariant differential operators on 
$X=G/H$; see e.g.~Prp.~IV.2.4 of Helgason \cite{Hel}.
\end{enumerate}
\end{remarks}

We next study relations between  $\alpha$ and $\phi$ for
  multiplicative pairs.

\begin{proposition}\label{mult-char}
If  $(\alpha,\phi)$  is multiplicative on $(X,D, K)$ with  $\phi\not\equiv 0$, 
then  $\alpha$  is a multiplicative function
of  $(D,*)$, i.e., for all $h_1,h_2,h\in D$,
$\alpha(h_1*h_2)=\alpha(h_1)\alpha(h_2)$.
Moreover, if $(\alpha,\phi)$ is in addition a semicharacter or character, the so is  $\alpha$ on   $(D,*)$.

Conversely,  for each  nontrivial multiplicative function  $\alpha\in C(D)$  on   $(D,*)$
there exist functions $\phi\in C(X)$  with  $\phi\not\equiv 0$ such
 that $(\alpha,\phi)$  is multiplicative on $(X,D, K)$.
\end{proposition}

\begin{proof} Let  $(\alpha,\phi)$ be multiplicative as described.
For $x\in X$,  $h_1,h_2\in D$ we have
\begin{align}
\phi(x)\cdot \alpha(h_1*h_2)&=\phi(x)\cdot \int_D  \alpha(h)\> d(\delta_{h_1}*\delta_{h_1})(h)=
 \int_D T_h\phi(x)\> d(\delta_{h_1}*\delta_{h_1})(h)\notag\\
&= T_{h_1}\circ T_{h_1}\phi(x) =\phi(x)\cdot \alpha(h_1)\alpha(h_2)
\notag
\end{align}
Taking $x\in X$ with $\phi(x)\ne0$ leads to the first claim. The second statement is clear.

For the last statement we conclude from   Lemma \ref{Tmult} that for each $g\in C_c(X)$ and each
nontrivial multiplicative function  $\alpha\in C(D)$  on   $(D,*)$,
 the function $\phi:=T_{\alpha^-}g$ with $\alpha^-(h):=\alpha(\bar h)$
satisfies $T_h\phi=\alpha(h)\cdot\phi$. We still have to check that we can choose $g$ such that 
$\phi\not\equiv 0$ holds. As $\alpha(e)=1$, we find a neighborhood $U\subset D$ of
 $e$ on which $\Re \alpha\ge0$ holds.
Now fix some $x\in X$ and a neighborhood $W\subset X$ of $x$ with $\pi(x,W)\subset U$. Now choose $g\in C_c(G)$ with
$g\ge 0$,  $g\not\equiv 0$, and $supp\> g\subset W$. Then 
$$\Re\phi(x)=\int_U \int_W g(y)\> K_h(x,dy)\> \Re\alpha(h)\> d\omega_D(h)>0$$
as claimed. 
\end{proof}

 (T1) and (T2) lead to a further standard construction of multiplicative pairs:

\begin{lemma}\label{char-on-X}
Let $(X,D,K)$ be a commutative CAS with  (T1) or  (T2).
 Then for all multiplicative function $\alpha\in C(D)$ on $(D,*)$ and $z\in X$, the pair
 $(\alpha, \tilde\alpha  :=\alpha\circ\pi_z)$  is  multiplicative on $(X,D, K)$.
\end{lemma}

\begin{proof} Assume first that (T1) holds.
Then  for  $x\in X$, $h\in D$,
\begin{equation}\label{T1-multpair}
T_h\tilde\alpha(x)= T_h(\alpha\circ \pi_z)(x)=\alpha(\pi_z(x)*h)= \alpha(\pi_z(x))\alpha(h)=
\tilde\alpha(x)\alpha(h)
\end{equation}
as claimed.

Assume now that (T2) holds. We conclude from the last assertion of  Proposition \ref{mult-char}
that for all $g\in C_c(X)$, the function $\phi_g:=T^{\alpha^-\circ\pi}g=T_{\alpha^-}g$ satisfies
$T_h\phi_g = \alpha(h)\cdot\phi_g$ for $h\in D$. On the other hand, as $\alpha\circ\pi$ 
is uniformly continuous on compact subsets of $X\times X$, we find relatively
 compact neighborhoods $U_n\subset X$ of $z$ with $U_{n+1}\subset U_n$ for all $n$ and 
$\bigcap_n U_n=\{z\}$
such that for all $g_n\in C_c(X)$ with $supp\> g_n\subset U_n$, 
$g_n\ge0$, and $\int_X g_n\> d\omega_X=1$, we have
$$\phi_{g_n}(x)=\int_X \alpha^-(\pi(x,y))\> g_n(y)\> d\omega_X(y)\longrightarrow 
 \alpha^-(\pi(x,z))=\alpha(\pi(z,x))=\alpha(\pi_z(x))$$
uniformly on compacta w.r.t.~$x$. Hence, the limit $\tilde\alpha  :=\alpha\circ\pi_z$ also satisfies
$T_h\tilde\alpha = \alpha(h)\cdot\tilde\alpha$ for $h\in D$ as claimed.
\end{proof}

The arguments of the preceding proof and in particular in (\ref{T1-multpair})
 can be combined in a different way and lead finally to:

\begin{theorem}\label{commT2T1} For each commutative CAS, (T2) implies (T1).
\end{theorem}

\begin{proof}
Let $(X,D,K)$ be a commutative CAS with (T2) and with associated commutative hypergroup $(D,*)$.
Let $\alpha\in\hat D$ be a character, and let  $z\in X$ and $\tilde\alpha:=\alpha\circ\pi_z$. Then, by property
(T2) and Lemma \ref{char-on-X}, $T_h\tilde\alpha = \alpha(h)\cdot\tilde\alpha$ for $h\in D$. Hence, for $x\in X$
and $h\in D$,
$$T_h(\alpha\circ \pi_z)(x)=T_h\tilde\alpha(x)=\tilde\alpha(x)\alpha(h)= \alpha(\pi_z(x))\alpha(h)=
\alpha(\pi_z(x)*h).$$
As this equation is linear in $\alpha$, we obtain that
$$T_h(f\circ \pi_z)=f_h\circ \pi_z$$
for all $h\in D$, $z\in X$ and $f\in C_b(D)$ of the form $f=\check\mu$ with $\mu\in M_b(\hat D)$
 in the notation of Section \ref{deffacts}(4). On the other hand, it is well-known from hypergroup theory
that, again with the notion of Section \ref{deffacts}(4), the set
$\{\check g: \> g\in C_c(\hat D)\}$ is a $\|.\|_\infty$-dense subspace of $C_0(D)$; see e.g. Theorem 2.2.32(vii)
of \cite{BH}.
We thus conclude that $T_h(f\circ \pi_z)=f_h\circ \pi_z$ for all $f\in C_c(D)$, $h\in D$ and  $z\in X$
as claimed.
\end{proof}

Theorem \ref{commT2T1}  and Proposition \ref{compact-t2}   imply:

\begin{corollary} Each compact commutative CAS is strong.
\end{corollary}

\begin{remark}
\begin{enumerate}
\item[\rm{(1)}] Notice that in the proof of the first part of 
Lemma \ref{char-on-X}, (T1) is needed for the specific $z\in X$
with $\tilde\alpha  :=\alpha\circ\pi_z$ only. 
\item[\rm{(2)}]   If $(X=G/H, D=G//H, K)$ is a commutative CAS
 associated with the Gelfand pair $(G,H)$, then (T1) and (T2) hold by 
\ref{group-to-strong-cont-asso-scheme}. Hence, for a
 multiplicative function  $\alpha$ of $(D,*)$ and $z\in G$, 
we may take $\tilde\alpha\in C(X)$ with
 $\tilde\alpha(xH):=\alpha(\pi_{zH}(xH))= \alpha(Hz^{-1}xH)$ for
$x\in G$.
\item[\rm{(3)}] Let $(X=G/H, D=G//H, K)$ be a commutative CAS which
 comes from some  symmetric space $G/H$.
 Then for each  multiplicative 
$\alpha\in C(D)$,  $E_\alpha\subset C(X)$ is the joint eigenspace of the algebra $D(G/H)$ of 
 all $G$-invariant differential operators on $X=G/H$ (with suitable eigenvalues).
In this case  $E_\alpha$ is completely known by Kashiwara et.~al.~\cite{Kash}; see also the description 
of the results in Section II.4.1 of  \cite{Hel}. As this goes beyond the scope of this paper, we skip details.
 For some interesting concrete examples of functions in $E_\alpha$ on hyperbolic planes 
we refer to the introduction of  \cite{Hel}.
\item[\rm{(4)}] We expect that the well-established representation theory of compact hypergroups and the arguments of the proof of 
Theorem \ref{commT2T1}  yield that for all compact CAS, (T2) implies (T1). Proposition \ref{compact-t2} then would imply that each
compact CAS is strong.
\end{enumerate}
\end{remark}

We now restrict our attention to positive semicharacters   $(\alpha_0,\phi_0)$  of some commutative CAS
$(X,D,K)$. It is  well known from  \cite{V1} or Section 2.3 of \cite{BH} that then the positive semicharacter
$\alpha_0$ of  $(D,*)$ leads to a deformed commutative hypergroup $(D,\tilde *)$
with the deformed convolution of point measures
\begin{equation}\label{def-convo-ass}
\delta_{h_1} \tilde* \delta_{h_2}:= \frac{\alpha_0}{\alpha_0(h_1)\cdot \alpha_0(h_2)}
\cdot( \delta_{h_1}*\delta_{h_2}) \in M^1(D) \quad\quad (h_1,h_2\in D)
\end{equation}
where the identity and involution of $(D,*)$ are not changed. Moreover,  if $\omega_D$ is a  Haar measure  of  $(D,*)$,
then  
\begin{equation}\label{def-haar-ass}
\tilde \omega_D:= \alpha_0^2 \cdot\omega_D\in  M^+(D)
\end{equation}
is a Haar measure of  $(D,\tilde *)$ by \cite{V1}. We now show that  $(\alpha_0,\phi_0)$ also
leads to a deformed commutative CAS $(X,D,\tilde K)$ which is associated with
  $(D,\tilde *)$.

\begin{proposition}\label{deformed-ass-scheme}
Let $(\alpha_0,\phi_0)$ be a positive  semicharacter as above. Define the deformed
 kernel $\tilde K$ from $X\times D$ to $X$ with
$$\tilde K_{h}(x, A):= \frac{1}{\alpha_0(h)\phi_0(x)}\int_A\phi_0(y)\>  K_{h}(x, dy)
 \quad\quad (h\in D,\> x\in X, \> A\in\cal B(X)).$$
Then $(X,D,\tilde K)$ is a  commutative CAS which is associated with
 $(D,\tilde *)$ above. Moreover:
\begin{enumerate}
\item[\rm{(1)}] If $\omega_X\in  M^+(X)$ is an invariant measure of  $(X,D,K)$, then
$\tilde\omega_X:= \phi_0^2\cdot\omega_X \in  M^+(X)$ is an invariant measure of  $(X,D,\tilde K)$.
\item[\rm{(2)}] Assume that $(X,D,K)$ has property (T1), and let $\tilde\alpha_0\in C(D)$ 
be a positive semicharacter and $z\in X$. Consider the positive semicharacter
  $(\alpha_0,\tilde\alpha_0:=\alpha_0\circ\pi_z)$ of $(X,D,K)$ according to Lemma \ref{char-on-X}.
Then for all $f\in C(D)$, $h\in D$ and $y\in X$,
\begin{equation}\label{prop-T1-schwach}
f(h\tilde* \pi_z(y))= (\tilde T_h(f\circ\pi_z))(y)
\end{equation}
where $\tilde T_h$ is the operator associated with the kernel $\tilde K_h$. This means that (T1) holds 
for $(X,D,\tilde K)$ for the specific $z\in X$.
\end{enumerate}
\end{proposition}

\begin{proof}
Notice that $\tilde K$ satisfies
$\tilde K_h(x,X)=1$ for $h\in D$, $x\in X$.
This normalization and the continuity of $K$ show that $\tilde K$ is a 
continuous Markov kernel from $X\times D$ to $X$. Moreover,  $\tilde K$ clearly  satisfies the conditions of
 \ref{central-def}(1) and (2) with the projection $\pi$ of the scheme $(X,D,K)$.
 We next check axiom \ref{central-def}(3). For $h_1,h_2\in D$, $x\in X$ and $A\in \cal B(X)$ we have
 \begin{align}
\tilde  K_{h_1}\circ\tilde  K_{h_2} (x, A) &=
\int_X \tilde   K_{h_2} (y, A) \> \tilde  K_{h_1}(x,dy)\notag\\
&=\frac{1}{\alpha_0(h_1) \phi_0(x)} \int_X \tilde    K_{h_2} (y, A) \> \phi_0(y)\>  K_{h_1}(x,dy)\notag\\
&=\frac{1}{\alpha_0(h_1)\alpha_0(h_2) \phi_0(x)} \int_X \int_A  \phi_0(z) \>
 K_{h_2} (y,dz) \>  K_{h_1}(x,dy)\notag\\
&=\frac{1}{\alpha_0(h_1)\alpha_0(h_2) \phi_0(x)}  \int_D \int_A  \phi_0(z) \> K_{h} (x,dz)
 \> d(\delta_{h_1} *\delta_{h_2})(h)\notag\\
&=\frac{1}{\alpha_0(h_1)\alpha_0(h_2) } \int_D \tilde K_{h} (x,A) \>   \alpha_0(h)
 \> d(\delta_{h_1} *\delta_{h_2})(h)\notag\\
&=\int_D \tilde K_{h} (x,A) \>d(\delta_{h_1} \tilde *\delta_{h_2})(h).
\notag
\end{align}
For the adjoint relation in
 \ref{central-def}(4), we observe for $f_1,f_2\in C_c(X)$ and $h\in D$ that
 \begin{align}\label{connection-deformed-adjoint}
\int_X f_1 \cdot \tilde T_{h}f_2 \> d\tilde\omega_X &= 
\int_X\int_X f_1(x)  f_2(y) \tilde K_{h}(x,dy) \> d\tilde\omega_X(x)
\\
&=\int_X\int_X  \frac{\phi_0(y)\> f_1(x) \>   f_2(y)}{\alpha_0(h) \>\phi_0(x)} K_{h}(x,dy) \>
\phi_0(x)^2\> d\omega_X(x)
\notag\\
&= \frac{1}{\alpha_0(h)} \int_X \phi_0f_1 \cdot  T_{h}(\phi_0f_2) \> d\omega_X.
\notag
\end{align}
This, $\alpha_0(\bar h)=\alpha_0(h)$, and the  adjoint relation \ref{central-def}(4)
for $(X,D,K)$  now lead to  \ref{central-def}(4) in the deformed case.
This completes the proof of the main statement and of part (1).

Finally, for the proof of (2) we observe that (T1) for $(X,D,K)$ and (\ref{def-convo-ass}) imply that
\begin{align}
f(h\tilde*\pi_z(y))  &= \frac{1}{\alpha_0(h)\alpha_0(\pi_z(y) )}
\cdot (\delta_h*\delta_{\pi_z(y) })(\alpha_0f)\notag\\
&=\frac{1}{\alpha_0(h)\tilde\alpha_0(y)} (\alpha_0f)_h(\pi_z(y))
\notag\\
&=\frac{1}{\alpha_0(h)\tilde\alpha_0(y)} T_h(\alpha_0f)(\pi_z(y))
\notag\\
&= \int_X f(\pi_x(w)) \> \tilde K_h(y,dw) \quad =\quad (\tilde T_h(f\circ\pi_z))(y).
\notag
\end{align}
\end{proof}

\begin{remark}
Let  $(X,D,K)$ be a commutative strong CAS, $\alpha_0\in C(D)$ 
 a positive semicharacter of $(D,*)$, and  $z\in X$. Then 
  $(\alpha_0, \alpha_0\circ\pi_z)$ is a positive semicharacter of $(X,D,K)$ by 
\ref{char-on-X}. Consider the associated deformed CAS $(X,D,\tilde K)$ 
which has property (T1) for  $z\in X$ by  \ref{deformed-ass-scheme}(2).
A short computation similar to the proof  of  \ref{deformed-ass-scheme}(2) shows that
we also have 
\begin{equation}
\tilde T_fg(x)= \tilde T^{f\circ \pi}g(x) \quad\quad\text{for all}\quad f\in C_c(D),\> g\in C_c(X),
\end{equation}
i.e.,  (T2) also holds for $(X,D,\tilde K)$ and  $z\in X$.
\end{remark}

In the setting of Prp.~\ref{deformed-ass-scheme}, the semicharacters of $(X,D,K)$ and 
 $(X,D,\tilde K)$ are closely related. This is well-known for hypergroup deformations
 from  \cite{V1} or Section 2.3 of \cite{BH}.

\begin{lemma}\label{mod-semichar} Let  $(X,D,K)$ and $(X,D,\tilde K)$ 
be related as in  Proposition \ref{deformed-ass-scheme}. Then
$$\{ (\alpha/\alpha_0, \phi/\phi_0): \quad (\alpha, \phi) \quad 
\text{a semicharacter of $(X,D,K)$}\}$$
is the set of all  semicharacters of $(X,D,\tilde K)$.
\end{lemma}

\begin{proof}
Let  $(\alpha, \phi)$ be a semicharacter of  $(X,D,K)$.
Then, for $h\in D$, $x\in X$,
\begin{equation}
\tilde T_h(\phi/\phi_0)(x)=\frac{1}{\alpha_0(h)\phi_0(x)}\int_X \phi(y) K_h(x,dy)
=\frac{\alpha(h)\phi(x)}{\alpha_0(h)\phi_0(x)}.
\end{equation}
This shows  (\ref{def-asso-semi}) for $(\alpha/\alpha_0, \phi/\phi_0)$ and thus the first part of the lemma.

For the converse statement we notice that $(1,1)$ is a  positive character of  $(X,D,K)$. Hence,
$(1/\alpha_0,1/\phi_0)$ is a positive semicharacter  of $(X,D,\tilde K)$ by the first part
of the lemma. We now  apply the first part
of the lemma to $(1/\alpha_0,1/\tilde\alpha_0)$ where the rules of  $(X,D,K)$ and $(X,D,\tilde K)$ are  interchanged.
 This readily proves that each semicharacter of  $(X,D,\tilde K)$ has the form as stated in the lemma.
\end{proof}

\section{Orbit schemes and their deformations}

We now study  examples of  semicharacters $(\alpha,\phi)$ and associated deformations
beyond the case
 $\phi=\alpha\circ\pi_z$ for $z\in X$ under condition (T1) or (T2).
We know from \ref{group-to-cont-asso-scheme} and \ref{group-to-strong-cont-asso-scheme} that
Gelfand pairs $(G,H)$ lead to commutative, strong CAS $(G/H, G//H,K)$. parametrized
Typical examples  are given by the following  orbit construction; see  \cite{J} for the  background.

\begin{orbit-schemes}\label{orbit-schemes-sec}
Let $G$ be a locally compact abelian group and $H\subset Aut(G)$ a compact group of
  automorphisms which acts continuously.
 Form the semidirect product $G\rtimes H$ which contains $H$ as a compact subgroup canonically. 
Then  $(G\rtimes H, H)$ is a Gelfand pair, and we may identify $(G\rtimes H)/H$ with $G$ via $(g,h)H\sim g$
($g\in G, h\in H$), and  $(G\rtimes H)//H$ with the space $G^H$ of all $H$-orbits in $G$ via 
 $H(g,h)H\sim g^H:=\{h(g): \> h\in H\}$ where all spaces carry the quotient topology.
 Consider the associated commutative
strong CAS
$$\Lambda:= ( X:= (G\rtimes H)/H=G, \> D:=(G\rtimes H)//H=G^H, \> K)$$
where the double coset hypergroup $(D,*)$ has the identity $\{e\}$ ($e$ the identity of $G$), the involution 
$\overline{ g^H}= (g^{-1})^H$, and the convolution
$$\delta_{g_1^H}*\delta_{g_2^H}:= \int_H\delta_{(g_1\cdot h(g_2))^H}\> d\omega_H(h) \quad\quad (g_1,g_2\in G)$$
for the normalized Haar measure $\omega_H$ of $H$. The Markov kernel $K$ is given by 
$$K_{g_1^H}(g_2, A):=\int_H\delta_{g_2\cdot h(g_1)}(A)\> d\omega_H(h) =\omega_H(\{h\in H:\> g_2\cdot h(g_1)\in A\})$$
for $g_1,g_2\in G$, $A\in \cal B(G)$. By the proof of
 Proposition \ref{group-to-cont-asso-scheme},
the map $\pi:X\times X\to D$ is given by $\pi(g_1H, g_2H):=(g_1^{-1}g_2)^H$. Moreover, 
 if $\omega_G$ is some Haar measure of $G$, then  $\omega_G\times\omega_H$ 
is a Haar measure of $G\rtimes H$, and we may choose the measures $\omega_X,\omega_D$ of  $\Lambda$
as $\omega_X:=\omega_G$ and 
$\omega_D:=\phi(\omega_G)$ for the orbit map $\phi:G\to G^H$ with $\phi(g)=g^H$.

We call $\Lambda$ the orbit scheme associated with $(G,H)$.
\end{orbit-schemes}

In this setting we have multiplicative pairs as follows:

\begin{lemma}\label{orbit-deformed-ass}
Let  $\phi\in C(G)$ be  multiplicative, i.e.,  
$\phi(g_1g_2)=\phi(g_1)\phi(g_2)$ for $g_1,g_2\in G$, with $\phi\not\equiv0$.
 Form   $\alpha\in C(D)$ with
\begin{equation}\label{orbit-char-def}
\alpha(g^H):=\int_H \phi(h(g))\> d\omega_H(h) \quad\quad (g\in G).
\end{equation}
Then  $(\alpha,\phi)$ is multiplicative on $\Lambda$. Moreover:
 \begin{enumerate}
\item[\rm{(1)}] If  $\phi\in C(G)$ is a character of $G$, then  so is the pair $(\alpha,\phi)$ on $\Lambda$. 
\item[\rm{(2)}] Let $\phi_0\in C(G)$ be positive and multiplicative such that the associated 
 $\alpha_0$ satisfies
$\alpha_0(\overline{g^H})=\alpha_0({g^H})$ for $g\in G$. Then $(\alpha_0,\phi_0)$ is a positive semicharacter on  $\Lambda$. 
Hence, with the associated 
kernel $\tilde K$
from Proposition \ref{deformed-ass-scheme}, $(X,D,\tilde K)$ is a commutative CAS.
\end{enumerate}\end{lemma}

\begin{proof} $\phi$ satisfies  $\phi(e)=1$ and $\phi(g)\ne0$ with
$\phi(g^{-1})=\phi(g)^{-1}$ for $g\in G$. In particular, we obtain $\alpha\not\equiv0$. Moreover, as 
for $g_1,g_2\in G$,
\begin{align}
T_{g_1^H} \phi(g_2)&=\int_G \phi(y)\> K_{g_1^H}(g_2,dy)= \int_H \phi(g_2 \cdot h(g_1))\> d\omega_H(h) 
\notag\\
&=
\int_H\phi(g_2 )  \cdot \phi( h(g_1))\> d\omega_H(h) = \phi(g_2 ) \cdot \alpha(g_1^H),
\notag
\end{align}
the first statement is clear. Parts (1) and (2) are  then  clear.
\end{proof}

\begin{remark}
Let  $\phi_0\in C(G)$ and $\alpha_0\in C(D)$ be as in the setting of Lemma \ref{orbit-deformed-ass}(2).
\begin{enumerate}
\item[\rm{(1)}] We are mainly interested in  nontrivial  $\phi_0$, i.e., $\phi_0\not\equiv 1$,
 as otherwise \ref{orbit-deformed-ass}(2)
does not lead to a nonidentical deformation.
This clearly works for noncompact groups $G$ only.
\item[\rm{(2)}] The push forwards 
$\pi_x(\tilde\omega_X)\in M^+(D)$ of  invariant measures $\tilde\omega_X$ as in \ref{deformed-ass-scheme}(2)
for $x\in X$ usually will not be  Haar measures on $(D,\tilde*)$; for examples see
below. Therefore, by Lemmas \ref{proj-haar-m-T1} and \ref{proj-haar-m-T2},
the deformed CAS of  Proposition \ref{orbit-deformed-ass}
 usually do not have  (T1) and (T2).
\item[\rm{(3)}] Consider the original orbit scheme $\Lambda$  as in Section \ref{orbit-schemes-sec}.
Let $\hat G$ be the dual group on which $H$ also acts via $h(\alpha)(g):=\alpha(h(g))$. Consider the orbit maps
$\Phi:G\to G^H=D, \> g\mapsto  g^H$ and $\hat\Phi:\hat G\to ( \hat G)^H, \> \alpha\mapsto  \alpha^H$.
It is well-known that $( \hat G)^H$ can be identified with the dual $(D,*)^\wedge$ via
 $\hat\Phi(\alpha)(g^H):=\int_H\alpha(h(g))\> d\omega_H(h)$, and that
for a Haar measure $\omega_G$ and its associated Plancherel measure   $\omega_{\hat G}$ 
(which is a Haar measure of $\hat G$), the push forwards $\Phi(\omega_G)\in M^+(D)$ and 
 $\hat\Phi(\omega_{\hat G})\in M^+((D,*)^\wedge)$ are a Haar measure of  $(D,*)$ and its Plancherel measure respectively;
 see \cite{J}.
 In particular, the support $S$ of the Plancherel measure is equal to  $(D,*)^\wedge$ for these examples.

 On the other hand, if  $\phi_0\in C(G)$ is a positive multiplicative function with
 $\alpha_0\not\equiv 1$ as in  \ref{orbit-deformed-ass}(2),
 then for the associated deformed hypergroup $(D,\tilde *)$
the support $\tilde S$ of the Plancherel measure  is a proper subset of 
 $(D,\tilde*)^\wedge$ for $\tilde\alpha_0\not\equiv 1$. In fact, we even have ${\bf 1}\not\in \tilde S$ by
 \cite{V1}.
 \end{enumerate}\end{remark}
 
It is an interesting  problem whether $(D,\tilde*)^\wedge$ or  $\tilde S$ carry  dual positive convolutions.
Generally, the answer is negative for $(D,\tilde*)^\wedge$; see  below.
On the other hand, for $\tilde S$ there exist some positive results.
In fact, for $\tilde S$, this problem is closely related with a property of $\alpha_0$:

\begin{lemma}\label{equiv-ex-dual-convo}
 In the setting of Lemma \ref{orbit-deformed-ass}(2), the following statements are equivalent:
\begin{enumerate}
\item[\rm{(1)}] $1/\alpha_0\in C_b(D)$ is  positive definite on the orbit hypergroup $(D,*)$;
\item[\rm{(2)}] for all
 $\tilde\alpha,\tilde\beta\in \tilde S$ there exists
$\tilde\mu_{\tilde\alpha,\tilde\beta}\in M^1(\tilde S)$ with
 $\tilde\alpha(x)\tilde\beta(x)
=\int_{\tilde S} \gamma(x)\> d\tilde\mu_{\tilde\alpha,\tilde\beta}(\gamma)$ for $x\in D$.
\item[\rm{(3)}] each character $\tilde\alpha\in\tilde S$ is positive definite on  $(D,*)$.
\end{enumerate}
\end{lemma}

\begin{proof} For $(1)\Longrightarrow(3)$ assume  that $1/\alpha_0$ is positive definite on $(D,*)$. For
$\tilde\alpha\in \tilde S$ we find a unique  $\alpha\in (D,*)^\wedge$ with
$\tilde\alpha=\alpha/\alpha_0$ by \cite{V1}. As $(D,*)^\wedge$  carries a
 dual positive convolution on  $(D,*)^\wedge$, we see that  $\tilde\alpha= \alpha\cdot (1/\alpha_0)$ is 
 positive definite on $(D,*)$ as claimed.

 $(3)\Longrightarrow(2)$ follows from (T2) for $(X,D,K)$, $S=(D,*)^\wedge$, and Corollary \ref{COR-MAIN-5}   

Finally, for  $(2)\Longrightarrow(1)$ we take $\tilde\alpha:=\tilde\beta:=1/\alpha_0\in \tilde S$ in (2). 
The homeomorphism $S\to\tilde S, \> \alpha\mapsto \alpha/\alpha_0$ 
 then yields that $1/\alpha_0$ is the inverse Fourier transform of some $\mu\in  M^1( S)$ as claimed.
\end{proof}

\begin{remark}
Consider some example in the setting of Lemma \ref{orbit-deformed-ass}(2) with $\alpha\not\equiv {\bf 1}$ such that one and thus all 
statements of Lemma \ref{equiv-ex-dual-convo} hold. Then $(X,D,K)\ne(X,D,\tilde K)$ and $(D,*)\ne (d,\tilde *)$.
 This shows that the technical condition $ {\bf 1}\in \tilde S$ in Theorem \ref{COR-MAIN-6} is essential.
\end{remark}

We now present some examples for the theory of Sections  \ref{orbit-schemes-sec}-\ref{equiv-ex-dual-convo}.

\begin{examples}\label{examples-orbit-ass}
 Fix an integer $d\ge1$ and put  $G:=(\mathbb R^d,+)$ and $H:=O(d)$ as orthogonal group acting on $G$.
 We use the canonical identification $D=[0,\infty[$. Then $(D,*)$ is the so-called
 Bessel-Kingman hypergroup of index $\alpha=d/2-1$; see e.g.~\cite{BH}, \cite{J}, and \cite{Ki}.

The multiplicative functions on $G$ have the form
 $\phi(x)=\phi_z(x):=e^{i\langle z,x\rangle}:=e^{i\sum_{k=1}^d z_kx_k}$ with 
 $z\in \mathbb C^d$. Then $\phi_z$ is a character precisely for $z\in  \mathbb R^d$,
 and $\phi_z$ is positive 
precisely for  $z\in  i\cdot\mathbb R^d$.

In the first case, the character $\alpha_z\in (D,*)^\wedge$ associated 
with $\phi_z$ ($z\in  \mathbb R^d$)
 according to
 (\ref{orbit-char-def}) is given by $\alpha_z(w)=j_\alpha(w\cdot \|z\|_2)$
with the the modified Bessel functions
\begin{equation}
j_\alpha(y):=_0F_1(\alpha+1;-y^2/4) \quad\quad(y\in\mathbb C).
\end{equation}
Please be careful with the different meaning of the parameter $\alpha$  and the functions $\alpha_z$.

In the second case the positive multiplicative function $\alpha_z\in C(D)$ 
associated with $\phi_z$ ($z\in  i\cdot \mathbb R^d$) is given by
 $\alpha_z(w)=j_\alpha( iw\cdot \|z\|_2)$.
In particular, as  the hypergroup $(D,*)$ is symmetric, all conditions of \ref{orbit-deformed-ass}
are satisfied in this case, i.e.,  $(\alpha_z,\phi_z)$ is a positive semicharacter of our orbit  CAS 
$(\mathbb R^d, [0,\infty[,K)$, and  $(\alpha_z,\phi_z)$ 
leads to a deformation  for 
each  $z\in  i\cdot\mathbb R^d$.  We now study examples for the equivalent conditions of 
\ref{equiv-ex-dual-convo}, and we discuss  whether the complete dual
$(D,\tilde*)^\wedge$ carries a dual positive convolution.

Before doing this, we notice that for each $c>0$, the map $x\mapsto cx$ is a hypergroup automorphism
on $(D=[0,\infty[,*)$. This ensures that we may restrict our attention to $z\in i\cdot\mathbb R^d$ with
$\|z\|_2=1$ without loss of generality.
\end{examples}

\begin{examples}\label{examples-orbit-ass2}
\begin{enumerate}
\item[\rm{(1)}] Let $d=1$, i.e., $\alpha=-1/2$ and $j_{-1/2}(x)=\cos x$.  Let 
$z=\pm i$. The deformed hypergroup $(D=[0,\infty[,\tilde*)$ 
 is then the so-called cosh-hypergroup; see \cite{Z} and Sections  3.4.7 and 3.5.72 of \cite{BH}. The
characters  are given by
$$\alpha_{\lambda}(x):=\frac{\cos(\lambda x)}{\cosh x } \quad(x\in [0,\infty[,\>\>
\lambda\in [0,\infty[\cup i\cdot[0,1])$$
where in this parameterization, $\alpha_{\lambda}$ is in the support  of the 
Plancherel measure precisely for $\lambda\in [0,\infty[$. Using
$$\frac{\cos(\lambda x)}{\cosh x}=
\frac{1}{2}\int_{-\infty}^\infty \frac{\cos(tx)}{\cosh((t+\lambda)\pi/2)}\>
  dt
\quad \text{for}\quad \lambda\in\mathbb C, \> |\Im \lambda|<1$$
(see (1) in \cite{Z} and references there), we see that the first condition of 
\ref{orbit-deformed-ass} holds. Hence, by \ref{orbit-deformed-ass}, the support of the 
Plancherel measure of $(D,\tilde *)$ carries a positive dual convolution.
 This convolution was computed explicitly in \cite{Z}.
We remark, that by \cite{Z}, there does not exist a  positive dual convolution on the complete dual space.
\item[\rm{(2)}] Let $d=3$, i.e.,  $\alpha=1/2$ and $j_{1/2}(x)=\frac{\sin x}{x}$. Let 
$\|z\|_2=1$. In this case, the deformed hypergroup $(D,\tilde*)^\wedge$ 
 is the Naimark hypergroup with convolution
$$\delta_x\tilde*\delta_y = \frac{1}{\sinh x \sinh y}\int_{|x-y|}^{x+y} \sinh t \> \delta_t\> dt 
\quad\quad (x,y\in [0,\infty[);$$
see \cite{J}, \cite{BH} and \cite{Z}. This example is also isomorphic
 with the double coset hypergroup $SU(1,1)//SU(2)$, and it follows from the work of
 Flensted-Jensen and Koornwinder \cite{FK1}, \cite{FK2} that all bounded hermitian spherical 
functions are positive definite on  $SU(1,1)$ and that thus  the complete dual 
$(D,\tilde*)^\wedge$ carries a dual positive convolution.
 The dual convolution was computed explicitly in \cite{Z}.
\end{enumerate}
\end{examples}

Here is a short list of further  examples for   \ref{orbit-schemes-sec}-\ref{equiv-ex-dual-convo}:

\begin{examples}\label{examples-orbit-ass3}
\begin{enumerate}
\item[\rm{(1)}] Put  $G:=(\mathbb Z,+)$ and $H:=\{\pm1\}$ which acts multiplicatively on $G$.
Then  $D=\mathbb N_0$ in a canonical way, and  $(D,*)$ is the so-called discrete polynomial hypergroup associated
with T-polynomials of the first kind; see e.g.~\\ref{La} and cite{BH}. 
The associated
  transition matrices are given by
$$S_0=I_{\mathbb Z}, \quad S_k(x,y)=\frac{1}{2}\delta_{k,|x-y|} 
\quad(k\in \mathbb N, \> x,y\in \mathbb Z)$$
with the Kronecker-$\delta$.

Similar to Examples \ref{examples-orbit-ass} and  \ref{examples-orbit-ass2}(1), we consider
 $\phi_z(k):=e^{zk}$ for $k\in \mathbb Z, z\in\mathbb R$. Then $\alpha_z(n)=\cosh(zn)$ for 
$n\in\mathbb N_0$, and we obtain deformed CAS similar to \ref{examples-orbit-ass2}(1).
 For further details on this discrete example see also Example 5.11 of \cite{Voit16}.

\item[\rm{(2)}] Fix integers  $p\geq q\geq 1 $ as well as one of the division algebras
 $\mathbb F:=\mathbb R, \mathbb C,$ or quaternions $ \mathbb H$.
Take $G:=(M_{p,q}(\mathbb F),+)$ as additive group of of $p\times q$ matrices over $\mathbb F$ on which the 
 unitary group $H:= U_p(\mathbb F)$ acts from the left.
$G$ is a real Euclidean vector space of dimension $dpq$ with real scalar product
$(x|y) = \mathfrak R \text{tr}(x^*y)$ where $x^* = \overline x^t$, 
$\mathfrak R t = \frac{1}{2}(t+ \overline t)$ is the real part of $t\in \mathbb F$,  and
$\text{tr}$  the trace in $M_{q,q}(\mathbb F)$.
The  action of $H$ is  orthogonal w.r.t.~this scalar product, and   $x,y\in G$
 are in the same $H$-orbit if and only if $x^*x = y^*y$.
Thus the  space of $H$-orbits is naturally parametrized by the 
cone $D:=\Pi_{q}(\mathbb F)$ of positive semidefinite $q\times q$-matrices over $\mathbb F$.

For $q=1$ and  $\mathbb F= \mathbb R$, we  just have $\Pi_1=[0,\infty[$, and we end up with
 the one-dimensional examples in Section  \ref{examples-orbit-ass}.
For $q\ge2$, the associated orbit hypergroup structures were discussed in \cite{R1} where the associated
multiplicative functions are  Bessel functions of matrix argument.

Similar to  Section  \ref{examples-orbit-ass}, we now fix  $z\in G$, and consider the positive multiplicative function
 $\phi_z(x):=e^{(x|z) }$  on $G$. The associated positive semicharacter $\alpha_z$ on 
 $(D, *)$ can be  written down explicitly in terms of   Bessel functions of matrix argument.
The associated deformed CAS $(M_{p,q}(\mathbb F), \Pi_{q}(\mathbb F), \tilde K)$ may now be written down explicitly.

\item[\rm{(3)}] We mention a further example.  Fix an integer  $ q\geq 1 $ as well as $\mathbb F$ as above.
Let  $G:=(H_{q}(\mathbb F),+)$ be the vector space of all  $\mathbb F$-hermitian $q\times q$-matrices  on which 
the  unitary group $H:= U_p(\mathbb F)$ acts by conjugation. Here, two matrices   $x,y\in G$
 are in the same $H$-orbit if and only if  $x$ and $y$ have the same (ordered) spectrum, i.e., we may identify
 the  space of $H$-orbits with the Weyl chamber 
$$C_q:=\{(x_1,\ldots,x_q)\in\mathbb R^q: \> x_1\ge\ldots\ge x_q\}$$
 of type $A$.
 Again we may write down the multiplicative functions  $\phi_z$ on $G$ explicitely,
 where the associated positive multiplicative 
functions $\alpha_z$ on $(D=C_q, *)$ are Bessel functions of type A.
\end{enumerate}
\end{examples}

\begin{remark}\label{counterex-72}
The example in \ref{examples-orbit-ass2}(1) shows that the condition
 ${\bf 1}\in\tilde S$ in Corollary \ref{COR-MAIN-1} is necessary. To explain this,
 define $(X,D,K)$ as the orbit CAS from \ref{examples-orbit-ass}(1) for $d=1$. Then (T2) holds
for $(X,D,K)$.
 Now let  $(X,D,\tilde K)$ be the deformation of this CAS considered in  \ref{examples-orbit-ass}(2).
If  Corollary \ref{COR-MAIN-1} would be correct, we would find some $\mu\in M^1([0,\infty[)$ with
$$\cos x= \int_0^\infty \frac{\cos(\lambda x)}{\cosh x}\> d\mu(\lambda )
\quad\quad\text{for all }\quad x\in [0,\infty[.$$
It we write the factor $\cosh x$ on the left hand side, it becomes clear that such a measure $\mu$
does not exist.
\end{remark}

\begin{remark}\label{example-unbounded-harmonic} 
Let $d\ge1$ be an integer. Consider the orbit
 CAS $(X=\mathbb R^d, D=[0,\infty[, K)$  from \ref{examples-orbit-ass}.
Fix vectors  $z_1,z_2\in i\cdot\mathbb R^d$ with
$\|z_1\|_2=\|z_2\|_2=1$ and $z_1\ne z_2$. Consider the associated multiplicative pairs
 $(\alpha_{z_j},\phi_{z_j})$ ($j=1,2)$ with $\phi_{z_j}:=e^{-i\langle z_j,x\rangle}$ 
and $\alpha_{z_1}=\alpha_{z_2}$. 

Now consider the deformation $(X,D,\tilde K)$ of $(X,D, K)$ associated with $(\alpha_{z_1},\phi_{z_1})$.
Then, by Lemma \ref{mod-semichar}, 
$( \phi_{z_2}/\phi_{z_1} ,   \alpha_{z_2}/\alpha_{z_1}={\bf 1})$ is a positive semicharacter of $(X,D,K)$, i.e.,
$\phi_{z_2}/\phi_{z_1}$ is a positive, unbounded harmonic function of  $(X,D,\tilde K)$.

This construction of  positive, unbounded harmonic functions 
can be extended to other classes of examples like the discrete ones in the next section.
\end{remark}

\section{Examples associated with infinite distance-transitive graphs}

The set of all infinite distance transitive graphs of finite valency
can be parametrized by two parameters as follows by
Macpherson \cite{Mp}.

\begin{cas-dist-gr}
 Let $a,b\ge2$ be integers. Let  $C_b$ the complete undirected graph
graph with $b$
 vertices, i.e.,  all vertices of $C_b$ are connected.
Consider the infinite
   graph $\Gamma:=\Gamma(a,b)$  where
 precisely $a$ copies of the graph $C_b$ are tacked together at each vertex 
in a tree-like way, i.e., there are no other cycles in  $\Gamma$ than
those in a  copy  of $C_b$. For $b=2$, $\Gamma$ is the homogeneous tree of
valency $a$. We denote the natural distance function on $\Gamma$ by $d$.

After drawing a picture, it is  clear
 that the group $G:=Aut(\Gamma)$
of all graph automorphisms acts on $\Gamma$ in a  distance-transitive way, i.e.,
for all $v_1,v_2,v_3,v_4\in\Gamma$ with
$d(v_1,v_3)=d(v_3,v_4)$ there exists  $g\in\Gamma$
 with $g(v_1)=v_3$ and $g(v_2)=v_4$.
 $Aut(\Gamma)$ is a totally disconnected, locally compact group w.r.t.~the
 topology of pointwise convergence, and
 the stabilizer subgroup $H\subset G$ of any fixed vertex
$e\in\Gamma$ is  compact and open. We  identify $G/H$ with
$\Gamma$, and $G//H$ with $\mathbb N_0$  by distance transitivity.
We now study the association scheme
  $\Lambda=(\Gamma\simeq G/H,\mathbb N_0=G//H, (R_i)_{i\in\mathbb N_0}) $
with  (T1) and (T2)
as well as  the  associated double coset hypergroup $(\mathbb N_0\simeq G//H, *)$.
As in the case of finite  distance-transitive graphs in \cite{BI}, 
 $\Lambda$ and  $(\mathbb N_0, *)$ are symmetric and associated with a
sequence of orthogonal polynomials in the Askey scheme \cite{AW}.

More precisely, it can be seen by some counting (see \cite {V1a}) that the hypergroup convolution is given by
\begin{equation}\label{faltung}
\delta_m*\delta_n = \sum_{k=|m-n|}^{m+n} g_{m,n,k} \delta_k\in  M^1(\mathbb
N_0)
\quad\quad  (m,n\in\mathbb N_0)
\end{equation}
with
$$ g_{m,n, m+n}= \frac{a-1}{a}>0, \quad
g_{m,n, |m-n|}=  \frac{1}{a(a-1)^{m\wedge n-1} (b-1)^{m\wedge n}}>0,$$
$$g_{m,n,|m-n|+2k+1}= \frac{b-2}{ a(a-1)^{m\wedge n-k-1}(b-1)^{m\wedge n-k}}\ge0
\quad  (k=0,\ldots ,m\wedge n-1), $$
$$g_{m,n,|m-n|+2k+2}= \frac{a-2}{a(a-1)^{m\wedge n-k-1}(b-1)^{m\wedge n-k-1}}\ge0
\quad (k=0,\ldots,m\wedge n-2).$$
The Haar weights are  given by $\omega_0:=1$, 
$\omega_n=a(a-1)^{n-1}(b-1)^n \quad (n\ge1)$.
Using
$$g_{n,1,n+1}=\frac{a-1}{a }, \quad 
g_{n,1,n}=  \frac{b-2}{a(b-1)}, \quad
 g_{n,1,n-1}=\frac{1}{a(b-1)},$$
we define a sequence of orthogonal polynomials
$(P_n^{(a,b)})_{n\ge0}$ by  
 $$P_0^{(a,b)}:=1, \quad\quad 
P_1^{(a,b)}(x):= \frac{2}{a}\cdot\sqrt{\frac{a-1}{b-1}}\cdot x +
\frac{b-2}{a(b-1)},$$
and the three-term-recurrence relation
\begin{equation}\label{recu}
P_1^{(a,b)}P_n^{(a,b)}= \frac{1}{a(b-1)}P_{n-1}^{(a,b)}
 + \frac{b-2}{a(b-1)} P_n^{(a,b)} +
\frac{a-1}{ a }P_{n+1}^{(a,b)} \quad\quad(n\ge1) .  
\end{equation}
 Then, 
\begin{equation}\label{prodcartier}\textstyle
 P_m^{(a,b)}P_n^{(a,b)}= \sum_{k=|m-n|}^{m+n} g_{m,n,k}P_k^{(a,b)} \quad (m,n\ge0). 
\end{equation}

We discuss some  properties of the  $P_n^{(a,b)}$ from 
\cite{V1a}, \cite{V4}. 
Eq.~(\ref{recu}) yields
\begin{equation}\label{potenz}
P_n^{(a,b)}\bigl(\frac{z+z^{-1}}{2}\bigr)= \frac{c(z)z^n +c(z^{-1})z^{-n}}
{((a-1)(b-1))^{n/2}} \quad\quad\text{for}\>\> z\in\mathbb C\setminus\{0, \pm 1\}
\end{equation}
with
\begin{equation}\textstyle
c(z):=\frac{(a-1)z -z^{-1} +(b-2)(a-1)^{1/2}(b-1)^{-1/2}}
{a(z-z^{-1})}.
\end{equation}
We  define
\begin{equation}\label{xnull}
 s_0:=  s_0^{(a,b)}:= \frac{2-a-b}{ 2\sqrt{(a-1)(b-1)}}, 
\quad\quad
s_1:= s_1^{(a,b)}:=\frac{ ab -a-b+2}{ 2\sqrt{(a-1)(b-1)}}.
\end{equation}
Then
\begin{equation}\label{reinerpofall}
P_n^{(a,b)}\left(s_1\right)=1,\quad\quad
 P_n^{(a,b)}\left(s_0 \right)=(1-b)^{-n} \quad\quad(n\ge0).
\end{equation}
It is shown in \cite{V4}
 that the
 $P_n^{(a,b)}$ fit into the Askey-Wilson
 scheme (pp.~26--28 of 
 \cite{AW}).
 By the orthogonality relations in \cite{AW}, the 
normalized orthogonality measure $\rho=\rho^{(a,b)}\in M^1(\mathbb R)$
 is 
\begin{equation}\label{orthohne}
d\rho^{(a,b)}(x)= w^{(a,b)}(x)dx\Bigr|_{[-1,1]} 
\quad\quad{\rm for}\quad a\ge b\ge 2
\end{equation}
and 
\begin{equation}\label{orthmit}
d\rho^{(a,b)}(x)=w^{(a,b)}(x)
dx\Bigr|_{[-1,1]} + \frac{b-a}{b} d\delta_{s_0}
  \quad{\rm for}\quad b> a\ge 2
\end{equation}
with
$$w^{(a,b)}(x):=\frac{a}{2\pi} \cdot  \frac{(1-x^2)^{1/2}}{(s_1-x)(x-s_0)}.$$
For $a,b\in\mathbb R$ with $a,b\ge 2$, the numbers $s_0,s_1$ satisfy 
$$-s_1\le s_0\le -1 < 1\le s_1.$$
By
Eq.~(\ref{potenz}), we have the 
dual space
 $$\hat D \simeq
\{x\in\mathbb R:\> (P_n^{(a,b)}(x))_{n\ge0} \quad\text{ is 
 bounded} \}= [-s_1, s_1].$$
This interval  contains the support 
\begin{equation}\label{char-supp-pi}
S:=supp\> \rho^{(a,b)}=
\left\{\begin{array}{r@{\quad\quad}l}
[-1,1] & \text{for}\quad a\ge b\ge 2\\
\{s_0\}\cup [-1,1] & \text{for}\quad b>a\ge2
\end{array}
\right.
\end{equation}
of the orthogonality measure, which is also  the Plancherel measure.; see \cite{La}.
We have $S=\hat D $ precisely for $a=b=2$.
The following theorem from \cite{V4} shows that for these examples several interesting
 phenomena appear, and that  Theorem \ref{main-cont} 
cannot be extended considerably from $S$ to a bigger subset of $\hat D$.
\end{cas-dist-gr}

\begin{theorem}\label{podi} In the setting above the following
 statements are equivalent for $x\in\mathbb R$:
\begin{enumerate}
\item[\rm{(1)}]  $x\in [s_0, s_1]$.
\item[\rm{(2)}] The kernel
$\Gamma\times\Gamma\to\mathbb R$, 
 $(v_1,v_2)\longmapsto P_{d(v_1,v_2)}^{( a, b)}(x)$ is 
 positive definite;
\item[\rm{(3)}] The mapping $g \longmapsto  P_{d(gH,e)}^{( a, b)}(x)$
is  positive definite on $G$.
\end{enumerate}
Moreover, for all $x,y\in[s_0, s_1]$
 there exists a unique 
 $\mu_{x,y}\in M^1([-  s_1,  s_1])$
 with
\begin{equation}\label{special-dist-prof}
P_n^{(a, b)}(x)\cdot P_n^{(a, b)}(y)= \int_{- s_1}^{s_1}
P_n^{(a, b)}(z )\> d\mu_{x,y}(z) \quad { for\>\> all}
\quad n\in\mathbb N_0.\end{equation}
Finally, there are $b>a$ and
$x,y\in[-s_1^{(a, b)},s_0^{(a, b)}[$ for which no 
$\mu_{x,y}\in M^1(\mathbb R)$ exists with
(\ref{special-dist-prof}).
\end{theorem} 

We next construct examples of positive semicharacters of
 $\Lambda$ and study the associated
deformed CAS. The approach  will be similar to  \cite{Voit17} for homogeneous trees.
However, we shall use the results of Section 8 which will simplify some computations.
 
\begin{exposmult}
Fix some constant $c\in\mathbb R$ as well as some point $B$ in the boundary
$\partial\Gamma$, i.e., $B$ is a sequence $(v_n)_{n\in\mathbb N_0}\subset
\Gamma$
of vertices with $d(v_{n+m}, v_n)=m$  for $n,m\in\mathbb N_0$ where $v_0$ 
is the vertex above which is stabilized by $H$.
We  define some kind of  ``distance'' $d(v,B)\in\mathbb Z$ of a vertex $v\in\Gamma$ from $B$
as follows: For $v\in\Gamma$ there is a unique index $n_0\in\mathbb N_0$ such that
$d(v,v_n)\in\mathbb N_0$ is minimal for $n=n_0$. We then put
 $d(v,B):=d(v,v_{n_0})-n_0$. We  in particular have the normalization  $d(v_0,B)=0$.
We now define the function $\phi:=\phi_{B,c}:\Gamma\to]0,\infty[$ with
$$\phi(v):=e^{c \cdot d(v,B)}.$$
\end{exposmult}

\begin{proposition} $\phi$ is a joint eigenfunction of all transition operators $T_h$, $h\in\mathbb N_0$.
More precisely, for all $v\in\Gamma$,
\begin{equation}\label{mult-dist}
T_h\phi(v)=\frac{1}{|\{w\in\Gamma:\> d(v,w)=h\}|} \sum_{w\in\Gamma:\> d(v,w)=h} \phi(w)=
P_h^{(a,b)}(x_c)\cdot \phi(v)
\end{equation}
with 
$$x_c:=\frac{1}{2}\Bigl( e^{c}\sqrt{(a-1)(b-1)} +\frac{1}{e^{c}\sqrt{(a-1)(b-1)}}\Bigr)\in [1,\infty[.$$
\end{proposition}

\begin{proof} The assertion is trivial for $h=0$. 

Assume now that $h\ge1$. We first observe by counting that 
 \begin{equation}\label{number-neighbor}  
 |S(v,h)| =a(a-1)^{h-1}(b-1)^h  \quad\quad\text{for}\quad\quad
S(v,h):=\{w\in\Gamma:\> d(v,w)=h\}.
\end{equation}
Moreover, again by counting we have the following facts:

\noindent There is 1 vertex  $w\in S(v,h)$ with $\phi(w)=e^{c \cdot d(w,B)}=e^{c \cdot (d(v,B)-h)}$;

\noindent there are  $b-2$ vertices  $w\in S(v,h)$ with $\phi(w)=e^{c \cdot( d(v,B)-h+1)}$;

\noindent there are $(a-2)(b-1)$ vertices  $w\in S(v,h)$ with $\phi(w)=e^{c \cdot( d(v,B)-h+2)}$ and so on.

\noindent In general, we see that for $k=0,1,\ldots,h-1$, there are $(b-2)(a-1)^k(b-1)^k$  vertices  $w\in S(v,h)$  
 with $\phi(w)=e^{c \cdot( d(v,B)-h+2k+1)}$,

\noindent  and  for $k=0,1,\ldots,h-2$,
 there are $(a-2)(b-1)^{k+1}(a-1)^k$ vertices  $w\in S(v,h)$ with  $\phi(w)=e^{c \cdot( d(v,B)-h+2k+2)}$.

\noindent Finally, there are $(a-1)^h(b-1)^h$ vertices  $w\in S(v,h)$ with  $\phi(w)=e^{c \cdot( d(v,B)+h)}$.

If we insert these facts into the definition of $T_h$ and use the formula for a geometric sum twice, we arrive at
$$T_h\phi(v)= \alpha(h,a,b,c)\cdot \phi(v)$$
with
\begin{align}\label{eigenvalue-dist-gr}
\alpha(h,a,b,c):=& \frac{1}{a(a-1)^{h-1}(b-1)^h }\Biggl(
1\cdot e^{-ch} + (a-1)^{h}(b-1)^he^{ch}+   \\
&\quad   +(b-2) e^{-ch}\sum_{k=0}^{h-1}(a-1)^k(b-1)^k e^{(2k+1)c} +\notag  \\
&\quad   + (a-2)(b-1) e^{-ch} \sum_{k=0}^{h-2}(a-1)^k(b-1)^k e^{(2k+2)c}\Biggr)\notag \\
=& \frac{e^{-ch}}{a(a-1)^{h-1}(b-1)^h }\Biggl(
1+ (a-1)^{h}(b-1)^he^{2ch}+ \notag \\
&\quad   + \frac{(b-2)e^{c}\bigl[(a-1)^h(b-1)^he^{2ch}-1\bigr]}{(a-1)(b-1)e^{2c}-1} + \notag \\
&\quad   +(a-2)(b-1) e^{2c} \frac{(a-1)^{h-1}(b-1)^{h-1}e^{2c(h-1)}-1}{(a-1)(b-1)e^{2c}-1}\Biggr).
\notag
\end{align}
In particular,  $\phi$ is a joint  eigenfunction of all $T_h$. We now conclude from Proposition 
\ref{mult-char}  that
the mapping  $\mathbb N_0\to]0,\infty[$, $h\mapsto \alpha(h,a,b,c)$ is multiplicative on the symmetric polynomial
hypergroup $(\mathbb N_0,*)$ which implies that this mapping is a positive semicharacter.  On the other hand, it
follows from the theory of polynomial hypergroups (see \cite{BH}) that the positive semicharacters on 
$(\mathbb N_0,*)$ have the form  $h\mapsto P_h^{(a,b)}(x)$ with some unique $x\in [1,\infty[$. In order to find 
the correct $x$, we compare the explicit representation (\ref{potenz}) of  $ P_h^{(a,b)}(x)$ with the eigenvalues
$\alpha(h,a,b,c)$ in (\ref{eigenvalue-dist-gr}) for large values of $h$. This leads readily to
$$x=x_c= \frac{1}{2}\Bigl( e^{c}\sqrt{(a-1)(b-1)} +\frac{1}{e^{c}\sqrt{(a-1)(b-1)}}\Bigr)$$
and thus $    \alpha(h,a,b,c)=   P_h^{(a,b)}(x_c) $  as claimed.
\end{proof}

We can now use the results of Section 8 to define a deformed CAS $(\Gamma,\mathbb N_0,\tilde K)$
according Proposition \ref{deformed-ass-scheme} with
$\tilde K_{0}(x, \{y\})=\delta_{x,y}$ and, for $h\ge1$,
$$\tilde K_{h}(x, \{y\}):= \frac{e^{c(d(y,B)-d(x,B))}}{P_h^{(a,b)}(x_c)}   K_{h}(x, \{y\})
= \frac{e^{c(d(y,B)-d(x,B))}}{P_h^{(a,b)}(x_c)\cdot a(a-1)^{h-1}(b-1)^h }   \quad\quad (x,y\in X).$$
For $b=2$, i.e., for homogeneous trees, this result fits with the results in \cite{Voit17} 
where the same kernels were obtained in a different, more computational  way.

We point out that in particular the case $x_c=1$ is interesting which appears precisely for
$$ c:=-\frac{1}{2}\ln((a-1)(b-1)).$$
In this case, the associated 
deformed polynomial hypergroup $(\mathbb N_0,\tilde *)$
has the functions
 $$\mathbb N_0\to\mathbb R, \quad 
n\mapsto \widetilde{P_n}^{(a,b)}(x):= P_n^{(a,b)}(x)/ P_n^{(a,b)}(1) \quad(x\in\mathbb R)$$
 as semicharacters.
It can be easily derived from  (\ref{potenz}) that here the dual space $(D,\tilde*)^\wedge$ corresponds to
$S$ from (\ref{char-supp-pi}).

We point out that the ideas of this section may be used to study 
deformations of infinite commuative association schemes associated with affine buildings e.g.~of type $\tilde A_n$.

\begin{remark}\label{example-no-t12}
Consider a homogeneous tree $\Gamma$ of valeny $a$, i.e., we take $b=2$ above.
Fix a point  $B\in\partial \Gamma$ in the boundary 
 and a constant $c\in\mathbb R\setminus \{0\}$ as before, and consider the associated multiplicative pair
 and the associated deformed CAS $(\Gamma,\mathbb N_0, \tilde K)$ as above. Let $v_0\in\Gamma$ as above.
Then, by  Proposition \ref{deformed-ass-scheme}(1), 
an invariant measure $\omega_\Gamma$ of $(\Gamma,\mathbb N_0, \tilde K)$ is given by
 $\omega_\Gamma(v)=e^{2c\cdot d(v,B)}$ for $v\in\Gamma$. Its  push forward
$\pi_{v_0}(\omega_\Gamma)\in M^+(\mathbb N_0)$ then satisfies
$$\pi_{v_0}(\omega_\Gamma)(0)=1, \quad \pi_{v_0}(\omega_\Gamma)(1)= e^{-2c}+(a-1)e^{2c},$$
and for $n\ge2$,
$$ \pi_{v_0}(\omega_\Gamma)(n)= e^{-2cn} +\sum_{l=1}^{n-1}e^{2(-n+2l)}  (a-2)(a-1)^{l-1}+e^{2cn} (a-1).$$
This measure $\pi_{v_0}(\omega_\Gamma)$ is usually not equal to the Haar measure  $\omega_c$ 
of the (deformed) polynomial hypergroup $(\mathbb N_0,*_c)$ with normalization $\omega_c(0)=1$, 
as we have
$$\omega_c(1)=\frac{1}{g_{1,1,0}^c}= \frac{((a-1)e^{2c}+1))^2}{a e^{2c}}$$
which is usually different from  $\pi_{v_0}(\omega_\Gamma)(1)$ above.
We thus in particular conclude from Lemmas \ref{proj-haar-m-T1}(1) and \ref{proj-haar-m-T2}(1)
that for $c\ne0$, the  deformed CAS $(\Gamma,\mathbb N_0, \tilde K)$ usually do not have the properties (T1) and (T2).
\end{remark}

\section{Further constructions of continuous association schemes}

In this section we  present further  constructions of CAS from  given ones.
We start with the direct product.

\begin{directp}
Let $(X_1,D_1, K^1)$ and  $(X_2,D_2, K^2)$ be CAS with associated hypergroups $(D_1,*_1)$ and  $(D_2,*_2)$.
 We then form $X:=X_1 \times X_2$ and   $D:=D_1 \times D_2$. We define the convolution of point measures
 via the direct product
of measures by
$$\delta_{(x_1,x_2)}*\delta_{(y_1,y_2)}:= (\delta_{x_1}*_1\delta_{y_1})\times (\delta_{x_2}*_2\delta_{y_2})
\quad\quad(x_1,y_1\in D_1, x_2,y_2\in D_2).$$
 It is well-known  that the unique bilinear, 
weakly continuous extension of this convolution leads to the so-called direct product hypergroup $(D,*)$; see
Section 10.5 of \cite{J} or Section 1.5.28 of \cite{BH}.
We now put
\begin{equation}\label{prod-kernel}
K_{(h_1,h_2)}((x_1,x_2), A_1\times A_2):= K_{h_1}^1(x_1,A_1)\cdot K_{h_2}^2(x_2,A_2)
\end{equation}
for $h_i\in D_i$, $x_i\in X_i$, and Borel sets $A_i\subset X_i$ with $i=1,2$. It is well-known and can be 
 easily checked that (\ref{prod-kernel}) leads to a unique Markov kernel $K_{(h_1,h_2)}$ on $X$
for $(h_1,h_2)\in D_1\times D_2$. Moreover, it can be easily seen that these kernels can be combined 
 to a  Markov kernel $K$
 from $X\times D$ to $X$.
We have:
\end{directp}

\begin{proposition}
$(X:=X_1 \times X_2, D:=D_1 \times D_2, K)$ is a CAS; it will be called the direct
 product of $(X_1,D_1, K^1)$ and  $(X_2,D_2, K^2)$.

If $(X_1,D_1, K^1)$ and  $(X_2,D_2, K^2)$ are commutative, symmetric, compact, or discrete, then so is $(X,D,K)$.
 Moreover, the properties (T1) and (T2) are also preserved.
\end{proposition}

\begin{proof} First of all, one has to check that $K$ is a continuous Markov kernel
 in the notion of the beginning of Section 4. 
For this we first notice that for $g\in C_c(X)$ of the form $g(x_1,x_2)=g_1(x_1)g_2(x_2)$ with $g_i\in C_c(X_i)$,
the map $((x_1,x_2),h)\mapsto T_h(g)(x_1,x_2)$ is continuous by the product structure in  (\ref{prod-kernel}).
As by the theorem of Stone-Weierstrass, the linear span of functions on $X$ of the form  $g(x_1,x_2)=g_1(x_1)g_2(x_2)$ 
for $g_i\in C_c(X_i)$ is $\|.\|_\infty$-dense in $C_c(X)$, 
the map above is continuous even for all  $g\in C_c(X)$. Taking Lemma \ref{compact-union} into account, we see that
the map above is continuous even for all  $g\in C_b(X)$.

 All properties in \ref{central-def}(1) and (2) can be checked easily. We only mention that for the projections $\pi_1$
and $\pi_2$ of the given CAS, the projection $\pi:X\times X\to D$ satisfies
 $\pi((x_1,x_2),(y_1,y_2))=(\pi_1(x_1,y_1),\pi_2(x_2,y_2))$ by  (\ref{prod-kernel}). 
 \ref{central-def}(3) can be also checked easily by the product structure in  (\ref{prod-kernel}).
 Moreover, if we define the  product measure $\omega_X:=\omega_{X_1}\times\omega_{X_2}$,
then Eq.~(\ref{ess-hg-eq}) in \ref{central-def}(4) can be also checked easily.

The statement about  $(X,D,K)$ being commutative, symmetric, compact, or discrete is also trivial.

We next check property (T2). We here first notice that for Haar measures $\omega_{D_i}$ of $(D_i,*_i)$ ($i=1,2$), the product
 $\omega_D:=\omega_{D_1}\times\omega_{D_2}$ is a Haar measure of $(D,*)$. Using (T2) for the given schemes, we readily obtain 
for $f_i\in C_c(D_i)$, $g_i\in C_c(X_i)$ and $x_i\in X_i$ ($i=1,2$) that
\begin{align}\label{prod-kernel2}
\int_X f_1(\pi_1(x_1,x_2))&f_2(\pi_2(y_1,y_2))\>  g(y_1)g(y_2) \> d(\omega_{X_1}\times\omega_{X_2})(y_1,y_2)=\notag\\
&=\int_D\int_X g(y_1)g(y_2)\> K_{h_1,h_2}((x_1,x_2), d(y_1,y_2)) \> f_1(h_1)f_2(h_2) d\omega_D(h_1,h_2).
\end{align}
Again, the theorem of Stone-Weierstrass shows that
$T_fg=T^{f\circ\pi}g$ holds for $f\in C_c(D)$ and $g\in C_c(X)$ as claimed.

Property (T1) can be derived in the same way by a Stone-Weierstrass argument.
\end{proof}

We next turn to joins.

\begin{joinconst}
We first recapitulate the join of hypergroups from Section 10.5 of \cite{J}.
Let $(D_1,*_1)$ be a discrete hypergroup with identity $e_1\in D_1$, and
 let  $(D_2,*_2)$ be a compact hypergroup
with normalized Haar measure $\omega_{D_2}$. We  form the disjoint union 
 $$D:=D_1\vee D_2:=(D_1\setminus\{e_1\})\cup D_2$$
with
$D_1\setminus\{e_1\}$ and $ D_2$ as open subsets. $D$  is locally compact.
 We  define the convolution of point measures
on $D$ via
$$\delta_x*\delta_y:= \left\{\begin{array}{r@{\quad\quad}l}
\delta_x*_2\delta_y & \text{for}\quad x,y\in D_2\\
(\delta_x*_1\delta_y)|_{D_1\setminus\{e_1\}}+ (\delta_x*_1\delta_y)(\{e_1\})\omega_2 & \text{for}\quad x,y\in D_1  \\
\delta_x & \text{for}\quad x\in D_1,y\in D_2 \\
\delta_y & \text{for}\quad y\in D_1,x\in D_2 
\end{array}
\right.
$$
Assume now that $(D_1,*_1)$ is associated with some discrete CAS $(X_1,D_1, K^1)$
and $(D_2,*_2)$ with some compact CAS $(X_2,D_2, K^2)$. 
We assume that  $\omega_{D_2}$ and $\omega_{X_2}$ are normalized such that they are both probability measures.
We form the join $(D,*)$ as above and put $X:=X_1\times X_2$.
Moreover, for $h\in D$, $x_1\in X_1$, $x_2\in X_2$, Borel sets $B\subset X_2$ and sets $A\subset X_1$ we put
\begin{equation}\label{prod-kernel-join}
K_{h}((x_1,x_2), A\times B):=  \left\{\begin{array}{r@{\quad\quad}l}
K_h^2(x_2,B)\cdot \delta_{x_1}(A) & \text{for}\quad h\in D_2\\
\omega_{X_2}(B)\cdot K_{h}^1(x_1,A) & \text{for}\quad h\in D_1\setminus\{e_1\}
\end{array}
\right.
\end{equation}
Clearly, each $K_h$ can be extended uniquely to a Markov kernel on $X$, 
and we can combine the $K_h$ to some Markov kernel $K$ from $X\times D$ to $X$.
\end{joinconst}

\begin{proposition}\label{join-def}
$(X, D, K)$ is a CAS which will be called the join of the CAS  $(X_1,D_1, K^1)$ and  $(X_2,D_2, K^2)$.
We shall write the join as $(X_1,D_1, K^1)\vee (X_2,D_2, K^2)$.

If $(X_1,D_1, K^1)$ and  $(X_2,D_2, K^2)$ are commutative or symmetric, then so is $(X,D,K)$.
 Moreover,  if  $(X_1,D_1, K^1)$ has property (T2), then  so has $(X,D,K)$.
\end{proposition}

\begin{proof}
By the topological structure of $X$ and $D$ it is clear
 that our kernel $K$  from $X\times D$ to $X$ is continuous.
Moreover, the properties \ref{central-def}(1) and (2) are obviously satisfied
with the  projection $\pi:X\times X\to D$ with
$$\pi((x_1,x_2),(y_1,y_2)):= \left\{\begin{array}{r@{\quad\quad}l}
\pi_2(x_2,y_2)\in D_2\subset D &\text{for}\quad x_1=y_1\\
\pi_1(x_1,y_1)\in D_1\setminus\{e_1\}\subset D &\text{for}\quad x_1\ne y_1
\end{array}
\right.$$
To check \ref{central-def}(3), we fix
 $x_i\in X_i$ and  Borel sets  $A_i\subset X_i$ ($i=1,2$). 
Then for $h_1,h_2\in D_2$,
\begin{align}
K_{h_1}\circ K_{h_2}&((x_1,x_2), A_1\times A_2)= \delta_{x_1}(A_1)\cdot K_{h_1}^2\circ K_{h_2}^2(x_2, A_2)
\notag\\
&=\delta_{x_1}(A_1)\cdot\int_{D_2}  K_{h}^2(x_2, A_2)\> d(\delta_{h_1}*_2\delta_{h_2})(h)
\notag\\
&=\int_{D_2}  K_{h}((x_1,x_2), A_1\times A_2)\> d(\delta_{h_1}*_2\delta_{h_2})(h)
\notag\end{align}
as claimed. Moreover, for  $h_1\in D_1$, $h_2\in D_2$, 
\begin{align}
K_{h_2}\circ K_{h_1}&((x_1,x_2), A_1\times A_2)= \int_{X_2}K_{h_1}((x_1,y_2), A_1\times A_2)\> K_{h_2}^2(x_2,dy_2)
\notag\\
&=\int_{X_2} \omega_{X_2}(A_2) K_{h_1}^1(x_1,A_1)\> K_{h_2}^2(x_2,dy_2) 
= K_{h_1}^1(x_1,A_1) \omega_{X_2}(A_2)\notag\\
&=K_{h_1}((x_1,x_2), A_1\times A_2)
\notag\end{align}
as claimed. Moreover, for  $h_1\in D_1$, $h_2\in D_2$ we conclude from Eq.~(\ref{invariant-omega}) that
\begin{align}
K_{h_1}\circ K_{h_2}&((x_1,x_2), A_1\times A_2)=
\int_{X_1}\int_{X_2} K_{h_2}^2(y_2,A_2) \delta_{y_1}(A_1)\> d\omega_{X_2}(y_2)\> K_{h_1}^1(x_1,dy_1)
\notag\\
&=K_{h_1}^1(x_1,A_1)\cdot  \omega_{X_2}(A_2) = K_{h_1}((x_1,x_2), A_1\times A_2)
\notag\end{align}
as claimed. Finally, for  $h_1,h_2\in D_1$ we conclude from Lemma \ref{unique-invar-measure}(3) with 
our normalizations that $\int_{D_2} K_{h}^2(x_2,A_2)\> d\omega_{D_2}(h)=\omega_{X_2}(A_2)$ and thus
\begin{align}
K_{h_1}\circ K_{h_2}&((x_1,x_2), A_1\times A_2)=
\int_{X_1}  \omega_{X_2}(A_2)    K_{h_2}^1(y_1,A_1) \> K_{h_1}^1(x_1,dy_1)\cdot \omega_{X_2}(X_2)
\notag\\
&= \int_{D_1} K_h^1(x_1, A_1)\> d(\delta_{h_1}*_1\delta_{h_2})(h)\cdot\omega_{X_2}(A_2)
\notag\\
&=\int_{D_1\setminus \{e_1\}}\omega_{X_2}(A_2)K_h^1(x_1, A_1)\>d(\delta_{h_1}*_1\delta_{h_2})(h) +\notag\\
& \quad\quad\quad\quad\quad+
(\delta_{h_1}*_1\delta_{h_2})( \{e_1\}) \delta_{x_1}(A_1) \int_{D_2} K_{h}^2(x_2,A_2)\> d\omega_{D_2}(h)
\notag\\
&=\int_{D} K_h((x_1,x_2), A_1\times A_2)\>  d(\delta_{h_1}*\delta_{h_2})(h)
\notag\end{align}
which completes the proof of \ref{central-def}(3). 

In order to check the adjoint relation \ref{central-def}(4), we put
$\omega_X:=\omega_{X_1}\times\omega_{X_2}$. Let  $f_1,g_1\in C_c(X_1)$, $f_2,g_2\in C_c(X_2)$ and consider
 $f,g\in C_c(X)$ with  $f(x_1,x_2):=f_1(x_1)f_2(x_2)$ and $g(x_1,x_2):=g_1(x_1)g_2(x_2)$.
 We conclude from (\ref{prod-kernel-join}) that then
for all  $h\in D$ and  $(x_1,x_2)\in X$, 
\begin{equation}
T_{h}f(x_1,x_2)=  \left\{\begin{array}{r@{\quad\quad}l}
f_1(x_1)\cdot T_h^2f_2(x_2) & \text{for}\quad h\in D_2\\
 T_h^1f_1(x_1) \cdot\int_{X_2} f_2\> d\omega_{X_2}
 &  \text{for}\quad h\in D_1\setminus\{e_1\}
\end{array}
\right.
\end{equation}
and hence
$$
\int_X T_{h}f \cdot g\> d\omega_X =  \left\{\begin{array}{r@{\quad\quad}l}
\int_{X_1}  f_1g_1\>  d\omega_{X_1} \cdot \int_{X_2} T_{h}^2f_2\cdot g_2\>d\omega_{X_2}
 & \text{for}\quad h\in D_2\\
 \int_{X_2} f_2 d\omega_{X_2} \cdot\int_{X_2} g_2 d\omega_{X_2} \cdot \int_{X_1} T_{h}^1f_1\cdot g_1\>d\omega_{X_1}
 &  \text{for}\quad h\in D_1\setminus\{e_1\}
\end{array}
\right.$$
This leads to the adjoint relation \ref{central-def}(4) for our particular functions $f,g$.
The assertion for general $f,g\in C_c(X)$ finally follows from a Stone-Weierstrass argument.
This completes the proof of the first part of the proposition.

Clearly, commutativity and symmetry are preserved under joins.

We next turn to property (T2). Recapitulate that the compact CAS $(X,D,K)$ has (T2) by Proposition \ref{compact-t2},
and  that
$$\omega_D:= \omega_{D_1}(\{e_1\}) \cdot \omega_{D_2}+ \omega_{D_1}|_{D_1\setminus\{e_1\}}\in M^+(D)$$
is a Haar measure of the join $(D,*)$ with the normalization  $\omega_{D_2}\in M^1(D_2)$.
We  check (T2) via Lemma \ref{t2-reformulation}. For this fix $(x_1,x_2)\in X$ and Borel sets
$A_1\subset X_1$, $A_2\subset X_2$. Then, by \ref{t2-reformulation},
\begin{align}
\int_D& K_h((x_1,x_2), A_1\times A_2)\> d\omega_D(h)=
\int_{D_2} \ldots  d\omega_D(h) + \int_{D_1\setminus\{e_1\}} \ldots  d\omega_D(h)
\notag\\
&= \omega_{D_1}(\{e_1\})\> \delta_{x_1}(A_1) \int_{D_2} K_h^2(x_2,A_2) \> d\omega_{D_2}(h)
+
\notag\\
&\quad\quad\quad\quad+
\omega_{X_2}(A_2)
\Bigl[ \int_{D_1}K_h^1(x_1,A_1) \> d\omega_{D_1}(h)-\omega_{D_1}(\{e_1\})\> \delta_{x_1}(A_1)\Bigr]
\notag\\
&=\omega_{D_1}(\{e_1\})\> \delta_{x_1}(A_1)\omega_{X_2}(A_2)+\omega_{X_2}(A_2)\Bigl[
\omega_{X_1}(A_1)-\omega_{D_1}(\{e_1\})\> \delta_{x_1}(A_1)\Bigr]
\notag\\
&=\omega_{X_1}(A_1)\omega_{X_2}(A_2)=\omega_X( A_1\times A_2).
\notag
\end{align}
As the Borel measures $\int_D K_h((x_1,x_2), .)\> d\omega_D(h) $ and
 $\omega_X$ on $X_1\times X_2$ are determined uniquely by their values on cylinder sets,
 it follows that both measures are equal independent of $(x_1,x_2)$.
  (T2) now follows from  Lemma \ref{t2-reformulation}.
\end{proof}

We now consider iterated joins of finite CAS. For this we fix a sequence $(\Lambda_n)_{n\in\mathbb N}$
of finite CAS. We form the iterated joins
$$ \Lambda_n^i:= \Lambda_n\vee(\ldots(\Lambda_3\vee(\Lambda_2\vee\Lambda_1))\ldots) \quad\quad(n\in\mathbb N)$$
and
$$ \Lambda_n^p:= (\ldots((\Lambda_1 \vee \Lambda_2)\vee \Lambda_3)\ldots)\vee\Lambda_n \quad\quad(n\in\mathbb N).$$
We now form the inductive limit of the sequence $(\Lambda_n^i)_{n\in\mathbb N}$
of
discrete CAS  as well as the  projective  limit of  the sequence $(\Lambda_n^p)_{n\in\mathbb N}$
of compact CAS  in an informal way and obtain as limits a discrete CAS $ \Lambda^i$ and a compact 
 CAS $\Lambda^p$. We here do not work out any theory of these limits which are well-defined on the level of hypergroups; see 
\cite{V3aa}. We here just present these limits as examples of CAS.
We start with the inductive limit:

\begin{example} 
We start with a sequence of finite CAS $(\Lambda_n=(X_n,D_n,K^n))_{n\in\mathbb N}$ with the associated hypergroups $(D_n,*_n)$ with the identities $e_n\in D_n$
and with the normalized Haar measures $\omega_{D_n}\in M^1(D_n)$, and with the normalized adjoint measures  $\omega_{X_n}\in M^1(X_n)$.
We also fix some sequence $(z_n)_{n\in\mathbb N}$ with $z_n\in X_n$ for each $n$. Assume that $|D_n|\ge 2$ for all 
$n$.
 We now define the discrete inductive limit  CAS $ \Lambda=(X,D,K)$ as follows:
$D$ is the discrete, disjoint union
$$D:=D_1 \cup \bigcup_{n\ge2} (D_n\setminus \{e_n\}),$$
and $X$ is the discrete, countable set
$$X:=\{ (x_n)_{n\in\mathbb N}:\> x_n\in X_n \>\>\text{for all}\>\> n, \>\>\text{and}\>\> x_n=z_n \>\>\text{for all except for at most finitely many}\>\>
n\}.$$
The convolution of point measures on the hypergroup $(D,*)$ is given by
$$\delta_h*\delta_l:=\left\{\begin{array}{r@{\quad\quad}l}
\delta_l &\text{for}\quad  l\in D_n, \quad h\in D_1\cup\bigcup_{k=1}^{n-1}D_k\setminus\{e_k\},\quad n\ge2\\
\delta_h &\text{for}\quad  h\in D_n, \quad l\in D_1\cup\bigcup_{k=1}^{n-1}D_k\setminus\{e_k\},\quad n\ge2\\
\delta_h*_1\delta_l  &\text{for}\quad  h,l\in D_1\end{array}
\right.$$
and, for $h,l\in D_n\setminus\{e_n\}$ with $n\ge2$, by the probability measure
\begin{align}
\delta_h*\delta_l:=&(\delta_h*_n\delta_l)|_{D_n\setminus\{e_n\}} + \notag\\
+&(\delta_h*_n\delta_l)(\{e_n\})\Bigl[
\sum_{k=2}^{n-1}\Bigl( \prod_{m=k+1}^{n-1}\omega_{D_m}(\{e_m\})\Bigr)
\cdot \omega_{D_k}|_{D_k\setminus\{e_k\}}
+\Bigl( \prod_{m=2}^{n-1}\omega_{D_m}(\{e_m\})\Bigr)\cdot \omega_{D_1}\Bigr].\notag
\end{align}
A Haar measure on $(D,*)$ is given by the measure
$$\omega_D:= \omega_{D_1}+\sum_{n=2}^\infty \Bigl( \prod_{m=1}^{n}\omega_{D_m}(\{e_m\})\Bigr)^{-1}\cdot
 \omega_{D_n}|_{D_n\setminus\{e_n\}}.$$

The kernels $K_h$ on $X$ are given by
$$K_h((x_n)_{n\in\mathbb N}, \{(y_n)_{n\in\mathbb N}\})=
\left\{\begin{array}{r@{\>\>}l}
K_h^1(x_1, \{y_1\})\cdot \delta_{(x_n)_{n\ge2},(y_n)_{n\ge2}} &\text{for}\> h\in D_1\\
\prod\limits_{n=1}^{l-1}\omega_{X_n}(\{y_n\}) \>
K_h^l(x_l, \{y_l\})\> \delta_{(x_n)_{n\ge l+1},(y_n)_{n\ge l+1}}&\text{for}\> h\in D_l\setminus\{e_l\},\\
 & \> l\ge 2
\end{array}
\right.$$
  where we again use the Kronecker-$\delta$. 
These $K_h$ are in fact Markov kernels on $X$ where
 the properties \ref{central-def}(1) and (2) obviously hold
with the  projection $\pi$ with
$$\pi((x_n)_{n\in\mathbb N},(y_n)_{n\in\mathbb N} )= 
\left\{\begin{array}{r@{\quad\quad}l}
\pi_1(x_1,y_1)\in D_1\subset D &\text{for}\quad (x_n)_{n\ge2}=(y_n)_{n\ge2}\\
\pi_n(x_n,y_n)\in D_n\setminus\{e_n\}\subset D &\text{for}\>\>  n:= \max\{l:  x_l\ne y_l\}\ge 2
\end{array}
\right.$$
Notice that the maximum exists by the definition of $X$.
 \ref{central-def}(3) can now be  checked by using the computations in the proof of Proposition \ref{join-def}.
 These computations also show that 
$$\omega_X(\{(x_n)_{n\in\mathbb N}\}):=\prod_{n=1}^{l}\frac{\omega_{X_n}(\{x_n\})}{\omega_{X_n}(\{z_n\})}
 \quad\text{for}\quad l \quad\text{with}\quad x_n=z_n\quad\text{for}\quad n>l $$
is an adjoint measure, and that (T2) holds by Lemma \ref{t2-reformulation}. We omit the details of proofs.
\end{example}

We next turn to the projective  limit:

\begin{example}\label{ex-projective} 
We start with the same sequence $(\Lambda_n=(X_n,D_n,K^n))_{n\in\mathbb N}$ with the  notations as before.
We define the  compact projective limit  CAS $ \Lambda=(X,D,K)$ as follows:
$$D:= \bigcup_{n\ge1} (D_n\setminus \{e_n\}) \cup \{e\}$$
is the one-point compactification of the discrete disjoint union $\bigcup_{n\ge1} (D_n\setminus \{e_n\})$
where the additional non-discrete limit point $e$  will be the neutral element $e$ of $(D,*)$.
The hypergroup convolution of point measures is given by
$$\delta_h*\delta_l:=\left\{\begin{array}{r@{\quad\quad}l}
\delta_l &\text{for}\quad  l\in D_n\setminus\{e_n\}, \quad h\in \{e\}\cup\bigcup_{k>n} (D_k\setminus\{e_k\})\\
\delta_h &\text{for}\quad  h\in D_n\setminus\{e_n\}, \quad l\in \{e\}\cup\bigcup_{k>n} ( D_k\setminus\{e_k\})
\end{array}
\right.$$
and, for $h,l\in D_n\setminus\{e_n\}$  by the probability measure
$$
\delta_h*\delta_l:=(\delta_h*_n\delta_l)|_{D_n\setminus\{e_n\}} + (\delta_h*_n\delta_l)(\{e_n\})
\sum_{k>n}\Bigl( \prod_{m=n+1}^{k-1}\omega_{D_m}(\{e_m\})\Bigr)
\cdot \omega_{D_k}|_{D_k\setminus\{e_k\}}.
$$
 The normalized Haar measure on the compact hypergroup $(D,*)$ is given by 
$$\omega_D:= \sum_{n=1}^\infty \Bigl( \prod_{m=1}^{n-1}\omega_{D_m}(\{e_m\})\Bigr)\cdot
 \omega_{D_n}|_{D_n\setminus\{e_n\}}\in M^1(D).$$
Moreover,  $X$ is the compact usual direct product $\prod_{n\ge1} X_n$,
 and the adjoint measure is the infinite product measure
$$\omega_X:= \prod_{n\ge0} \omega_{X_n}.$$
The kernel $K_e$ on $X$ will be the identity, and for $h\in D_n\setminus \{e_n\}\subset D$ with $n\in \mathbb N$,
we have
$$K_h((x_n)_{n\in\mathbb N}, \{(y_1,\ldots y_n)\}\times A)=
\delta_{(x_1,\ldots,x_{n-1}), (y_1,\ldots,y_{n-1})} \> K_h^n(x_n, \{y_n\}) \>  \prod_{l\ge n+1}\omega_{X_l}(A)$$
for all $(x_n)_{n\in\mathbb N}\in X$, $ (y_1,\ldots y_n)\in X_1\times \ldots\times X_n$ and Borel sets
$A\subset \prod_{l\ge n+1} X_l$.
These formulas clearly determine unique  Markov kernels $K_h$ on $X$ where
 the properties \ref{central-def}(1) and (2) obviously hold
with the  projection $\pi:X\times X\to D$ with
$$\pi((x_n)_{n\in\mathbb N},(x_n)_{n\in\mathbb N} )=e\in D$$ and
$$\pi((x_n)_{n\in\mathbb N},(y_n)_{n\in\mathbb N} )= \pi_n(x_n,y_n)\in D_n\setminus\{e_n\}\subset D$$
for $(x_n)_{n\in\mathbb N}\ne(x_n)_{n\in\mathbb N}  $ and $n:=\min\{l:  x_l\ne y_l\}$.
 \ref{central-def}(3) and (4) can now be  checked by using the computations
 in the proof of Proposition \ref{join-def}. Finally, (T2) is clear by compactness and Proposition 
\ref{compact-t2}.

We notice that for  given  sequences  $(\Lambda_n)_n$ of finite CAS which are not coming from groups,
 the construction  \ref{ex-projective} leads to examples of compact, non-discrete stron  CAS  
which are not associated with groups according
 to Proposition \ref{group-to-cont-asso-scheme}.
\end{example}

There exist several generalizations of the join on the level of hypergroups
 like substitutions of open hypergroups in \cite{V3a} or Section 8.1 of \cite{BH} or
further constructions used in several papers of H.~Heyer, S.~Kawakami, and others; 
see e.g. \cite{HKKK}, \cite{HKY}, and papers cited there. 
We expect that these constructions should also have a meaning on the level of CAS.

\section{Random walks on continuous association schemes}

In this section we  introduce and investigate random walks on $X$ associated with some given CAS 
$(X,D,K)$. Before doing so we briefly recapitulate some facts on random walks on the hypergroup $(D,*)$. 
For simplicity we restrict our attention to the time-homogeneous case.

\begin{random-walk-hy} Let $(D,*)$ be a second countable hypergroup with identity $e$.
 Let either $T:=\mathbb N_0$ or $T:=[0,\infty[$. 
A family  $(\mu_t)_{t\in T}\subset M^1(D)$ of probability measures is called a (discrete or continuous)
convolution semigroup, if $\mu_0=\delta_e$, and if for all $s,t\in T$, $\mu_{s+t}=\mu_s*\mu_t$,
 and if in the continuous case, the mapping $[0,\infty[\to  M^1(D)$, $t\mapsto \mu_t$ is weakly 
continuous in addition.
 It can be easily checked and is well-known that for each $t\in T$ we may form
 the Markov kernel $K_t$ on $D$ via
$$K_t(h,A):=(\delta_x*\mu_t)(A) \quad\quad (h\in D,\> A\in\cal B(D),\> t\in T).$$
The $K_t$ are in fact Feller kernels (see the beginning of Section 4) with $K_e$ as trivial kernel and 
$K_s\circ K_t=K_{s+t}$ for all  $s,t\in T$; see the beginning of Section 4 for the notations. 
In particular, $(K_t)_{t\in T}$ is a semigroup of transition kernels.
As $(D,*)$ is second countable and locally compact,
it is  a matter of fact that for the starting distribution $\delta_e$ and this semigroup
 there exists a time-homogeneous
 Markov process $(Y_t)_{t\in T}$ with the transition probabilities
$$ P(Y_{s+t}\in A| Y_s=h)=K_t(h,A)=(\delta_x*\mu_t)(A) \quad\quad (x\in D, \> A\in \cal B(D), \> s,t\in T).$$
 These Markov processes are called random walks on $(D,*)$ associated with $(\mu_t)_{t\in T}$.

Notice that for a locally compact group $D=G$, this just means that  $(Y_t)_{t\in T}$ just consists
 of group products of
iid 
$G$-valued random variables in the discrete case, and that
  $(Y_t)_{t\in T}$ is a Levy process in the continous case.

For a detailled study of  random walks on hypergroups including
 limit theorems for special classes we refer to \cite{BH} and references there.
\end{random-walk-hy}

We now use these ideas to define random walks  $(Z_t)_{t\in T}$ on $X$ for a CAS $(X,D,K)$:

\begin{random-walk-CAS} Let  $(X,D,K)$ be a CAS with associated hypergroup $(D,*)$ with identity $e$.
 Let $T:=\mathbb N_0$ or $T:=[0,\infty[$, and let $(\mu_t)_{t\in T}\subset M^1(D)$ be 
 a (discrete or continuous)
convolution semigroup of probability measures on $D$ as before.
For each $t\in T$ we now define the Markov kernel $K^X_{\mu_t}$ on $X$ via
$$K^X_{\mu_t}(x,A):= \int_D K_h(x,A)\> d\mu_t(h) \quad\quad (x\in X,\> A\in\cal B(X),\> t\in T)$$
which is associated with the transition operator $T_{\mu_t}$ of Lemma \ref{op-t-mu}, i.e., the 
$K^X_{\mu_t}$ are  in fact Feller kernels by \ref{op-t-mu}. Moreover, precisely as in Proposition 
\ref{alg-hom} we see that $(K^X_{\mu_t})_{t\in T} $ is a  semigroup of transition kernels.
If we now fix some starting point $x_0\in X$, we again find some associated  time-homogeneous
 Markov process $(Z_t)_{t\in T}$. These processes are called random
 walks on $X$ with start in $x_0$ associated with $(\mu_t)_{t\in T}$.
\end{random-walk-CAS}

For  $T:=[0,\infty[$, we show that  $(Z_t)_{t\in T}$ is a so-called Feller process, i.e., that the operators $T_{\mu_t}$
associated with the Feller kernels $K^X_{\mu_t}$ satisfy in addition the condition
\begin{equation}\label{Feller-semigr-cond}
\|T_{\mu_t}g-g\|_\infty \longrightarrow 0 \quad\quad \text{for}\quad\quad t\to0 
\quad\quad \text{and all}\quad\quad g\in C_0(X).
\end{equation}

\begin{proposition} Let $(X,D,K)$ be a CAS and $T:=[0,\infty[$.
Each random walk $(Z_t)_{t\in T}$ on $X$ with start in any $x_0\in X$ associated with any 
convolution semigroup $(\mu_t)_{t\in T}$ on $(D,*)$ is a Feller process.

In particular, $(Z_t)_{t\in T}$ admits a modification such that all paths of this modification are RCLL, i.e.,
they are right continuous such that for all $t>0$ the left limits exist.
For  $T:=[0,\infty[$, we thus shall assume in addition that the random walk $(Z_t)_{t\in T}$ on $X$ has the
RCLL property.
\end{proposition}

\begin{proof} In order to check (\ref{Feller-semigr-cond}), we fix some $g\in C_0(X)$ and $\epsilon>0$.
By Lemma \ref{uniform-cont}, we find some open neighborhood $U\subset D$ of $e$ such that $|g(x)-g(y)|\le\epsilon$
holds for all $x,y\in X$ with $\pi(x,y)\in U$. Hence, for each $x\in X$ and $t\ge0$,
$$\int_U \int_X |g(x)-g(y)|\> K_h(x,dy)\> d\mu_t(h)\le \epsilon.$$
Thus,
\begin{align}
|T_{\mu_t}g(x)-g(x)|&= \Bigl|\int_D \int_X (g(x)-g(y))\> K_h(x,dy)\> d\mu_t(h)\Bigr|
\notag\\
&=\epsilon + 2\|g\|_\infty\> \mu_t(D\setminus U)\le 2\epsilon
\notag
\end{align}
whenever $t$ is small enough. This proves (\ref{Feller-semigr-cond}).
 The second statement is well-known; see e.g., Section 17 of \cite{Ka}.
\end{proof}

The random walks  $(Z_t)_{t\in T}$ on $X$ may be  studied by using known results for random walks
 $(Y_t)_{t\in T}$ on $(D,*)$. This follows from the  following result which
 seems obvious at a first glance and can be seen easily in the group cases $X=G/H$, $D=G//H$.

\begin{theorem}\label{proj-random-walk}
Let $(X,D,K)$ be a CAS and $x_0\in X$ fixed and consider the 
continuous, open, and closed map
$\psi:X\to D$, $\psi(x):=\pi(x_0,x)$. 

 Let  $(Z_t)_{t\in T}$ be a  random
 walk on $X$ with start in $x_0$ associated with some convolution semigroup $(\mu_t)_{t\in T}$ on $(D,*)$
as described above where we assume that all paths are RCLL for $T=[0,\infty[$.
 Then the process $(\psi(Z_t))_{t\in T}$ is a random
 walk on $(D,*)$ with start in $e$ associated with $(\mu_t)_{t\in T}$.
\end{theorem}

For the proof we first consider the case $T=\mathbb N_0$ and 
 check that the projected process $(\psi(Z_t))_{t\in T}$
is a  Markov process. 

 For this we fix $n\in \mathbb N_0$ as well as Borel sets $A_0,\ldots, A_n\in\cal B(D)$. We 
consider the sub-probability measures
$P_{A_0,\ldots, A_n}\in  M^{(1)}(D)$ and $\tilde P_{A_0,\ldots, A_n}\in  M^{(1)}(X)$ with
$$P_{A_0,\ldots, A_n}(A):= P(\psi(X_0)\in A_0, \ldots,\psi(X_n)\in A_n, \psi(X_{n+1})\in A) 
\quad\quad (A\in\cal B(D))$$
and
$$\tilde P_{A_0,\ldots, A_n}(B):= P(\psi(X_0)\in A_0, \ldots,\psi(X_n)\in A_n, X_{n+1}\in B) 
\quad\quad (B\in\cal B(X)).$$
Then clearly, $P_{A_0,\ldots, A_n}$ is the image measure of $\tilde P_{A_0,\ldots, A_n}$ under $\psi$.
We need the following reconstruction of $\tilde P_{A_0,\ldots, A_n}$ from $P_{A_0,\ldots, A_n}$:

\begin{lemma}\label{eq-proj-random}
For all  $n\in \mathbb N_0$,  $A_0,\ldots, A_n\in\cal B(D)$, and $B\in\cal B(X)$,
$$\tilde P_{A_0,\ldots, A_n}(B)=\int_D K_h(x_0,B)\> dP_{A_0,\ldots, A_n}(h).$$
\end{lemma}

\begin{proof}[Proof of Lemma \ref{eq-proj-random}]
We prove the lemma by induction on $n$.

In fact, for $n=0$ we have the two cases $\psi(x_0)\not\in A_0$ and  $\psi(x_0)\in A_0$.
In the first case the assertion is trivial, and in the second case we get the claim from
$$\tilde P_{A_0}(B)=P(Z_1\in B)=K_{\mu_1}^X(x_0,B)=\int_D K_h(x_0,B)\>d\mu_1(h)=\int_D K_h(x_0,B)\>dP_{A_0}(h).$$

We now turn to the step $n-1\to n$ for $n\ge1$.
As all spaces are second countable and locally compact, we may use the concept of regular conditonal
probabilities, and obtain from the Markov property of  $(Z_t)_{t\in T}$, the assumption of our induction step, 
and from the axioms of a CAS that
\begin{align}\label{distribution-eq-1}
\tilde P&_{A_0,\ldots, A_n}(B)= P(\psi(X_0)\in A_0, \ldots,\psi(X_n)\in A_n, X_{n+1}\in B)
\\ &=
\int_{\psi^{-1}(A_n)} P(X_{n+1}\in B|\> \psi(X_0)\in A_0, \ldots,\psi(X_{n-1})\in A_{n-1}, X_n=x_n) 
\> d\tilde P_{A_0,\ldots, A_{n-1}}(x_n)
\notag\\ &=
\int_{\psi^{-1}(A_n)} K^X_{\mu_1}(x_n,B) \> d\tilde P_{A_0,\ldots, A_{n-1}}(x_n)
\notag\\ &=
\int_{A_n}\int_X  K^X_{\mu_1}(x_n,B) \> K_h(x_0,dx_n)\> dP_{A_0,\ldots, A_{n-1}}(h)
\notag\\ &=
\int_{A_n} K_h\circ K^X_{\mu_1}(x_0,B)\> dP_{A_0,\ldots, A_{n-1}}(h)
\notag\\ &=
\int_{A_n} K_h\circ K^X_{\mu_1}(x_0,B)\> dP_{A_0,\ldots, A_{n-1}}(h)
\notag\\  &=
\int_{A_n}\int_D K_l(x_0,B)\> d(\delta_h*\mu_1)(l) \> dP_{A_0,\ldots, A_{n-1}}(h)
\notag\\  &=
\int_D K_l(x_0,B)\>  d(P_{A_0,\ldots, A_{n-1}}|_{A_n} *\mu_1)(l)
\notag
\end{align}
where $|_{A_n}$  means the restriction of a measure to $A_n$.
If we take $B=\psi^{-1}(A)$ for $A\in\cal B(D)$ in (\ref{distribution-eq-1}),  we obtain
from the axioms of a CAS that
 \begin{align}\label{distribution-eq-2}
 P_{A_0,\ldots, A_n}&(A)=  P(\psi(X_0)\in A_0, \ldots,\psi(X_n)\in A_n, ,\psi(X_{n+1})\in A)
\\ &=    \int_D K_l(x_0,\psi^{-1}(A))\>  d(P_{A_0,\ldots, A_{n-1}}|_{A_n} *\mu_1)(l)
\notag\\  &=(P_{A_0,\ldots, A_{n-1}}|_{A_n} *\mu_1)(A).
\notag
\end{align}
If we insert this identity in the end of (\ref{distribution-eq-1}), we get the claim
$$\tilde P_{A_0,\ldots, A_n}(B)=\int_D K_h(x_0,B)\> dP_{A_0,\ldots, A_n}(h).$$
\end{proof} 

\begin{proof}[Proof of Theorem \ref{proj-random-walk}]
We first consider the case $T=\mathbb N_0$  and
 check that  $(\psi(Z_t))_{t\in T}$
is a  Markov process. 
For this consider $n\in\mathbb N$ and   $A_0,\ldots, A_n,A\in\cal B(D)$. Let $(\cal F_t)_{t\in T}$
be the canonical filtration of  $(\psi(Z_t))_{t\in T}$ on the associated probability space $(\Omega,\cal A,P)$.
Then, by Lemma \ref{eq-proj-random} and the first lines of (\ref{distribution-eq-1}),
\begin{align}\label{distribution-eq-3}
&\int_{\{\psi(X_0)\in A_0, \ldots,\psi(X_n)\in A_n\}} {\bf 1}_{\{\psi(X_{n+1})\in A\}} \> dP= \\
&\quad\quad=
 P(\psi(X_0)\in A_0, \ldots,\psi(X_n)\in A_n, \psi(X_{n+1})\in A)
\notag \\ &\quad\quad=
\int_{A_n} K_h\circ K^X_{\mu_1}(x_0,\psi^{-1}(A))\> dP_{A_0,\ldots, A_{n-1}}(h)
\notag \\ &\quad\quad=
\int_{\{\psi(X_0)\in A_0, \ldots,\psi(X_n)\in A_n\}} K_{\psi(X_n)}\circ K^X_{\mu_1}(x_0,\psi^{-1}(A))\> dP.
\notag 
\end{align}
As the last integrand is  measurable w.r.t.~the $\sigma$-algebra $\sigma(\psi(X_n))\subset \cal F_n$, we obtain
 from (\ref{distribution-eq-3}) that a.s.
$$ P( \psi(X_{n+1})\in A|\> \cal F_n)=K_{\psi(X_n)}\circ K^X_{\mu_1}(x_0,\psi^{-1}(A))$$
and thus a.s.
\begin{equation}\label{distribution-eq-4}
 P( \psi(X_{n+1})\in A|\>\sigma(\psi(X_n)))=
K_{\psi(X_n)}\circ K^X_{\mu_1}(x_0,\psi^{-1}(A))= P( \psi(X_{n+1})\in A|\> \cal F_n).
\end{equation}
Therefore,   $(\psi(Z_t))_{t\in T}$ is a Markov process. Moreover,   a comparison of (\ref{distribution-eq-4}) with
the transition probabilities of  random walks on $(D,*)$ 
shows immediately that  $(\psi(Z_t))_{t\in T}$ is a random walk  on $(D,*)$  associated with $(\mu_t)_{t\in T}$ and start in $e$
as  claimed.

Let us now turn to the case $T=[0,\infty[$. A  change of the notations of the preceding 
proof shows readily that for all $n\in\mathbb N$, $0< t_1<t_2<\ldots<t_{n+1}$, and $ A\in\cal B(D)$,
\begin{align}\label{distribution-eq-4-cont}
 P( \psi(X_{t_{n+1}})\in A|\>\sigma(\psi(X_{t_n})))&=
K_{\psi(X_{t_n})}\circ K^X_{\mu_{t_{n+1}-t_n}}(x_0,\psi^{-1}(A))
\notag\\
&= P( \psi(X_{t_{n+1}})\in A|\>\sigma(\psi(X_{t_0}),\ldots
\psi(X_{t_n}) ))  \quad\quad\text{a.s..}
\notag
\end{align}
Standard arguments from the theory of Markov processes with RCLL paths (see Section 17 of \cite{Ka}) now show that 
\begin{equation}\label{distribution-eq-5-cont}
 P( \psi(X_{t_{n+1}})\in A|\>\sigma(\psi(X_{t_n})))= P( \psi(X_{t_{n+1}})\in A|\>\sigma(\psi(X_t); t\in[0,t_n])) 
\quad\quad\text{a.s.}.
\end{equation}
This is the Markov property, and the proof can be completed  as for $T=\mathbb N_0$.
\end{proof}

For many classes of commutative hypergroups $(D,*)$ there exist limit theorems for 
 random walks  on $(D,*)$ like (strong) laws of large numbers, central limit theorems, and so on;
see Ch.~7 of  \cite{BH} and references there.
 Theorem \ref{proj-random-walk} may be now used to transfer these results to random walks on $X$ for 
suitable commutative CAS $(X,D,K)$. This seems to be interesting 
 in particular for examples which appear as deformations of 
group CAS, as here random walks  on $X$ may be seen as ``radial random walks with aditional drift'' on the
homogeneous space $X$. This seems to be interesting in particular for such random walks on
affine buildings and on noncompact Grassmann manifolds.
We shall study these examples in a forthcoming paper.

\section{Open problems}

We finally collect some open problems which appeared in the preceding course on CAS:
\begin{enumerate}
\item[\rm{(1)}] Do there exist (commutative) CAS with (T1), but without (T2)?
\item[\rm{(2)}] Does (T2) always imply (T1)? 
Notice that this is the case in the discrete and in the commutative case, 
and that it is likely that it can be shown in the compact case.
\item[\rm{(3)}] Is each discrete CAS a generalized association scheme? 
The answer is positive in the finite case. The general problem is part of: 
\item[\rm{(4)}] Does each (commutative) CAS $(X,D,K)$ admit a further associated (commutative)
 CAS $(X,D,\tilde K)$ with the same spaces $X,D$, and the same projection $\pi$ such that  $(X,D,\tilde K)$
has property (T2).
\item[\rm{(5)}] Give examples of  (commutative) CAS $(X,D,K)$ with $X$,$D$ connected, which are not 
of the form $X=G/H$, $D=G//H$ for a locally compact group $G$ with a compact subgroup $H$ or deformations 
of such group cases.
\end{enumerate}

\bibliographystyle{amsplain}

\end{document}